\documentclass[reqno]{amsproc}

\usepackage{amstext,amsmath,amsthm,amssymb,amsfonts}
\usepackage[latin1]{inputenc}
\usepackage{epsfig}
\usepackage{hyperref}
\usepackage{color}

\usepackage{inputenc}
\usepackage{mathrsfs}

\usepackage{latexsym}

\usepackage{enumitem}
\usepackage{bbm}
\usepackage{geometry}
\usepackage{amsthm}

\definecolor{myurlcolor}{rgb}{0,0,0.7}

\newtheorem{theorem}{Theorem}[section]
\newtheorem{lemma}[theorem]{Lemma}
\newtheorem{corollary}[theorem]{Corollary}

\newtheorem{definition}[theorem]{Definition}

\newtheorem{proposition}[theorem]{Proposition}
\newtheorem{remark}[theorem]{Remark}

\newcommand{\vspi}{\vspace{0.4cm}}

\newcommand{\bee}{\mathbf{b}}

\newcommand{\bea}{\begin{eqnarray}}
\newcommand{\eea}{\end{eqnarray}}
\newcommand{\beq}{\begin{equation}}
\newcommand{\eeq}{\end{equation}}
\newcommand{\enn}{\nonumber \end{equation}}

\newcommand{\col}{\text{color}}

 \newcommand{\Z}{\mathbb{Z}}

 \newcommand{\cG}{\mathcal{G}}
 \newcommand{\cV}{\mathcal{V}}
 \newcommand{\cE}{\mathcal{E}}
 \newcommand{\cF}{\mathcal{F}}
 \newcommand{\cH}{\mathcal{H}}

 \newcommand{\cB}{\mathcal{B}}
 \newcommand{\cT}{\mathcal{T}}
 \newcommand{\cR}{\mathcal{R}}

\newcommand{\mT}{\mathfrak{T}}

  \newcommand{\crrG}{ \mathscr{G} }
  \newcommand{\crrtG}{ \mathscr{G}_{ \mathfrak{f}^0} }
  \newcommand{\cGcol}{ \mathscr{G}_{  \col; \mathfrak{f}^0} }

 \newcommand{\brrG}{ {\partial \mathscr{G}_{ \mathfrak{f}^0 } } }

\newcommand{\bG}{\partial{\mathcal{G}}_{\mf^0}}
\newcommand{\bV}{{\mathcal{V}}_{\partial}}
\newcommand{\bE}{{\mathcal{E}}_{\partial}}
\newcommand{\bC}{{\mathcal{C}}_{\partial}}

\newcommand{\wt}{\widetilde}

 \newcommand{\mf}{\mathfrak{f}}

 \newcommand{\sset}{\Subset}

 \newcommand{\inter}{{\rm int}}
 \newcommand{\ext}{{\rm ext}}

\title[Tutte and BR polynomial for stranded graphs]{
Extending the Tutte and Bollob\'as-Riordan Polynomials 
\vspace{0.1cm}\\ to Rank 3 Weakly-Colored Stranded Graphs}

\author{Remi C. Avohou}
\address[R.C.A.]{
University of Abomey-Calavi, 
International Chair in Mathematical Physics and Applications,
ICMPA-UNESCO Chair, 072BP50, Cotonou, Rep. of Benin}
\email{remi.avohou@cipma.uac.bj}

\author{Joseph Ben Geloun}
\address[J.B.G.]{Universit\'e Paris 13, Sorbonne Paris Cit\'e, 99, J.-B. Cl\'ement
LIPN, Institut Galil\'ee, CNRS UMR 7030, 93430, Villetaneuse, France, and
University of Abomey-Calavi, 
International Chair in Mathematical Physics and Applications,
ICMPA-UNESCO Chair, 072BP50, Cotonou, Rep. of Benin}
\email{bengeloun@lipn.univ-paris13.fr}

\author{Mahouton N. Hounkonnou}
\address[M.N.H.]{
University of Abomey-Calavi, 
International Chair in Mathematical Physics and Applications,
ICMPA-UNESCO Chair, 072BP50, Cotonou, Rep. of Benin}
\email{norbert.hounkonnou@cipma.uac.bj}

\begin{document}

\maketitle 

\begin{abstract}

The Bollob\'as-Riordan polynomial [Math. Ann. 323, 81 (2002)]  is a universal polynomial invariant for ribbon graphs. 
We find an extension of this  polynomial for a particular family of  combinatorial objects, called rank 3 weakly-colored stranded graphs.   Stranded graphs arise in the study of tensor models for quantum gravity in physics, and generalize graphs and ribbon graphs.  We present a  seven-variable polynomial invariant of these graphs, which obeys a contraction/deletion recursion relation similar to that of the Tutte and Bollob\'as-Riordan polynomials. 
 However, it is defined on a much broader class of objects, 
 and furthermore captures properties that are not encoded by the Tutte  
or Bollob\'as-Riordan polynomials. 

\vspi 

MSC(2010): 05C10, 57M15

Keywords: Graph Theory, Tutte polynomial, Bollob\'as-Riordan polynomial, Tensor Models.  

 \end{abstract}

\tableofcontents

\section{Background and main results}

 The Tutte polynomial is a universal polynomial invariant defined on  abstract graphs. This invariant obeys a particular recursion relation for the contraction and deletion of the  edges of a graph, see for instance \cite{tutte}. 
The Bollob\'as-Riordan (BR) polynomial, see \cite{bollo2} and \cite{bollo},  defines a  genuine extension of the Tutte polynomial  to ribbon graphs or neighborhoods of graphs embedded in surfaces. 
Before, in the context of quantum groups, 
Reshetikhin and Turaev  showed in  \cite{Reshetikhin:1990pr} the existence 
of a generalized Jones polynomial on graphs embedded in surfaces. 

In this work, we investigate  families of combinatorial objects generalizing abstract and ribbon graphs 
to which both the Tutte and BR polynomial invariants may be extended. 
 These generalized graphs were identified as Feynman graphs of matrix 
\cite{Di Francesco:1993nw} 
 and tensor models \cite{ambjorn}.  Such models in physics aim at defining a quantum spacetime, and their graphs represent discrete geometries.
The study of colored tensor graphs, i.e.  
the Feynman graphs of colored tensor models introduced in  \cite{Gurau:2009tw}
(see also \cite{Gurau:2009tz}, \cite{Gurau:2011xp} and  \cite{Geloun:2009pe}), 
has  acknowledged  a growing interest after the recent discovery, by Gurau in \cite{Gurau:2010ba}, of the main organizing tool of their partition function. 
  The story of colored tensor graphs has also been successful 
from the point of view of combinatorics and graph theory 
because of their tractable properties. For instance,  well-known combinatorial notions 
 defined on graphs such as polynomial invariants and Hopf-algebras 
 have been extended to tensor graphs in \cite{Gurau:2009tz}, \cite{Tanasa:2010me}, \cite{Avohou:2015sia}, and \cite{Raasakka:2013kaa}.

There are two main works\footnote{Note that, while this work was under submission, 
the present work have been generalized to arbitrary rank \cite{Avohou:2015}.}  addressing 
 a generalization of the BR polynomial  to a  family of generalized graphs similar
to that  presented here:
the works by Gurau in \cite{Gurau:2009tz} and 
	by Tanasa in \cite{Tanasa:2010me}. 
For tensor graphs, the vertex is stranded, and the usual way of contracting an edge immediately leads  to a different vertex. As a consequence, the simplest prescriptions of contracting and deleting edges are not well defined without further considerations.  Indeed, the two aforementioned 
works undertake the definition of the contraction and deletion of edges in an  unusual manner. 

There are several layers of difficulties   while defining  the type of graphs that ought to be considered. One should also distinguish
the type of topological polynomial defined for the family of graphs,
and the type of recursion relations that the invariant polynomial needs to satisfy. 
Our method is radically different from that of Tanasa, 
since we  use the Gurau colored prescription to keep under control  the type of generated graphs. 
Moreover, and in contrast with the work by Gurau, 
we propose to enlarge the family of  colored tensor
graphs to what is called   {\it weakly-colored (w-colored) stranded graphs},
 for which contraction, and a similar notion of deletion
 make sense. These w-colored stranded graphs 
are precisely the result of edge contractions of colored tensor graphs. 
 This class of graphs is closed under contraction.
From this stability,   we can define, 
for w-colored graphs, a polynomial which is invariant  under a linear recursion relation. 
 We must mention that there might be  room for improvement for the present definition of w-colored stranded graphs.
Although this is not needed in our analysis, there could be an equivalent description in terms of its basic constituents. 
Such a study is left for the future.

Our main result appears in  Theorem \ref{theo:contens} of Section \ref{sect:tgraph}. 
This statement  determines the contraction and cut (an operation  replacing the deletion)  recursion relation with respect to an edge, fulfilled by a new polynomial given in Definition \ref{def:topotens}. The  class of graphs on which the polynomial is defined has  several components. Some of these components turn out to be well defined for  abstract  and ribbon graphs, 
and therefore, they are introduced step by step in the sequel.   

A first component is captured by the notion of half-edge
 which can be introduced at the  abstract graph level. 
 Half-edges are of crucial importance in Physics, see \cite{riv} and \cite{krf},  for instance.
The analysis of half-edged graphs (HEGs) is the purpose of Section \ref{sect:gratu}.
The Tutte polynomial for   half-edged graphs satisfies 
   linear recursion identities  similar to  those of the Tutte polynomial for abstract
graphs (see Proposition \ref{theo:contutefla}),
and so it turns out to be an evaluation of it. 

A second key  notion  we  need is that of open and closed
faces of    a ribbon graph with ribbon half-edges. The consequences that open and closed faces in a ribbon graph have on the related  BR polynomial, are addressed in Section \ref{sect:ribbonbollo}.  Another important result of this work is set in  Theorem \ref{theo:BRext}, which establishes  a new recursion relation
obeyed by the BR polynomial for ribbon graphs with ribbon half-edges. 
 Then follow the concluding remarks in Section 5.

 Finally, a closing appendix provides  explicit examples of the polynomial obtained for some w-colored graphs.

\section{Half-edged graphs and the Tutte polynomial}
\label{sect:gratu}

 We start by setting up our graph theoretic notation.

\subsection{Basic definitions and  notation}

 We  use $G(\cV,\cE)$, or at times just $G$, to denote a graph (which allows loops and multiple edges)
with vertex set $\cV$ and edge set $\cE$.

 We denote a spanning subgraph $A(\cV,\cE_A)$ of a graph $G(\cV,\cE)$ by  $A \subset G$. 
Consider two disjoint graphs $G_1(\cV_1,\cE_1)$ and
$G_2(\cV_2,\cE_2)$.  
 The disjoint union of the graphs $G_1$ and $G_2$  is denoted by 
$G_1 \sqcup G_2$.  
 The one-point join graph $G_1 \cdot _{v_1,v_2} G_2$ 
  simply merges  the two graphs at the vertices $v_1 \in \cV_1$ and $v_2 \in \cV_2$, by removing $v_1$ and $v_2$, and 
inserting a new vertex $v$ which has all edges incident 
to $v_1$ and to $v_2$.

 Using the notation of  \cite{bollo3}, we have the following definition. 
\begin{definition}[Tutte polynomial]
\label{defTutte1}
Let $G(\cV,\cE)$ be a graph, then the \emph{Tutte polynomial} $T_{G}$ 
of $G$ is  defined as
\beq
    T_{G}(x,y)=\sum_{A\subset G}(x-1)^{r(G)-r (A)}(y-1)^{n(A)},
\label{tutexp}
\eeq
    where $k(A)$ is the number of connected components of
 the spanning subgraph $A$,
$r(A)=|\cV|-k(A)$ its \emph{rank}, 
and $n(A)=|\cE_A|+k(A)-|\cV|$ is its \emph{nullity} or \emph{cyclomatic number}.
\end{definition}

There is an equivalent definition of this polynomial  (in  notation of \cite{bollo3}): 

\begin{definition}[Tutte polynomial 2]
\label{defTutte2}
The Tutte polynomial $T_{G}$ of $G$ is defined as:

(a) If $G$ has no edges, then $T_{G}(x,y)=1.$

 Otherwise, for any edge $e\in \cE$,
    
(b) $T_{G}(x,y)=T_{G-e}(x,y) + T_{G/e}(x,y)$, if $e$ is an ordinary edge  (neither a bridge nor a loop).

(c) $T_{G}(x,y)=xT_{G-e}(x,y) = x T_{G/e}$, if $e$ is a bridge.

(d) $T_{G}(x,y)=yT_{G-e}(x,y)=y T_{G/e}(x,y)$, if $e$ is a loop.
\end{definition}

For a terminal form composed with $m$ bridges and $p$ loops, the Tutte polynomial evaluates  as $x^m y^p$. Let $G_1$ and $G_2$ be two disjoint graphs, then 
we have 
\bea 
T_{G_1 \sqcup G_2}= T_{G_1} T_{G_2} = T_{G_1 \cdot_{v_1,v_2} G_2}\,,
\label{cups0}
\eea
for any  vertices  $v_{1}\in G_{1}$ and $v_2 \in G_2$
 (for proofs, see \cite{bollo3}). 
These relations  hold  because of  the additivity property of the 
rank and the nullity with respect to the disjoint union and
 one-point join operation for graphs.  Note also that the 
same identities can be also easily derived from the contraction/deletion relations.

\subsection{Half-edged graphs (HEGs)}
\label{subsec:tufla}

 We now introduce the  family of graphs with half-edges, called henceforth HEGs (for half-edged graphs).

\begin{figure}[h]
 \centering
     \begin{minipage}[t]{.8\textwidth}
      \centering
\includegraphics[angle=0, width=1cm, height=0.1cm]{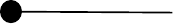}
\caption{ {\small A  half-edge incident to a unique vertex. }}
\label{fig:flag}
\end{minipage}
\end{figure}

A  half-edge  is represented by a line incident to a unique vertex and without forming a loop. 
(See Figure \ref{fig:flag}.) 

\begin{figure}[h]
 \centering
     \begin{minipage}[t]{.8\textwidth}
      \centering
\includegraphics[angle=0, width=2.5cm, height=1.5cm]{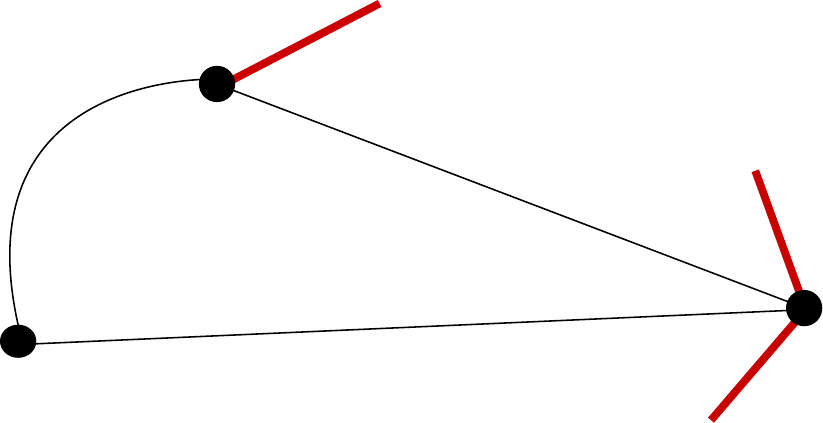}
\caption{ {\small A HEG and its  set $\mf^0$ of 
half-edges in red.}}
\label{fig:heg}
\end{minipage}
\end{figure}

\begin{definition}[Half-edged graph (HEG)]\label{heg}
A \emph{ HEG } $G(\cV,\cE,\mf^0)$ or more briefly $G_{\mf^0}$
is a graph $G(\cV,\cE)$ with a set $\mf^0$, whose elements are called \emph{half-edges} together with a relation which associates  with each  half-edge a unique vertex. 
(See Figure \ref{fig:heg}.)  The graph $G$ is called the \emph{underlying graph} of $G_{\mf^0}$.
\end{definition}
From the above definition, it is direct to observe that an abstract graph is a HEG with $\mf^0=\emptyset$.

\begin{figure}[h]
 \centering
     \begin{minipage}[t]{.8\textwidth}
      \centering
\includegraphics[angle=0, width=5cm, height=0.1cm]{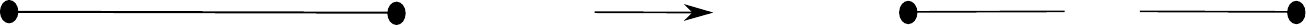}
\caption{ {\small Cutting an edge. }}
\label{fig:cute}
\end{minipage}
\end{figure}

\begin{figure}[h]
 \centering
     \begin{minipage}[t]{.8\textwidth}
      \centering
\includegraphics[angle=0, width=6cm, height=1.5cm]{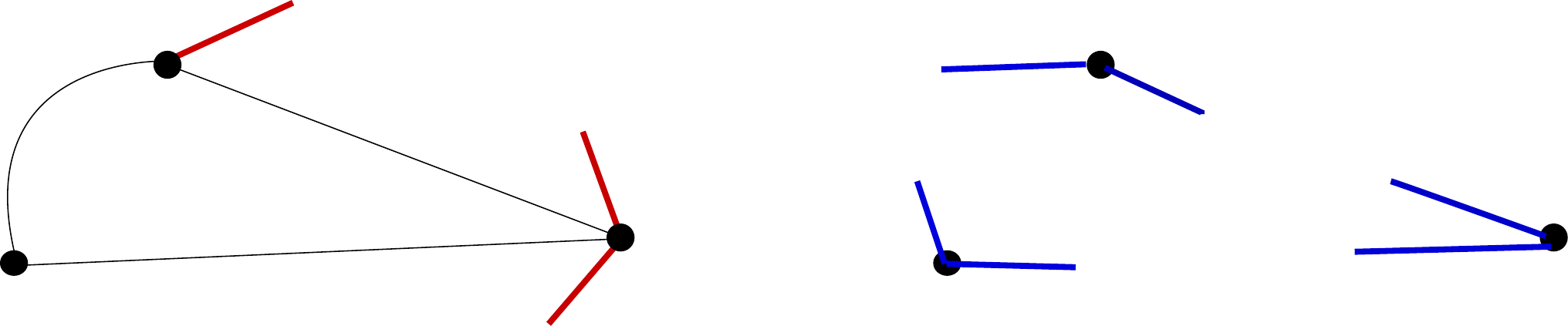}
\caption{ {\small A HEG (on the left), its set $\mf^0$ of 
half-edges in red, and set $\mf^1$ in blue  obtained by cutting} all its edges and removing
$\mf^0$. }
\label{fig:grafla}
\end{minipage}
\end{figure}

\begin{definition}[Cut of an edge \cite{riv}] \label{def:cut}
Let $e \in \cE$ be an edge of $G(\cV,\cE,\mf^0)$. The \emph{cut HEG}
$G_{\mf^0}\vee e$ is obtained from $G_{\mf^0}$ by \emph{cutting} $e$, 
namely by deleting $e$ and  attaching a half-edge to each of the end vertices 
of $e$.  Two half-edges are attached to the same vertex if $e$ is a loop. (See Figure \ref{fig:cute}.) 
We denote by $\mf^1$ the set of half-edges that result from cutting
all edges in $\cE$, and write $\mf$ for $\mf^0 \dot{\cup}\mf^1$. 
 (See an illustration in Figure \ref{fig:grafla}.) 
\end{definition}

\begin{figure}[h]
 \centering
     \begin{minipage}[t]{.8\textwidth}
      \centering
\includegraphics[angle=0, width=10cm, height=1.5cm]{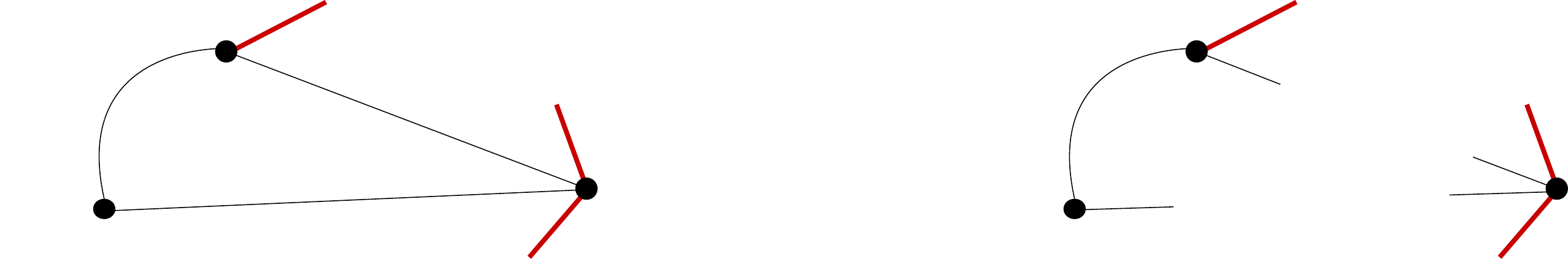}
\vspace{0.1cm}
\caption{ {\small{A HEG $G_{\mf^0}$ and the spanning c-subgraph $A$ associated with $e$}. }}
\label{fig:cgra}
\end{minipage}
\put(-270,-10){$G_{\mf^0}$}
\put(-306,17){$e$}
\put(-83,-10){$A$}
\put(-130,17){$e$}
\end{figure}

\begin{definition}[Spanning c-subgraphs of a HEG]\label{subfla}

  A spanning c-subgraph  $A$ of $G(\cV,\cE,\mf^0)$ is the result of taking a spanning subgraph 
of $G(\cV,\cE)$ viewing it as embedded in $G(\cV,\cE,\mf^0)$, 
then cutting all the edges of $\cE \setminus\cE_A$. More precisely, a \emph{spanning c-subgraph}  $A$ of $G(\cV,\cE,\mf^0)$ is defined as  a HEG $A(\cV_A,\cE_A,\mf^0_A)$, the edge set $\cE_A$ of which 
is a subset of $\cE$ with all vertices  and all additional half-edges of $G_{\mf^0}$.  Hence $\cE_A\subseteq \cE$, $\cV_A = \cV$, and $\mf^0_A = \mf^{0} \cup \mf^{1}_A(\cE_A)$, where $\mf^{1}_A (\cE_A)$ is the set of half-edges  
obtained by cutting all edges in $\cE\setminus\cE_A$. We denote it $A \sset G_{\mf^0}$. 
 (See $A$ as an illustration in Figure \ref{fig:cgra}.) 

\end{definition}

The isomorphism class of HEGs follows the same idea of that of abstract graphs. The sets of half-edges in the
same class of HEGs must be  of the same cardinality and satisfy the same incidence relation onto vertices.

 Note that a spanning c-subgraph  $A \sset G_{\mf^0}$, unless $A =G_{\mf^0}$, has always a 
greater number of half-edges than $G_{\mf^0}$. The rank and nullity of a HEG $G_{\mf^0}$ are defined
to be the rank and nullity of the underlying graph $G$: $r(G_{\mf^0})= r(G)$, $n(G_{\mf^0})=n(G)$. 
 Deleting an edge $e$ in a HEG $G_{\mf^0}$ is deleting it from its underlying graph
while keeping all half-edges and their incidence relation. We  denote it in the standard way $G_{\mf^0}-e$.

\begin{definition}[Edge contraction of HEG]
\label{contfla}
Let $G(\cV,\cE,\mf^0)$ be a HEG. 
We define the \emph{contraction} of a non-loop 
edge $e$ in $G_{\mf^0}$, i.e. $G_{\mf^0}/e$  to be the HEG obtained
from $G_{\mf^0}$ by removing $e$ and identifying 
the two end vertices into a new vertex having all their additional 
half-edges
and remaining incident lines.  For a loop $e$, 
contraction $G_{\mf^0}/e$ and deletion $G_{\mf^0}-e$ coincide. 
\end{definition}

Noting that the edge contraction does not 
change the number of additional half-edges in a HEG, the following proposition is straightforward.
\begin{proposition}\label{cutflag}
 Let $G_{\mf^0}=G(\cV,\cE, \mf^0)$ be a HEG and $e$ one of its edge. Then
$G_{\mf^0}/e=(G_{\mf^0}/e)(\cV', $ $\cE\setminus \{e\}, \mf^0)$, 
where $\cV'$ is the set of vertices of the underlying abstract graph $G/e$,
and  $G_{\mf^0}\vee e=(G_{\mf^0}\vee e)(\cV,\cE\setminus \{e\}, \mf^0 \cup \{e_1,e_2\})$,
where $e_{1}$ and $e_2$ are the two half-edges obtained by cutting $e$. 
\end{proposition}

The disjoint union  
 and one-point join of HEGs are defined as for graphs carrying the half-edges along. 
Spanning c-subgraphs of ${G_1}_{\mf_1^0} \sqcup {G_2}_{\mf_2^0}$ and
of ${G_1}_{\mf_1^0}  \cdot_{v_1,v_2} {G_2}_{\mf_2^0}$ are of the form 
$A_{1} \sqcup A_{2}$ and $A_{1}\cdot_{v_1,v_2}  A_{2}$, 
respectively, with $A_i \sset {G_i}_{\mf^0_i}$, $i=1,2$.  The spanning c-subgraphs
of ${G_1}_{\mf_1^0}\sqcup {G_2}_{\mf_2^0}$ and of ${G_1}_{\mf_1^0} \cdot_{v_1,v_2} {G_2}_{\mf_2^0}$   are thus in  one-to-one correspondence.

\subsection{A Tutte polynomial for HEGs}
\label{subsec:tutpoly}

We extend  the Tutte polynomial to HEGs and determine some properties of the new function. 
\begin{definition}
 Let  $G_{\mf^0}$ be a HEG and $G$ its underlying graph. Then 
\beq\label{ctt}
 \cT_{G_{\mf^0}}(x,y,t)=t^{|\mf^0|+2n(G_{\mf^0})} \,T_{G}(X, \frac{Y}{t^2}) \,,
\eeq
where $X=t^2(x-1) + 1$ and $Y=y +t^2-1$.
\end{definition}
The particular expressions for $X$ and $Y$ were determined to get
a simplified form of the following properties. 
\begin{proposition}\label{theo:contutefla}
The function $\cT$ has the following properties for a HEG 
$H= G_{\mf^0}$, analogous to those of the Tutte polynomial for 
graphs. 

(a) $\cT_{(E_n)_{\mf^0}}(x,y,t) = t^{|\mf^0|}$, where $(E_n)_{\mf^0}$ is
the HEG with $n$ vertices, no edges and $|\mf^0|$ half-edges. 

(b)  
\bea
\cT_{H}= \left\{
\begin{array} {lr}
 X\cT_{H-e}(x,y,t) &  \text{if $e$ is a bridge,}\\
 Y \cT_{H/e} (x,y,t) &  \text{if $e$ is a loop,} \\ 
 t^2 \cT_{H-e} (x,y,t) + \cT_{H/e}(x,y,t) & \text{if $e$ is ordinary;}
\end{array}
\right.  
\eea

(c) 
\beq
\cT_{H} =\cT_{H \vee e} (x,y,t) + \cT_{H/e}(x,y,t) \quad \text{if $e$ is ordinary;}
\eeq

(d) 
\beq
\cT_{H} (x,y,t) = 
\sum_{A\sset H}(x-1)^{r(H)-r(A)}(y-1)^{n(A)} \,t^{|\mf^0(A)|}\,.
\label{tuflaexp}
\eeq
\end{proposition}

\proof
Properties $(a)$ and $(b)$ follow immediately from Definition \ref{defTutte2},
using the usual additive properties of the rank for bridges, loops and ordinary
edges, and noting that $|\mf^0(H)|=|\mf^0(H-e)|=|\mf^0(H/e)|$. 

For property $(c)$, note that $H \vee e$ is just $H-e$ with the addition 
of two half-edges appended to the end points of $e$,  so that
$\cT_{H\vee e}(x,y,t)= t^2 \cT_{H-e}(x,y,t)$. 

To prove property $(d)$, we use Definition \ref{defTutte1}, 
substituting the expressions for $X$ and $Y$ given 
in \eqref{ctt}, and noting that there is a one-to-one correspondence
between the spanning subgraphs of $H=G_{\mf^0}$ and $G$,
and that  $|\mf^0(A)|=|\mf^0|+2|\cE|-2|\cE_A|$. 

\qed

\vspi 

 Note that property $(c)$ can replace the third term of property $(b)$
suggesting that cut may be a more natural operation on HEGs than
deletion.  The polynomial $\cT_{G_{\mf^0}}$ has always a factor of $t^{|\mf^{0}|}$
which may be removed by a new normalization. 
 The exponent of $t$ in $\cT_{G_{\mf^0}}(x,y,t)/t^{|\mf^0|}$ is always even since each c-spanning subgraph 
is defined via successive cut of edges yielding each two half-edges. 
 The polynomial of the terminal form of a HEG with $m$ bridges, $n$ loops
and $q$ additional half-edges is
\beq
 ( 1+ (x-1)t^2)^m\, (y-1 + t^2)^n \,t^q=X^m Y^n t^q\,.
\eeq
For $G_{1 \mf^0_1}$ and $G_{2 \mf^0_2}$  two HEGs, then 
it is immediate to get 
\bea 
\cT_{G_{1 \mf^0_1} \sqcup G_{2 \mf^0_2}}= \cT_{G_{1 \mf^0_1}} \cT_{G_{2 \mf^0_2}} = \cT_{G_{1 \mf^0_1} \cdot_{v_1,v_2} G_{2 \mf^0_2}}\,,
\label{optufla}
\eea
for any vertices $v_{i}$ in $G_{i \mf^0_i}$, $i=1,2$, 
from the  properties \eqref{cups0} and \eqref{ctt}, and 
because the number
of half-edges and nullity are additive on disconnected graphs. 
 Putting $t=1$ in $\cT_{G_{\mf^0}}(x,y,t)$ yields the ordinary 
Tutte polynomial 
\beq
\cT_{G_{\mf^0}}(x,y,1)=T_{G}(x,y)\,. 
\eeq

\section{Ribbon graphs and  the Bollob\'as-Riordan polynomial}
\label{sect:ribbonbollo}

In this section, we first recall the generalization of the Tutte polynomial to ribbon graphs known as the BR 
polynomial for ribbon graphs,  following the notation of \cite{bollo} and \cite{bollo2}. 
  Then, we investigate ribbon graphs with half-edges, or half-ribbons, analogous to the HEGs of Subsection \ref{subsec:tufla}. 

\subsection{Ribbon graphs}
\label{subsect:ribbon}

We recall basics of ribbon graphs and
 the BR polynomial according to conventions of \cite{bollo}.
 Note that the following axioms are equivalent to those  of \cite{Reshetikhin:1990pr}. 

\begin{definition}[Ribbon graphs]\label{def:ribbongraph}
	A \emph{ribbon graph} $\cG$ is  a (not necessarily orientable) surface with boundary represented as the union of two 
sets of closed topological discs called vertices $\cV$ and edges $\cE.$ These sets satisfy the following:

$\bullet$ vertices and edges intersect  in  disjoint line segments,

$\bullet$ each such line segment lies on the boundary of precisely one vertex and one edge,

$\bullet$ every edge contains exactly two such line segments.
\end{definition}

 Defining the class of ribbon graphs we are considering, we  follow conventions of \cite{bollo}, \cite{bollo2} and  \cite{ellis0}. 
The following description of ribbon graphs, known as rotation systems, 
is dated back to  \cite{Heffter}. 
A \emph{signed rotation system} is a graph $G$ together with a cyclic 
ordering of the edges at each vertex of $G$, and an assignment 
of a sign $+$ or $-$ to each edge of $G$. The \emph{flip} of a vertex $v$
is an operation on the signed rotation system which reverses
the cyclic order of the edges incident to $v$ and switches all signs
of edges incident to $v$,  except the signs of its  \emph{loops} (edges incident to the same vertex). Two  rotation systems  are  \emph{equivalent}  if 
they can be transformed by a sequence of vertex flips composed by 
graph isomorphisms. 

There is an equivalence between signed rotation systems and ribbon graphs. 
To make this equivalence manifest, consider a ribbon graph $\cG$,  
choose an orientation on each  vertex,
and assign to each edge an orientation according to
the fact that the orientation of its end vertices across the edge 
are consistent or not, respectively. The underlying graph $G$ of $\cG$
with the vertex and edge orientations is a signed rotation
system (the initial choice of orientations for the vertex and edges of $\cG$ does not matter in the construction of the signed rotation system because the latter is stable under vertex flips). Given now a rotation system $G$, we can construct a ribbon graph 
from $G$ by replacing its vertices by discs, giving each disc an orientation, 
and attaching ribbon edges to these discs in the order given by the cyclic order
of the vertices of $G$. A ribbon edge has an orientation: if this orientation
is $+$, the ribbon edge is attached with both the ends in a consistent way with the orientation of its end vertices, otherwise, if the orientation is $-$, 
then exactly one end is attached consistently with one end vertex. 
We call ribbon edges with a $+$ sign \emph{positive} and 
\emph{negative}, otherwise.

If the notions of   \emph{ordinary} ribbon edges and \emph{bridges} need  no comment,
that of  loops  in ribbon graphs needs  careful
attention.   A \emph{loop} is a ribbon edge incident to the same vertex.
 A loop $e$ at a vertex $v$ is called \emph{trivial} if there is no cycle in $\cG$ which can be contracted to form a loop $f$ at $v$ such that the ends of $e$ and $f$ alternate in the cyclic order at $v$  (see again \cite{bollo}). 

We can address now the notion of contraction and
deletion of ribbon edges.  The following operations on ribbon graphs  can be found detailed
and illustrated in  \cite{bollo}. 
Let $\cG$ be a ribbon graph and $e$ one of its edges. 
 We call $\cG-e$ the ribbon graph obtained from 
$\cG$ by \emph{deleting} $e$.  If $e$ is not a loop, consider its
end vertices $v_1$ and $v_2$. The graph $\cG/e$  obtained by \emph{contracting} $e$ is defined from $\cG$ by replacing $e$, $v_1$ and $v_2$ by a single vertex disc $e\cup v_1 \cup v_2$. We now suppose that $e$ is a loop and
is incident to a unique vertex $v$.
Two situations may occur: either $e\cup v$ forms a M${\rm \ddot{o}}$bius band (trivial negative loop) 
with a single boundary cycle,  or an annulus (trivial positive loop) with two boundary cycles. 
 If  $e\cup v$ forms a trivial negative loop, the graph $\cG/e$ is obtained from $\cG$ by deleting $e$ and $v$ and adding one new vertex  whose boundary is the boundary of $e \cup v$. If $e$ is a trivial positive loop, the contraction of $e$ is
the deletion of $e$ and $v$ and the addition of  two  vertices, the union of whose boundaries is the boundary of $e \cup v$. 
In this case, the operation at a vertex $v$ splits the cyclic order at $v$ into two cycles, which may represent two vertices into which $v$ is split.
We must emphasize that this operation of contraction should preserve the incidence and cyclic ordering of the remaining edges of the graph.

We  use the following terminology: 

\begin{definition}[Faces]
Consider $\cG$ a ribbon graph as a surface with boundary. 
A \emph{face} is a boundary component of $\cG$. 
\end{definition}
Ribbon graphs are known to be equivalent to graphs cellularly embedded in surfaces, see for instance \cite{ellis0}. 
 A face of a ribbon graph uniquely corresponds to a  face of the embedding.
In the next subsections, after introducing more combinatorial definitions,
the notion of face might refer to a combinatorial object different from the
usual connected component of the boundary of a ribbon graph. 

Spanning subgraphs in this context are denoted again as $A \sset \cG$. We are in position to define the BR polynomial. 

\begin{definition}[BR polynomial \cite{bollo}]
\label{defBRpoly}
Let $\cG$ be a ribbon graph. The  \emph{BR polynomial} of $\cG$ is an element of $\Z[x,y,z,w]$ quotiented by 
the ideal generated by $w^2-w$, given by:
\beq
R_{\cG}(x,y,z,w) =\sum_{A\sset \cG} (x-1)^{r(\cG)-r(A)}(y-1)^{n(A)}z^{k(A)-F(A)+n(A)} w^{o(A)},
\label{brpoly}
\eeq
where $F(A)$ is the number of faces of $A$, and where $o(A)=0$ if $A$ is orientable, and $o(A)=1$ if not.  
\end{definition}

 In \cite{bollo}, it is proved that the BR polynomial also obeys a contraction and deletion rule for ordinary edges as
\beq\label{cbr}
R_{\cG}=R_{\cG/e}+ R_{\cG-e}\,.
\eeq 
Also, for every bridge $e$ of $\cG$, we have 
\beq
R_{\cG}=x \, R_{\cG/e}\,,
\label{cbrid}
\eeq 
for a trivial positive loop, 
\beq
R_{\cG}=y \,R_{\cG-e}\,,
\label{csel}
\eeq
and for a trivial negative loop, the following relation holds
\beq
R_{\cG}=(1+(y-1)zw) \,R_{\cG-e}\,.
\label{ctsel}
\eeq
 The relations \eqref{cbrid}-\eqref{ctsel} 
are assimilated to boundary conditions  for the contraction deletion recursion relation.  
These boundary conditions were extended to a larger family of one-vertex ribbon graphs  in \cite{avo}. 
 Nevertheless, in contrast with the Tutte polynomial, the relations \eqref{cbr}-\eqref{ctsel}  and the like are not sufficient
for computing the BR polynomial for arbitrary ribbon graphs after a series of contractions and deletions.

  Take any positive integer valued function  $h: \mathbb{N}^3 \to \mathbb{N}$ and replace exponent of $z$ in Definition \ref{defBRpoly}
by \beq
h(k(A),F(A),n(A))\,.
\eeq
One can show that \eqref{cbr} still holds for an ordinary edge. 
The choice $k(A) - F(A) + n(A)$ for the exponent of $z$ is interesting because 
\beq
k(A) - F(A) + n(A)= 2k(A) - (|\cV| - |\cE_A| - F(A))
\eeq
 is nothing but the genus or twice the genus (for oriented surfaces) of the subgraph  $A$.  Furthermore, writing the exponent in this form also helps for the determination of the terminal forms because $k(A) - F(A)$ 
and $n(A)$ turn out to be additive quantities with respect to 
the  one-point join of disjoint graphs. 
These remarks will be exploited 
to find  a generalization of the BR polynomial to stranded graphs.

\subsection{Half-edged ribbon graphs (HERGs)}
\label{subset:ribfla}

We generalize the BR polynomial to ribbon graphs with half-edges or half-ribbons. 
For convenience, 
we use an idea by Chmutov in \cite{schmu} who considers a ribbon edge as a rectangle incident to the corresponding vertices along a pair of its opposite sides.

\begin{figure}[h]
 \centering
     \begin{minipage}[t]{.9\textwidth}
      \centering
\includegraphics[angle=0, width=3cm, height=1.3cm]{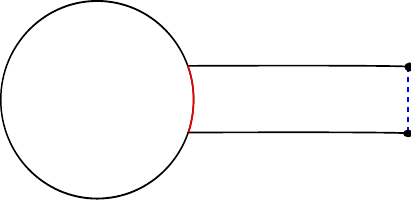}
\vspace{0.1cm}
\caption{ {\small A HR $h$ incident to a vertex disc. The two segments of $h$ are $s_1$ intersecting the vertex 
and $s_2$, the external segment with end points $a$ and $b$. The strands of the HR are the segments $[aa']$
and $[bb']$.}}
\label{fig:ribbonflag}
\end{minipage}
\put(-194,17){$s_1$}
\put(-149,17){$s_2$}
\put(-147,26){$a$}
\put(-147,4){$b$}
\put(-176,15){$h$}
\put(-196,28){$a'$}
\put(-196,2){$b'$}
\end{figure}
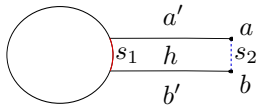

\begin{definition}[Ribbon half-edges and external points]
\label{ribflag}
A \emph{ribbon half-edge} (or simply half-ribbon, denoted henceforth HR) 
is a rectangle incident to a unique vertex of a ribbon graph by a unique line segment $s$ on the boundary and without forming a loop. The segment parallel to $s$ is called  \emph{the external segment}. The end-points of any external segment are called \emph{external points} of the HR.  The two boundary segments of a ribbon edge or of a HR 
that are neither external nor incident to a vertex are called \emph{strands}. A HR is always oriented consistently with the vertex it intersects. 
\end{definition} 
A HR incident to a vertex is drawn in Figure \ref{fig:ribbonflag}.
 The next definition was introduced in \cite{krf}. We reformulate
it according to our previous notation and definitions.
\begin{definition}[Cut of a ribbon edge]
\label{cutriedge}
Let $\cG$ be a ribbon graph and let $e$ be a ribbon edge of $\cG$.
The \emph{cut graph} $\cG \vee e$ is obtained from $\cG$ by 
deleting $e$ and attaching two HRs
at the same line segments where $e$ was incident
 to the end vertices, one at each of the end vertices of $e$. If $e$ is a loop, the two HRs are on the same  vertex.
\end{definition}

\begin{definition}\label{def:herg}
A HERG  $\cG(\cV,\cE,\mf^0)$ (or simply $\cG_{\mf^0}$) 
is a ribbon graph $\cG(\cV,\cE)$ (or shortly $\cG$) with a set $\mf^0$ of HRs
 such that each HR is attached to a unique vertex as in Definition \ref{ribflag}, and the segments where the HRs are attached 
are disjoint from each other and from the segments where any ribbon edges are attached. 
The ribbon graph $\cG$  is called the underlying ribbon graph of the HERG $\cG_{\mf^0}$. 
\end{definition}

  The cut of a ribbon edge $e$ in a HERG is the cut of $e$ in its underlying 
ribbon graph. Using Definitions \ref{subfla}
and \ref{def:herg}, spanning c-subgraphs of HERGs make sense.  A  spanning c-subgraph is then  formed by cutting some subset of the ribbon edges. We denote again the spanning c-subgraph inclusion as $A \sset \cG_{\mf^0}$. 
Note that a ribbon graph is a HERG with $\mf^0=\emptyset$.

Considered as geometric surfaces,  cutting an edge
in a ribbon graph or in a HERG modifies the boundary of that surface. 
 There are boundary components following the 
 boundary of the HR.  
 Combinatorially, we  want to distinguish this type of face and those which uniquely follow the boundary of ribbon edges.
 This will be also useful in the following section when 
considering the case of stranded graphs. 
The combinatorial objects that are defined below were introduced
in \cite{Gurau:2009tz}. We reformulate them using our conventions.

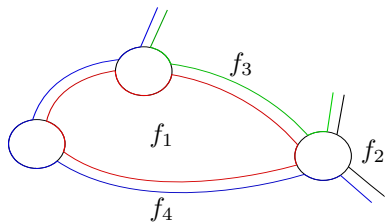
\begin{figure}[h]
 \centering
     \begin{minipage}[t]{.8\textwidth}
      \centering
\includegraphics[angle=0, width=5cm, height=2.6cm]{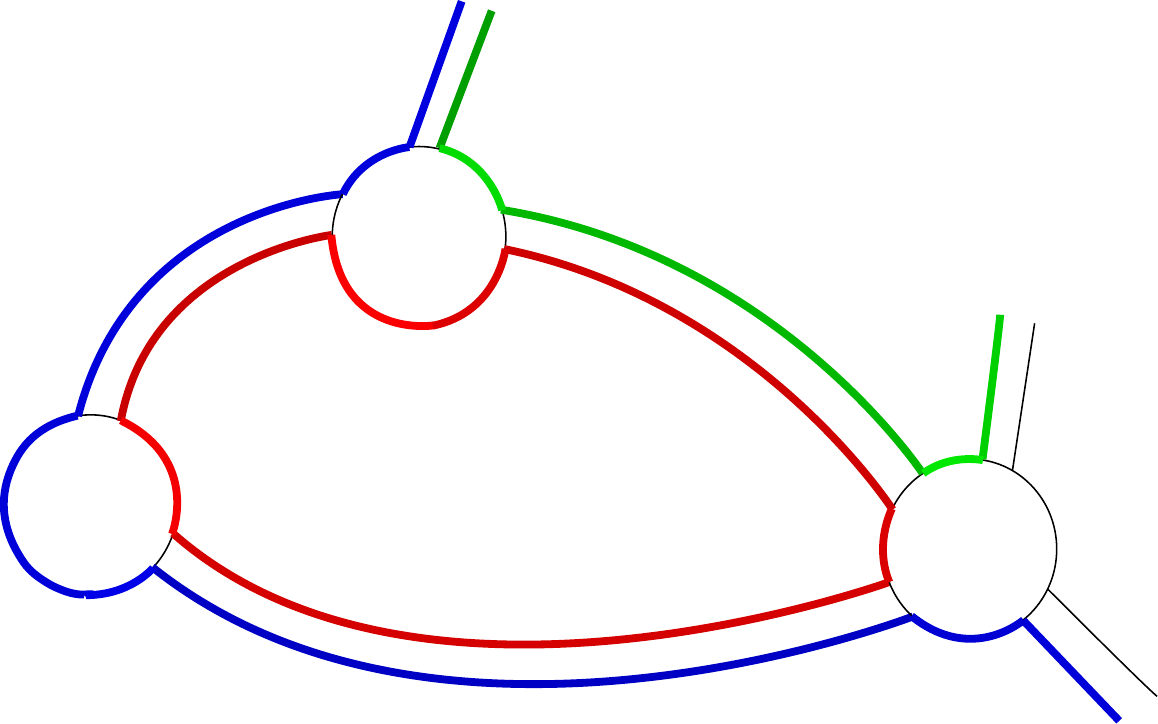}
\vspace{0.1cm}
\caption{ {\small A HERG with a closed face $f_1$ 
(in red) and open faces $f_{2}, f_{3},f_{4}$ (in black, green and blue, resp.).
}}
\label{fig:ribbonfla}
\end{minipage}
\put(-189,23){$f_1$}
\put(-109,17){$f_2$}
\put(-159,50){$f_3$}
\put(-189,-5){$f_4$}
\end{figure}

\begin{definition}[Closed, open faces]
\label{faces}
Consider  a HERG $\cG_{\mf^0}$.

$\bullet$ A \emph{closed face} is a boundary  component of $\cG_{\mf^0}$  which never passes through any external segment of a HR.  The set of closed faces is denoted $\cF_{\inter}$. (See the closed face $f_{1}$ in 
Figure \ref{fig:ribbonfla}.)

$\bullet$ An \emph{open  face}  is a boundary arc
leaving an external point of some HR rejoining another external point
 without passing through any external segment of a HR.  
 The set of open faces is denoted  $\cF_{\ext}$. 
  (Examples of open faces are provided in Figure \ref{fig:ribbonfla}.)

$\bullet$ The set of faces $\cF$ of a HERG  is defined by 
$\cF_{\inter} \cup \cF_{\ext}$.

\end{definition}

\begin{figure}[h]
 \centering
     \begin{minipage}[t]{.8\textwidth}
      \centering
\includegraphics[angle=0, width=5cm, height=2.6cm]{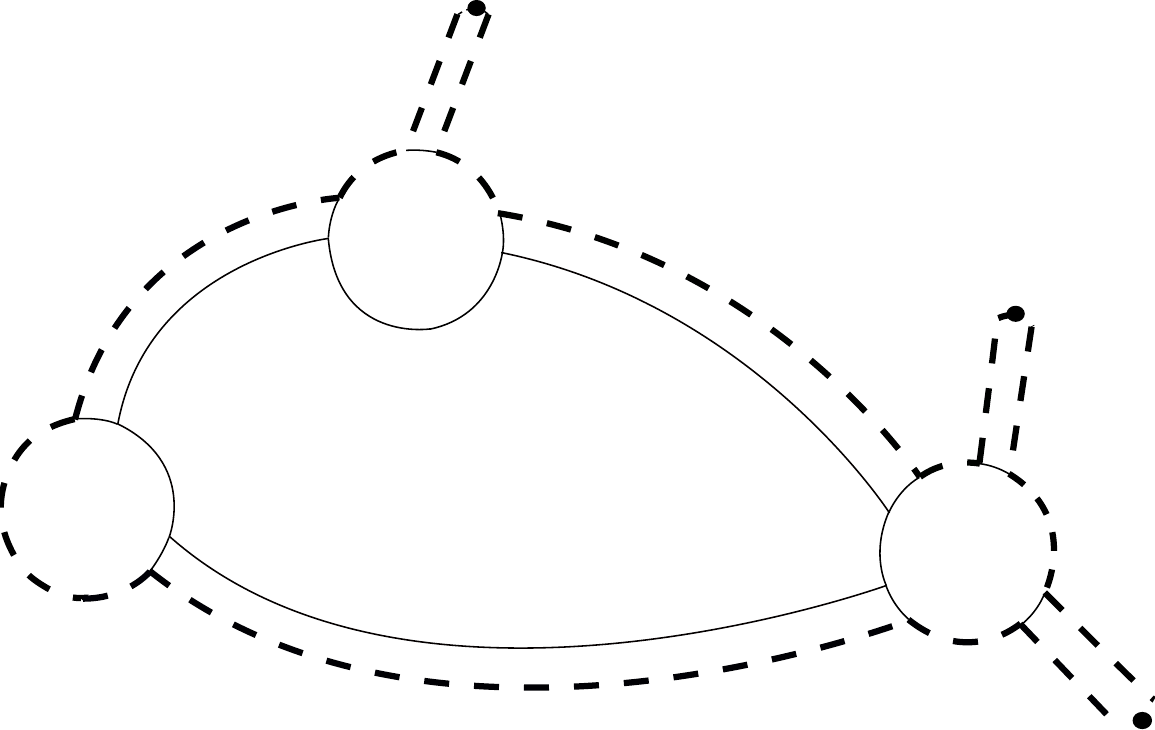}
\vspace{0.1cm}
\caption{ {\small The boundary graph $\bG$ of $\cG_{\mf^0}$ of Figure \ref{fig:ribbonfla} is the dashed cycle.}}
\label{fig:boundary}
\end{minipage}
\end{figure}

\begin{definition}[Boundary graph \cite{Gurau:2009tz}]
\label{bnd}
The \emph{boundary graph} $\bG$  of a HERG $\cG_{\mf^0}$ is an abstract graph $\bG(\bV,\bE)$ 
such that $\bV$ is in one-to-one correspondence with $\mf^0$,
and $\bE$ is in one-to-one correspondence with $\cF_{\ext}$. 
Consider an edge $e$ of $\bE$, its corresponding open face $f_e \in\cF_{\ext} $, a vertex $v$, and its corresponding 
HR $h_v$. The edge $e$ is incident to $v$ if and only if  $f_e$ has one end-point in $h_v$,  and,  if both end-points of $f_e$ are in $h_v$, then $e$ is a loop.   (The boundary of the graph given in Figure \ref{fig:ribbonfla}
is provided in  Figure \ref{fig:boundary}.)
\end{definition} 
The boundary  graph  of a ribbon graph is empty. 
Note that $\bG$ is a disjoint union of cycles and
also that the connected components of $\bG$ are in one-to-one 
correspondence with the faces of $\cG$ that do not correspond
to the closed faces of $\cG_{\mf^0}$.

The notion of edge contraction and deletion for HERGs can be simply understood as in the case of ribbon graphs. Let $\cG_{\mf^0}$ be a HERG and $\cG$ its underlying ribbon graph. 
In notation of Definition \ref{defBRpoly}, we now define $r(\cG_{\mf^0})=r(\cG),\, n(\cG_{\mf ^0})=n(\cG)$, $o(\cG_{\mf^0})=o(\cG)$.

\begin{definition}[BR polynomial for HERGs]
Let $\cG(\cV,\cE,\mf^0)$ be a HERG. We define the BR polynomial of $\cG_{\mf^0}$ to be
\beq
\cR_{\cG_{\mf^0}}(x,y,z,s,w,t)
=\sum_{A\sset \cG_{\mf^0}} (x-1)^{r(\cG_{\mf^0})-r(A)}(y-1)^{n(A)}
z^{k(A)-F_{\inter}(A)+n(A)} \, s^{C_\partial(A)} w^{o(A)}\,   t^{f(A)},
\label{brfla}
\eeq
where $F_{\inter}(A)= |\cF_{\inter}(A)|$, 
$C_\partial(A)= |\bC(A)|$ is the number 
of connected components of $\partial A$, the boundary graph of $A$, 
 $o(A)=0$ if $A$ is orientable, and 1 otherwise,
$f(A)=|\mf^0(A)|$ is the number of HRs of $A$,
and where $w^2=w$ holds.
\end{definition}

The polynomial $\cR$ \eqref{brfla} generalizes the BR polynomial 
$R$ \eqref{brpoly}.
 $R$
 can be recovered from $\cR$  using 
\beq\label{crcrp}
\cR_{\cG_{\mf^0}}(x,y,z,z^{-1},w,t=1) =  
R_{\cG}(x,y,z,w)\,. 
\eeq

The graph operations of disjoint union and one-point join extend to ribbon graphs  and to HERGs.
The one-point join $\cG_{1\mf_1^{0}} \cdot_{v_1,v_2} 
\cG_{2\mf_2^{0}}$ of disjoint HERGs $\cG_{1\mf_1^{0}}$ and $\cG_{2\mf_2^{0}}$ is obtained by choosing two vertices $v_1$ and $v_2$ of $\cG_{1\mf_1^{0}}$ and $\cG_{2\mf_2^{0}}$, respectively, and merging $v_1$ and  $v_2$ on an arc of each of these which does not contain any ribbon edges or HRs.
Combinatorially, we must respect the cyclic order of all ribbon edges and HRs on the previous vertices $v_1$ and $v_2$. It is well-known  that such a procedure might lead to non isomorphic  ribbon graphs,  see for instance \cite{bollo}.
 The fact that $R_{\cG_1 \sqcup \cG_2}
= R_{\cG_1} R_{\cG_2} = R_{\cG_1 \cdot_{v_1,v_2} \cG_2}$ holds
for ribbon graphs  can be 
extended  for HERGs under certain conditions. 

We come to the properties of  $\cR_{\cG_{\mf^0}}$. 
The following proposition holds.
\begin{proposition}[Union of HERGs] 
\label{opBRfla}
Let $\cG_{1\mf^0_1} $ and $\cG_{2\mf^0_2}$ be two disjoint HERGs, then 
\beq 
\cR_{\cG_{1\mf^0_1} \sqcup \cG_{2\mf^0_2}} = \cR_{\cG_{1\mf^0_1} }
 \cR_{\cG_{2\mf^0_2}}\,. 
\label{brcups1}
\eeq
\end{proposition}
\proof 
It suffices to check the exponents for any $A\sset \cG_{1\mf^0_1} \sqcup \cG_{2\mf^0_2}$.    A spanning c-subgraph $A$ of such a union of HERGs expresses as $A = A_1 \sqcup A_2 \sset  \cG_{1\mf^0_1} \sqcup \cG_{2\mf^0_2}$ 
with $A_{1} \sset \cG_{1\mf^0_1}$ and $A_{2} \sset \cG_{2\mf^0_2}$. It is straightforward to see that 
$k(A)$, $r(A)$, $F_{\inter}(A)$, $n(A)$, $C_\partial(A)$, and  $f(A)$  are
all additive quantities. Furthermore, $o(A_1 \sqcup A_2)=\max\{o(A_1),o(A_2)\}$,
but $w^2=w$. The result follows.   
 
\qed

\begin{theorem}[Contraction and cut for BR polynomial for HERGs]
\label{theo:BRext}
Let $\cH= \cG(\cV,\cE,\mf^0)$ be a HERG. Then,
for an   ordinary edge $e$, 
\beq
\cR_{\cH}=\cR_{ \cH\vee e} +\cR_{\cH/e}\,,
\label{retcondel}
\eeq 
for a bridge $e$, we have 
\beq
\cR_{\cH} =(x-1)\cR_{\cH \vee e}+ \cR_{\cH/e}\,;
\label{retbri}
\eeq 
 for a trivial  negative loop $e$, the following holds 
\beq
\cR_{\cH}=\cR_{\cH \vee e} + (y-1)zw \,\cR_{\cH/e}\,,
\label{retsel2}
\eeq 
whereas for a trivial positive loop $e$, we have
\beq
\cR_{\cH}=\cR_{\cH \vee e} + (y-1) \cR_{\cH/e}\,.
\label{retsel1}
\eeq 
\end{theorem}  
\proof  In the following proof, 
spanning c-subgraphs are simply called subgraphs because
no confusion can arise. 
Let us make two preliminary remarks. 
(A) The subgraphs $A$ of { $\cH$} which do not contain $e$ are precisely the subgraphs of { $\cH\vee e$}. 
(B) Also, if $e$ is not a loop then the map $A\mapsto A/e$ provides a bijection from the subgraphs of { $\cH$} which contain $e$ to the subgraphs of { $\cH/e$}  that preserves closed faces as well as components  of the boundary graph. 
 Note importantly that, although  $A$ and $A/e$ do not have
the same vertices and edges, they do have the  same
HRs and same faces.  

Let us prove \eqref{retcondel}. For an   ordinary edge $e$,
let $A \sset { \cH}\vee e$,  and let $A'$  be the corresponding subgraph  in { $\cH$} such that 
 $e \notin A'$, by remark (A). 
The fact that the monomial of $A$ in $\cR_{{\cH}\vee e}$
and the monomial corresponding  to $A'$ 
in $\cR_{{ \cH}}$ are identical is simple to check.  Then
 $\sum_{A\sset { \cH}; e \notin A }(x-1)^{r({ \cH})-r(A)}(y-1)^{n(A)}
z^{k(A)-F_{\inter}(A)+n(A)} \, s^{C_\partial(A)} w^{o(A)}\, t^{f(A)} = \cR_{{ \cH}\vee e}$. 

We now concentrate on the remaining  sum related to the contraction
of $e$. In  particular, we focus on the sets of faces during the contraction.  
Choose $A\sset { \cH}$ with $e \in A$ and, by remark (B), let 
 $A' \sset { \cH}/e$ be its corresponding subgraph.
One has
\bea
F_{\inter}(A)=F_{\inter}(A')\,,\;  C_\partial(A)=C_\partial(A')  \,,\;
o(A)= o(A')\,,\,\mbox{ and } \,\,  f(A) =  f(A') \,. 
\eea
The monomial of $\cR_{{ \cH}/e}$ related to {$A'$} is of the form 
 \bea
&&
(x-1)^{r({ \cH}/e) - r(A')} (y-1)^{n(A')}z^{k(A')-F_{\inter}(A')+n(A')} \, s^{C_\partial(A')} \, w^{o(A')} t^{f(A')} \crcr
&&= (x-1)^{r({ \cH})-1 - (r(A)-1)} (y-1)^{n(A)}
 z^{k(A)-F_{\inter}(A)+n(A)} w^{o(A)}\, s^{C_\partial(A)} \,  t^{f(A)}
\eea 
which achieves the proof that  $\cR_{{ \cH}/e}=\sum_{A \sset { \cH}; e\in A} (x-1)^{r({ \cH}/e) - r(A)} (y-1)^{n(A)}z^{k(A)-F_{\inter}(A)+n(A)} $
$ s^{C_\partial(A)}w^{o(A)} t^{f(A)}$, and then \eqref{retcondel}  holds.

We now focus on \eqref{retbri}. 
Let $e$ be a bridge in ${ \cH}$. Decompose $\cR_{{ \cH}}$ as
 \beq
\sum_{A\sset { \cH}; e \notin A }(x-1)^{r({ \cH})-r(A)}(y-1)^{n(A)}
z^{k(A)-F_{\inter}(A)+n(A)} \, s^{C_\partial(A)} w^{o(A)}\,  t^{f(A)}+ \cR_{{ \cH}/e}\,.
\eeq 
It remains to prove that the first sum corresponds to  
$(x-1)\cR_{{ \cH} \vee e}$ 
but this is straightforward from 
$r({ \cH}) = r({ \cH}\vee e)+1$ and since all other terms remain
unchanged. 

The proofs of relations \eqref{retsel2} and \eqref{retsel1} are now given. 
Consider a trivial (positive or negative) loop  $e$ in ${ \cH}$, 
then  $\sum_{A\sset { \cH}; e \notin A }(x-1)^{r({ \cH})-r(A)}(y-1)^{n(A)}
z^{k(A)-F_{\inter}(A)+n(A)} \, s^{C_\partial(A)} w^{o(A)}\,  t^{f(A)} =\cR_{{ \cH} \vee e}$
still  holds in any case. The mapping from the subgraphs of ${ \cH}$ containing $e$ to those of ${ \cH}-e$, or conversely, is just 
obtained by deleting $e$ or gluing $e$ to the corresponding subgraph.

 Consider  a trivial  negative loop  $e$ in $\cH$. To each $A\sset \cH$ such that $e \in A$ 
and its corresponding $A' \sset \cH/e$, 
we find that $n(A)=n(A')+1$, $F_{\inter}(A) = F_{\inter}(A')$, 
and $C_\partial(A) = C_\partial(A')$. Therefore,
we get the following relation between  the terms:
\bea
&&
(x-1)^{r(\cH) -r (A) } (y-1)^{n(A)} 
 z^{k(A)-F_{\inter}(A)+n(A)} \, s^{C_\partial(A)} \,w^{o(A)}\, t^{f(A)}
\crcr
&&= (x-1)^{r(\cH/ e) -r (A') } (y-1)^{n(A')+1} 
 z^{k(A')-F_{\inter}(A')+n(A')+1} \, s^{C_\partial(A')}\,w^{o(A')+1} \,  t^{f(A')}\,.
\eea
Thus \eqref{retsel2} is satisfied. 

We now suppose that $e$ is a trivial  positive loop. 
With $A\sset \cH$ such that $e \in A$, 
we associate a unique corresponding element $A'$ in $ \cH /e$. We can infer that $k(A)=k(A')-1$, $n(A)=n(A')+1$, $F_{\inter}(A) = F_{\inter}(A')$, 
and $C_\partial(A) = C_\partial(A')$. 
Thus, the following relation between the terms corresponding
to $A$ and $A'$ holds: 

\bea
&&
(x-1)^{r(\cH) -r (A) } (y-1)^{n(A)} 
 z^{k(A)-F_{\inter}(A)+n(A)} \, s^{C_\partial(A)}\, w^{o(A)} \,  t^{f(A)}
\crcr
&&  (x-1)^{r(\cH/e) -r (A') } (y-1)^{n(A')+1} 
 z^{(k(A')-1) - F_{\inter}(A')+(n(A')+1)} \, s^{C_\partial(A')} \,  w^{o(A')}\, t^{f(A')}\,,
\eea
so that \eqref{retsel1} is obtained.  

\qed

 Observe that, as in the case of the of the BR polynomial for ribbon  graphs, Theorem \ref{theo:BRext} is not a complete reduction, since one-vertex
HERGs with nontrivial loops must still be computed from \eqref{brfla}. 
 Given $\cG_{\mf^0}$ a HERG and $\cG$ its underlying ribbon graph, we also note the reduction
\beq\label{crpBR}
\cR_{\cG_{\mf^0}}(x,y,z,s=z^{-1},w,t) = t^{|\mf^0| + 2n(\cG)}
R_{\cG}(X,\frac{Y}{t^2},z,w)
\eeq
where, once again, $X=t^2(x-1)+1$ and $Y=y+t^2-1$.

\section{Rank $D$ 
half-edged stranded graphs and a generalized polynomial invariant}
\label{sect:tgraph}

This section investigates the definition of a new polynomial
for particular graphs which aims at generalizing the BR polynomial $\cR$ for HERGs. 
The main notion of graphs discussed below is
combinatorial and can be always pictured in a $3D$ space. 
 Our main result appears in Theorem \ref{theo:contens}
after defining the generalized graphs we are dealing with.

The  notion of rank $D$ colored tensor graphs considered here has been introduced by Gurau in \cite{Gurau:2009tw}.  Using a duality,  it is well-known that ribbon graphs can be mapped onto triangulations of surfaces. In a similar way, colored tensor graphs can be interpreted as simplicial complexes
or dual of triangulations of  topological spaces in any dimension. 
They are of particular interest in certain quantum field theories  defined with tensor fields hence the name 
tensor graphs (ribbon graphs are, in this sense, rank 2 or matrix graphs).
The importance of these graphs has been highlighted 
by Gurau in \cite{Gurau:2010nd}  where that author proved that  colored tensor graphs 
have a cellular structure which associates each colored tensor graph with a simplicial pseudo-manifold of $D$ dimensions.  It has been also proved 
that a colored tensor graph which is bipartite 
induces naturally an orientation of the 
dual simplicial complex in \cite{Caravelli:2010nh}. 

Some previous studies have addressed the  generalization of the BR polynomial for higher dimensional 
objects within the framework of such graphs. Mainly, two authors, Gurau in \cite{Gurau:2009tz} and Tanasa in \cite{Tanasa:2010me}, 
have defined two  distinct notions of generalized BR polynomials. 
 Let us give a brief review of their results and compare 
these to the one obtained in the present work.  

 Gurau defined a multivariate polynomial invariant for colored tensor  graphs associated with simplicial complexes with boundaries in any dimension $D$. The polynomial that we obtain in the present work can be put in a multivariate form which is related to Gurau's polynomial restricted to 3$D$. Our polynomial extends Gurau's polynomial  to a class larger than colored tensor graphs.   We emphasize 
that the difficulties encountered by Gurau (as well as Tanasa, see below) 
for defining a contraction procedure for such type of graphs
without destroying the entire graph structure  
is much improved in our scheme. 

In a different perspective, the work by Tanasa in \cite{Tanasa:2010me} 
deals with tensor graphs without colors  that are equipped with
 another  stranded vertex. The polynomial as worked out by this author is only valid
for graphs triangulating topological objects without boundary. The polynomial that we define is radically different from that one in several features, since mainly, it relies on a graph coloring.

In another close related setting dealing with simplicial complexes, Krushkal and Renardy in \cite{kru2} identified a 
4-variable polynomial invariant for triangulations and handle decompositions of orientable manifolds reducing to the BR polynomial.
This polynomial might agree with the one that we will define
for closed triangulations after specializing to 1 all
variables related to the boundary of the complex. This 
 deserves as well to be fully addressed elsewhere.

\subsection{Stranded, tensor, colored graphs with half-edges}
\label{subset:tgraphs}

Like ribbon graphs, colored tensor graphs have both a topological meaning (realized in the dual triangulation) and a combinatorial formulation that is  our main concern here.  
Before reaching the definition of colored tensor graphs, we must introduce the  
combinatorial concept of stranded graphs, the true backbone of this theory. 

\vspi

\noindent{\bf Stranded and tensor graphs.} 
We start by  basic notions of stranded objects.

\begin{definition}[Stranded vertex and edge]\label{def:svert}
 Let  $D$ be an integer, $D \ge 0$. 

\noindent $\bullet$   A \emph{stranded pre-vertex} of rank $D\ge 1$ consists of a set $S$ of $2n$ elements called \emph{(vertex) points}, together with two partitions of $S$: one into non-empty parts of size at most $D$, called \emph{pre-edges}, and one into pairs, called \emph{chords}. 
A rank $D=0$ stranded vertex consists of a set $S$ with a single or no element,
and no partition of $S$ (then there are no  pre-edges or chords).

The case $n=0$ is allowed (then there are no pre-edges, or chords),  
and the resulting stranded pre-vertex is called a  \emph{trivial circle} 
which, by convention, we consider of any rank $D \ge 0$.  

The \emph{coordination number} of $v$ is the number of pre-edges in $v$. 
We say that a stranded pre-vertex $v$ is \emph{disconnected} if we can partition the set $S$ into non-empty parts $S_1$ and $S_2$ such that no pre-edge or chord contains elements of both $S_1$ and $S_2$, otherwise, $v$ is \emph{connected}.  
A rank $D$ \emph{stranded vertex} is a rank $D$ stranded pre-vertex that is connected.

\noindent $\bullet$ The \emph{vertex graph} of a rank $D$ stranded pre-vertex
is the graph whose vertex set is the set of pre-edges and with one edge for each chord. 
The edge is incident to two vertices (or one in a loop case), corresponding to  the two pre-edges (or one in a loop case) that contain the two paired points. 

\noindent $\bullet$  A \emph{stranded edge of rank} $D\ge 1$  consists of a set $S$ of $2D$ elements called \emph{(edge) points},  together with two partitions, one into two sets of size $D$ called \emph{ends}, and one into $D$ sets of size $2$ called \emph{strands}, each strand containing one point from each end.

\end{definition}

 Note that a stranded vertex is a stranded pre-vertex whose vertex graph is connected. 

\begin{figure}[h]
 \centering
     \begin{minipage}[t]{1\textwidth}
      \centering
\includegraphics[angle=0, width=10cm, height=3cm]{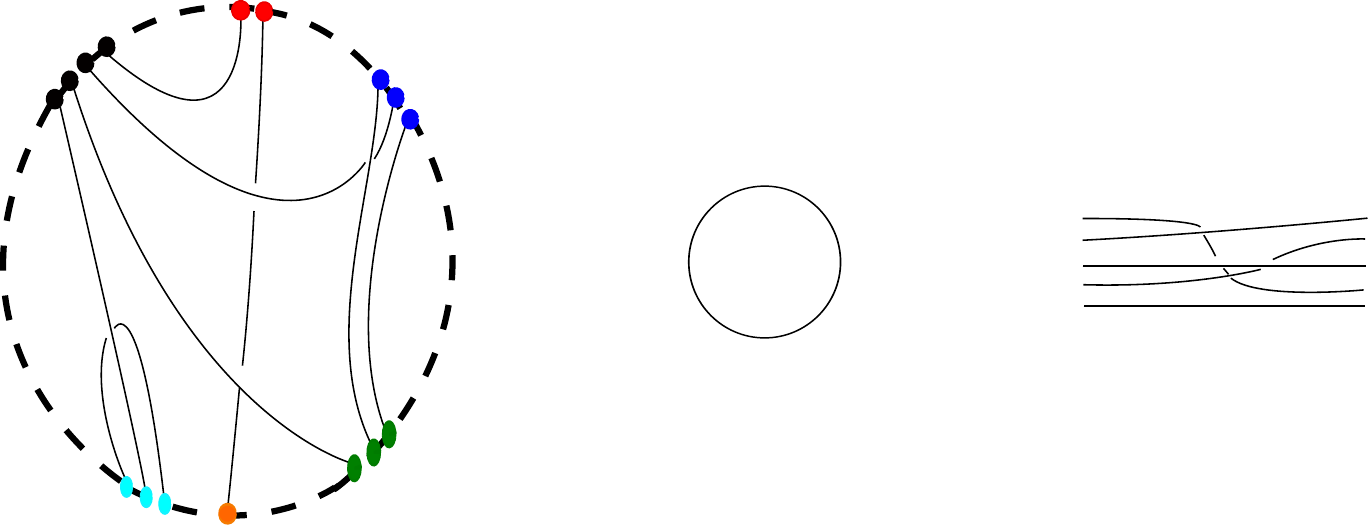}
\vspace{0.1cm}
\caption{ {\small Graphical representation of stranded
vertices and edges: a rank 4 stranded vertex $v$ of coordination 6 (left), with
 pre-edges (highlighted with different colors) with non intersecting chords; a trivial circle vertex $d$; a rank 5 
 stranded edge $e$ with non parallel strands (right).}}
\label{fig:vxedge}
\end{minipage}
\put(-315,-10){$v$}
\put(-200,5){$d$}
\put(-100,5){$e$}
\end{figure}

In the rest of this work,  we use a particular realization for drawing  stranded vertices and edges which does not play any 
important role.   The chords of stranded vertices (resp. strands
of stranded edges) are represented line segments  that do not intersect and their crossings are irrelevant.  The  end points of chords are partitioned
into sets representing the pre-edges with $1,2,\dots$ or $D$ points. These points are drawn on a single arc of a 
fictitious circle, called the  \emph{vertex frontier},  with no other end points on this arc.  
Drawn in dash or dots, or as a solid circle when no confusion can arise, 
the vertex frontier is used to separate in drawings stranded vertices and edges  in the (yet to be defined) stranded graphs. 
 The cyclic order around the dotted circle does not
matter, and this circle plays a different role from the vertex circles in ribbon graphs. 
Examples of stranded vertices and edges with rank $D=4$ and $D=5$ have been provided in Figure \ref{fig:vxedge}.

\begin{definition}[Stranded and tensor graphs]\label{def:tens}
$\bullet$  A \emph{rank $D$ stranded graph} $\crrG=\crrG(\cV,\cE)$
consists of a set $\cV$ of disjoint rank $D$ stranded vertices (i.e., with all vertex points distinct), a set $\cE$ of disjoint stranded edges
(i.e., with all edge points distinct)
of rank at most $D$, and an incidence map $\phi$ (not usually written in the notation) from the set of edge points to the set of vertex points satisfying the following conditions: $\phi$ is injective, and on each end of each stranded edge, $\phi$ acts as a bijection to some pre-edge of some stranded vertex in $\cV$.
In other words, each end of each stranded edge is attached to a distinct pre-edge of the correct size, with the $D$ strands in each rank $D$ edge attached bijectively to the $D$ points in the pre-edge.

$\bullet$ A rank $D$ tensor graph $\crrG$ is a rank $D$ stranded graph such that 
the stranded vertices of $\crrG$ have a fixed coordination $D+1$  and their  pre-edges have a fixed   cardinality $D$. 
 All vertex graphs are  $K_{D+1}$. 
The stranded edges of $\crrG$ are of rank $D$. 
\end{definition}

\begin{figure}[h]
 \centering
     \begin{minipage}[t]{1.1\textwidth}
      \centering
\includegraphics[angle=0, width=7cm, height=3cm]{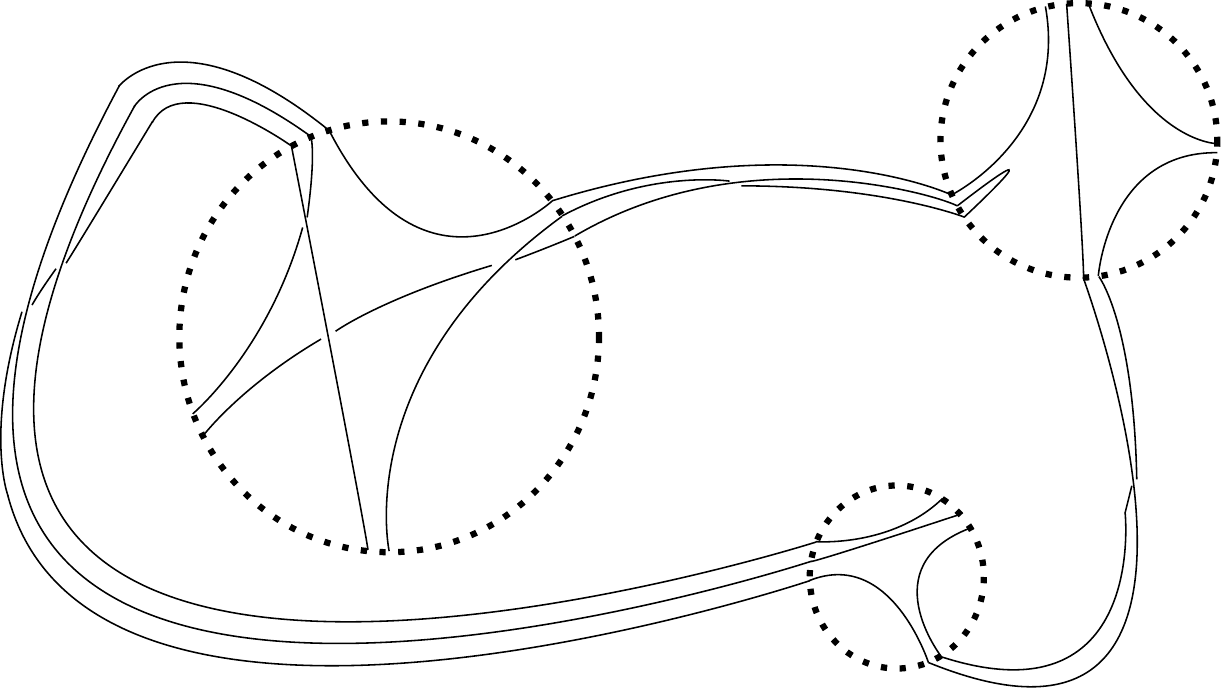}
\caption{ {\small A rank 3 stranded graph.}}
\label{fig:strangra}
\end{minipage}
\end{figure}

\begin{figure}[h]
 \centering
     \begin{minipage}[t]{1.1\textwidth}
      \centering
\includegraphics[angle=0, width=7cm, height=3cm]{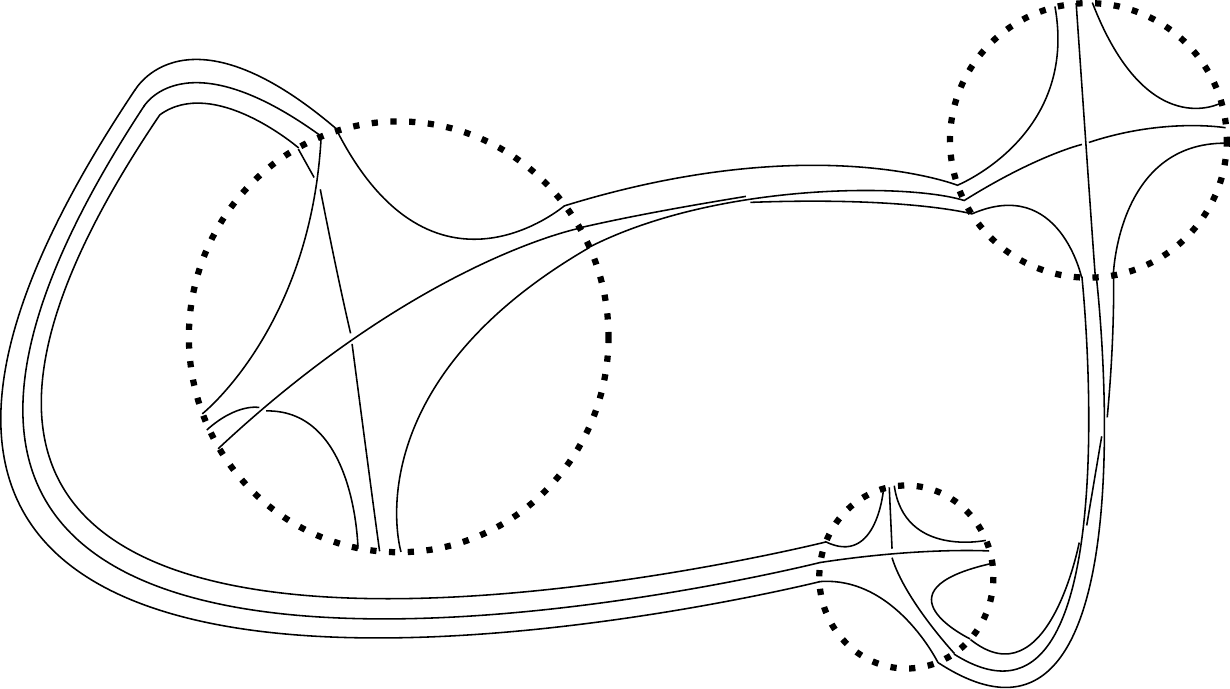}
\caption{ {\small A rank 3 tensor graph with rank 3 stranded vertices as  with 
fixed coordination 4, pre-edges with 3 points linked by chords according
to the pattern of $K_4$; stranded  edges are rank 3.}}
\label{fig:tensgra}
\end{minipage}
\end{figure}

Some illustrations of a rank 3 stranded and tensor graphs are given in Figure \ref{fig:strangra} and \ref{fig:tensgra}, respectively.
It should be pointed out that,  although the stranded vertices and edges of stranded and tensor graphs are drawn in a three dimensional space, we do not treat these as embedded graphs.

There are motivations for the introduction of tensor graphs. 
While stranded graphs can be treated as  generalized graphs, 
 tensor graphs, via a combinatorial duality, actually map to 
simplicial spaces. 
Consider $\crrG$ a rank $D>0$ tensor graph. 
A stranded vertex of $\crrG$ represents a $D$-simplex
and the pre-edges are precisely their boundary 
 $(D-1)$-simplices. 
A stranded edge of $\crrG$ represents the gluing of the  $D$-simplex
along one of their boundary  $(D-1)$-simplices. 
Hence, a tensor graph $\crrG$ maps to a simplicial complex. 
As an illustration of the above combinatorial duality,
a rank 3 tensor graph represents a simplicial 
complex in 3$D$ composed by tetrahedra (3-simplex) which are glued along their boundary triangles 
(2-simplex).
Such a duality has given a handle on the study of random
simplicial manifolds in physics topics like quantum gravity \cite{ambjorn}.
The interpretation of the circle vertices in the definition of tensor graphs
can be done in the following convention: 
a trivial circle vertex in a rank $D \geq 2$ tensor
graph represents a $D$ dimensional ball.

\vspi

\noindent{\bf  Underlying graph and connectedness.} 
Given a stranded graph
$\crrG$ and collapsing its stranded vertices to points and its stranded edges to simple lines,  
 one obtains an abstract graph called \emph{underlying graph} of $\crrG$. 
 A stranded graph is \emph{connected} if and only if its corresponding underlying graph is connected. 
For instance, both graphs of Figures \ref{fig:strangra}
and \ref{fig:tensgra}  are connected since  their underlying graphs
coincide with the closed cycle graph $C_3$. 

Any rank $D$ stranded graph 
  is a rank $D' $ stranded graph for any $D' \geq D$. This is however
not true for tensor graphs. In the following, 
the rank of a stranded graph refers to the minimal rank $D$ for which
this stranded graph is well defined.

\vspi

\noindent{\bf Equivalence class of stranded graphs.}
 Two stranded graphs are  \emph{isomorphic} if there are two bijections, one from the vertex points of the first to the vertex points of the second, and one from the edge points of the first to the edge points of the second, that preserve the structure, meaning the grouping into pre-edges, chords, edge ends and strands, as well as respecting the edge point to vertex point incidence maps. 
 Furthermore these two graphs should have the same number of trivial circles.

\begin{figure}[h]
 \centering
     \begin{minipage}[t]{1\textwidth}
      \centering
\includegraphics[angle=0, width=8cm, height=7.5cm]{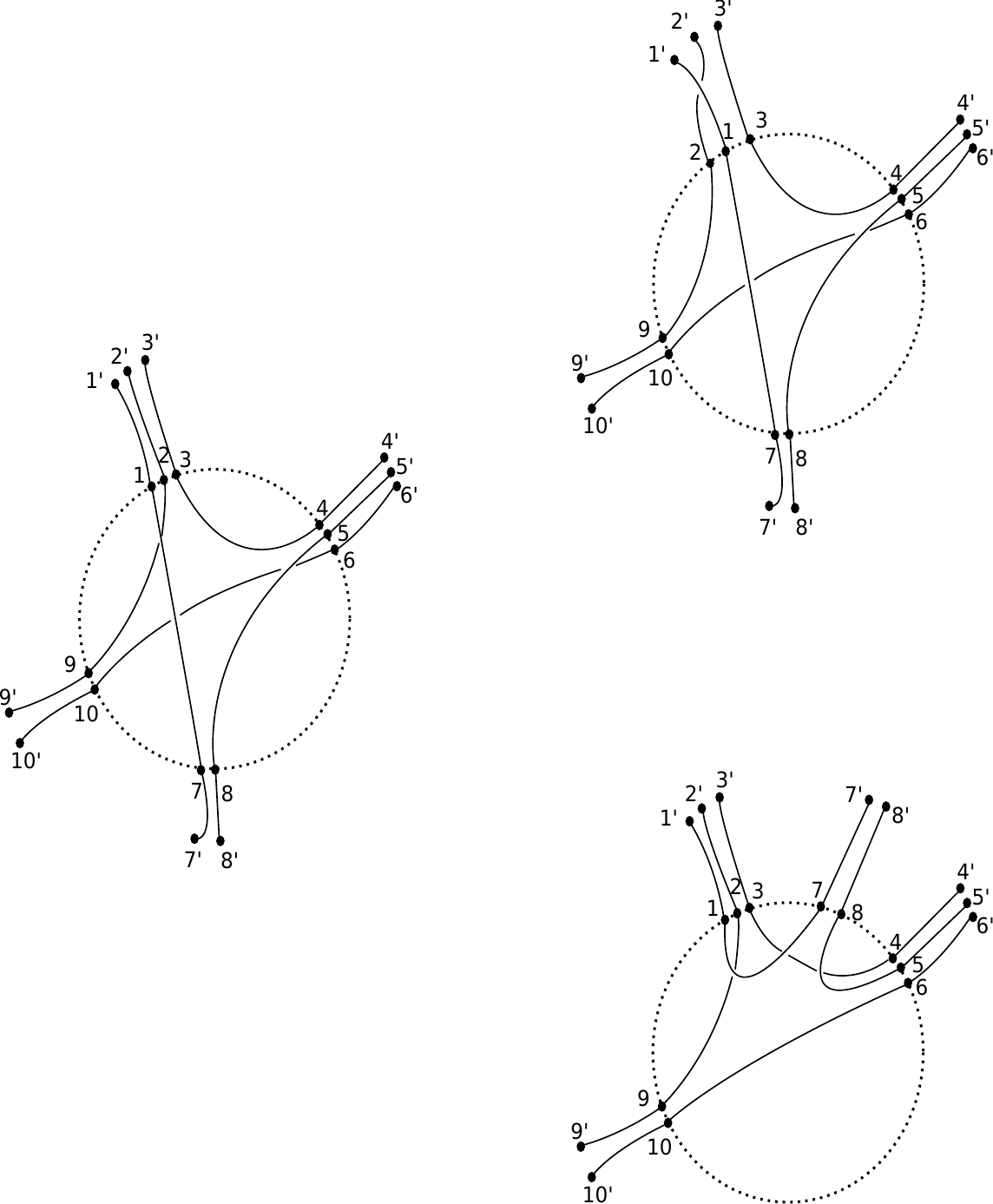}
\vspace{0.1cm}
\caption{ {\small Different drawings of the same part of a stranded graph. }}
\label{fig:samevert}
\end{minipage}
\put(-305,162){$e_1$}
\put(-230,135){$e_2$}
\put(-283,50){$e_3$}
\put(-340,75){$e_4$}
\put(-185,210){$e_1$}
\put(-98,190){$e_2$}
\put(-150,108){$e_3$}
\put(-208,135){$e_4$}
\put(-185,75){$e_1$}
\put(-98,60){$e_2$}
\put(-125,78){$e_3$}
\put(-208,5){$e_4$}
\end{figure}

 Figure \ref{fig:samevert} illustrates different ways of drawing part of the same stranded graph: the stranded vertex is defined by the pre-edges $\{1, 2, 3\}$, $\{4, 5, 6\}$, $\{7, 8\}$, $\{9, 10\}$ and chords $\{1, 7\}$, $\{2, 9\}$, $\{3, 4\}$, $\{5, 8\}$,  $\{6, 10\}$ and (in the edges) the strands are $\{i, i'\}$; $i=1, \cdots, 10$.

\vspi

\noindent{\bf Low rank stranded graphs.} 
As illustrations, we now discuss  the lowest rank stranded 
graphs.

\quad (0)  Rank 0 stranded or tensor graphs are abstract graphs 
made with  only vertex points as  rank 0 stranded vertices, no stranded edges, and possibly with trivial circles.

\quad (1) 
Rank 1 stranded/tensor graphs have both  stranded vertices and edges like segments (possibly with additional trivial circles). 
Note that the connectivity condition implies that a rank 1 stranded vertex can only 
have $n=0$ (hence is a trivial circle) or $n=1$ (hence have two pre-edges of size 1).
Such rank 1 stranded graphs cannot be 
directly identified with abstract graphs.
Hence, neither rank 0 nor rank 1 stranded graphs
define abstract graphs in general. 
Abstract graphs  
can be directly obtained from rank $D$ stranded
graph through
the collapsing procedure described above,
hence using the underlying graph.

\quad (2)  
Ribbon graphs, in the sense of Definition \ref{def:ribbongraph},  are in one-to-one 
 correspondence  with  certain rank 2 stranded graphs.
In the following, we show (a) how to construct the rank 2 stranded 
graph corresponding to a ribbon graph and (b) prove that
two equivalent ribbon graphs lead to two equivalent rank 2 stranded graphs in the sense identified above.

\begin{figure}[h]
 \centering
     \begin{minipage}[t]{.8\textwidth}
      \centering
\includegraphics[angle=0, width=6cm, height=2cm]{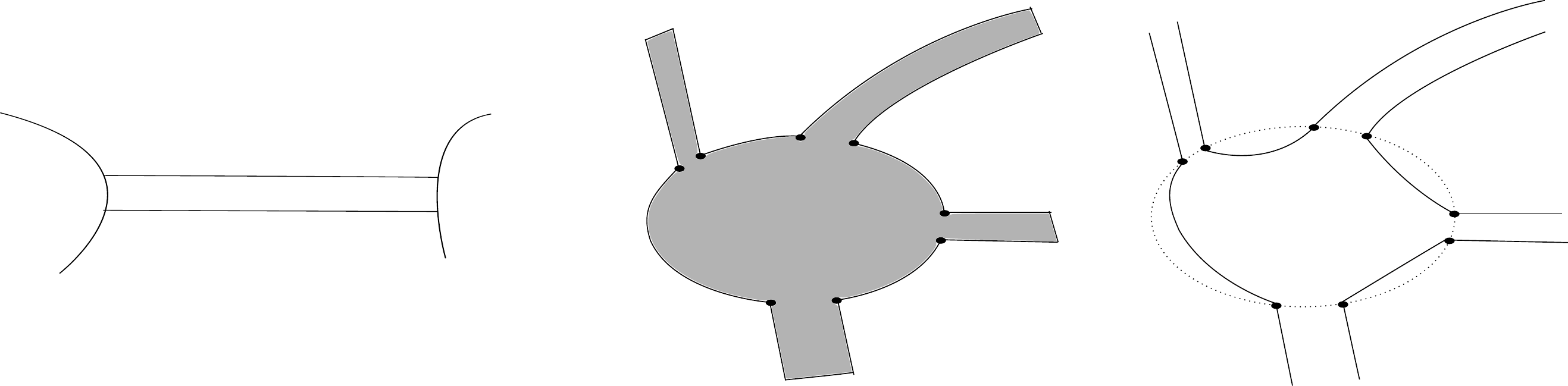}
\caption{ {\small Ribbon edge and vertex of a ribbon graph as a rank 2 stranded edge and vertex of a stranded graph (the frontier
vertex appears in dash is fictitious).}}
\label{fig:derib}
\end{minipage}
\put(-223,17){$s_1$}
\put(-223,35){$s_2$}
\put(-270,25){A}
\put(-70,25){B}
\end{figure}
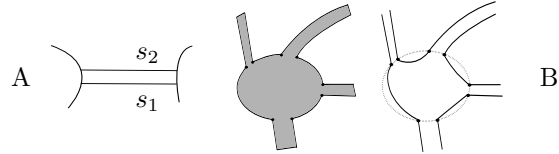

(a) Let $\cG$ be a ribbon graph  (Definition \ref{def:ribbongraph}) and let $e$ be an edge of $\cG$. The ribbon edge  $e$ is incident to
its end vertex (loop case) or vertices (ordinary case or bridge) in two segments $s$ and $s'$.  
Consider the  two distinct segments on the boundary 
of $e$ which are not $s$ and $s'$ (denoted by $s_1$ and $s_2$
in Figure \ref{fig:derib}A).  
These two segments define the two strands of a rank 2 stranded edge with 
$s$ and $s'$ as end segments. 
  By convention, a  circle is a stranded vertex in any rank, in particular $D= 2$,
and so it is a stranded graph with no incident stranded edge.
Consider now a vertex $v$ of $\cG$ with incident edges 
$e_1$, $e_2$, \dots, $e_p$, $p \in \mathbb{N}$.
We can construct a stranded vertex $v'$ of rank $2$  from the data of $v$
and of its incident ribbon edges  in the following manner. Draw a fictitious circle as the frontier vertex of $v'$. 
Insert on the frontier the end points of end segments of $e_k$, $1\leq k \leq p$, which define the pre-edges (with exactly 2 points) of $v'$. 
From these pre-edges, one defines the chords of $v'$ 
as the segments disposed cyclically between the end segments 
of $e_k$'s which lie in $v$. See Figure \ref{fig:derib}B.  
We stress again that the vertex frontier is totally virtual
and is useful only for determining the separation
between stranded edges and vertices. We can conclude
that a rank 2 stranded graph obtained from a ribbon
graph,  in this way,  is a combinatorial object defined on the boundary of the ribbon graph.

 (b) Let us address now the issue of equivalent classes. 
It is obvious that equivalent ribbon graphs give rise
to equivalent rank 2 stranded graphs. Now suppose 
that  $\crrG_1$ and $\crrG_2$ are equivalent rank 2 stranded
graphs arising from ribbon graphs $\cG_1$ and $\cG_2$, respectively. We must show that the corresponding 
ribbon graphs are equivalent. The strands, chords,
and points of $\crrG_1$ describe a 2-regular graph.
Now add an edge between each pair of vertices that
belong to the same pre-edge and obtain a 3-regular
graph. Introduce a proper edge coloring of this cubic
graph, according to the edge arises from 
a chord, a strand or a pre-edge. The resulting 
object is a description of $\cG_1$ as a graph encoded
map (see, for instance \cite{bonn}, page 30, for background on graph encoded maps). Since equivalence of vertex stranded graphs does not change the arising  graph encoded map,
$\cG_1$ and $\cG_2$ must be equivalent ribbon graphs.

From the above analysis,  ribbon graphs map to
rank 2 stranded graphs. Nevertheless, it must be made clear that 
the converse is not true.  To any rank 2 stranded graph,
we cannot necessarily assign a ribbon graph (for 
several reasons, one of which is that, in rank 2 stranded graphs,
 pre-edges with a single point are allowed and this does not make sense
in the context of ribbon graphs). Still from the previous analysis, 
 we can immediately establish that ribbon graphs 
with vertices with fixed coordination equals to 3 are 
rank 2 tensor graphs. In that case, the converse is also true
because, in any rank 2 tensor graph, we can
always make cyclic a rank 2 vertex by a sequence of point permutations. 
The cyclic vertex and its incident 3 stranded edges (2 in presence of a loop) are one-to-one 
mapped with a disc with three incident ribbon edges (2 in presence of a loop). 
Studying in the following sections rank 2  stranded graphs,
we will directly treat them as ribbon graphs.

\begin{figure}[h]
 \centering
     \begin{minipage}[t]{.8\textwidth}
      \centering
\includegraphics[angle=0, width=1.5cm, height=2cm]{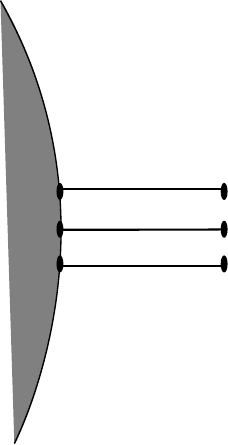}
\caption{ {\small A rank 3 sHE with its external points $a$,
$b$ and $c$.}}
\label{fig:tensflag}
\end{minipage}
\put(-153,38){$a$}
\put(-145,25){$b$}
\put(-153,13){$c$}
\end{figure}

\begin{definition}[Half-edged stranded graphs (HESGs) and stranded half-edge (sHE)]\label{def:tensflag}
$\bullet$ A rank $D$  \emph{half-edged stranded graph}  (HESG) $\crrG(\cV,\cE,\mf^0)$ (or more simply $\crrG_{\mf^0}$) consists of a set $\cV$ of disjoint 
rank $D$ stranded vertices, a set $\cE\cup\mf^0$, $\cE\cap\mf^0 = \emptyset$, of disjoint stranded edges of rank at most $D$, together with an incidence map $\phi$ from the 
set of edge points to the set of vertex points satisfying the following conditions: $\phi$ is defined on both ends of every edge in $\cE$ but only one end of each edge in $\mf^0$;
on this set of edge points $\phi$ is injective, and (as before) on each end of each stranded edge on which  $\phi$ is defined, it is a bijection to some pre-edge. Finally, $\phi$ is a surjection onto all vertex points so that every pre-edge of every stranded vertex of  has some end of some stranded edge of $\cE\cup\mf^0$ attached to it.  

\noindent$\bullet$ A rank $d$, $ 0 \le d \le D$,  stranded edge in $\mf^0$ is called rank $d$ \emph{stranded half-edge} (sHE) since they are attached at one end.  The edge points of the set of segment ends of a sHE which do not intersect the stranded vertex are called \emph{external points} of the rank $d$ sHE.  (See Figure \ref{fig:tensflag}.)

\noindent$\bullet$ Removing the set $\mf^0$ of sHEs of a rank $D$ HESG $\crrG_{\mf^0}$, and using the proper restriction of $\phi$ on the edge points from the ends of $\cE$, we  define a rank $D$ stranded graph $\crrG$ called the underlying stranded graph of $\crrG_{\mf^0}$. 
\end{definition}

\begin{figure}[h]
 \centering
     \begin{minipage}[t]{.8\textwidth}
      \centering
\includegraphics[angle=0, width=8cm, height=3cm]{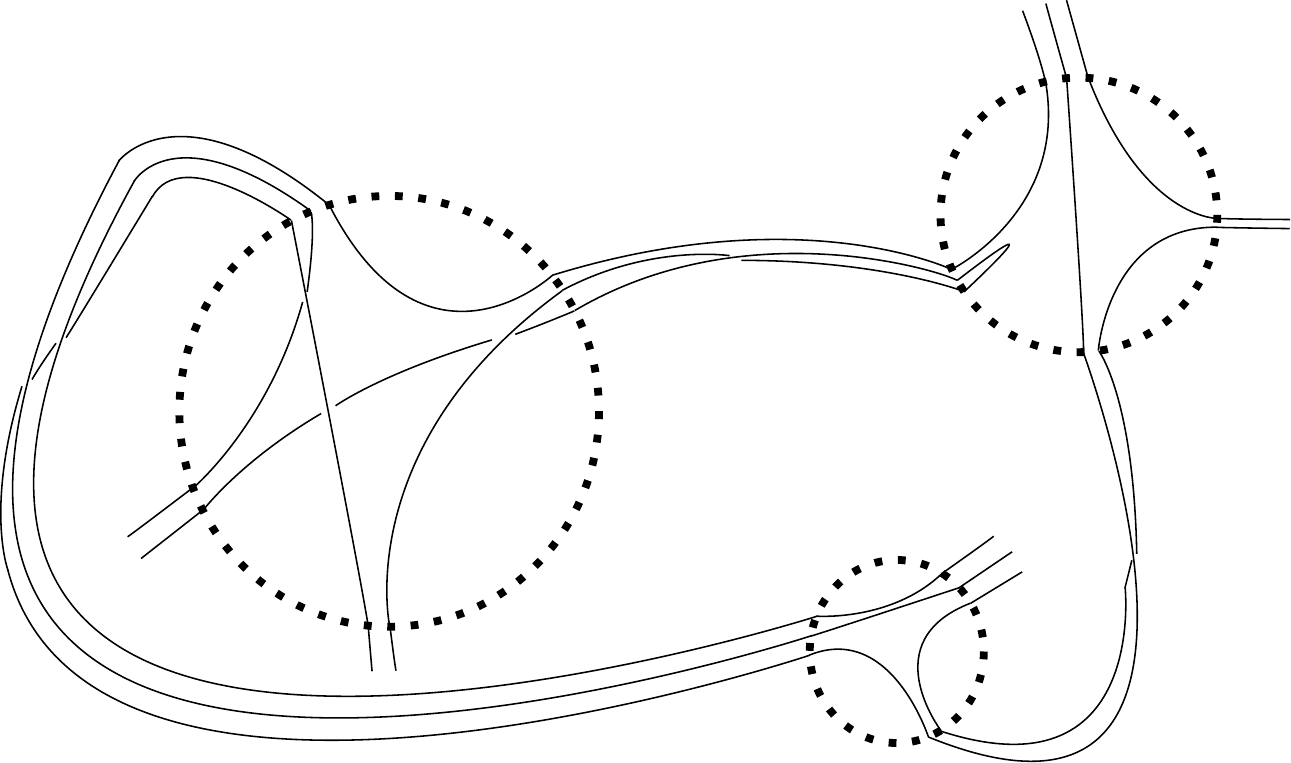}
\caption{ {\small A rank 3  HESG}.}
\label{fig:graphwfla}
\end{minipage}
\end{figure}

We represent a rank $D$ sHE by a set of $D$ (parallel) line segments (see Figure \ref{fig:tensflag}). In a HESG $\crrG_{\mf^0}$ with an empty set of sHEs, $\mf^0=\emptyset$,  every pre-edge  is  connected to a certain stranded edge. 
It is obvious that a HESG $\crrG_{\mf^0}$ may 
have trivial circles as stranded vertices if the stranded graph $\crrG$ does have ones.
An example of a HESG is given in  Figure \ref{fig:graphwfla}. 
A \emph{half-edged tensor graph} is nothing but a HESG with stranded vertices and edges satisfying also the conditions  of
a tensor graph, see Definition \ref{def:tens}. 

We can define the following edge operations on stranded graphs.

\begin{figure}[h]
 \centering
     \begin{minipage}[t]{.8\textwidth}
      \centering
\includegraphics[angle=0, width=5cm, height=1cm]{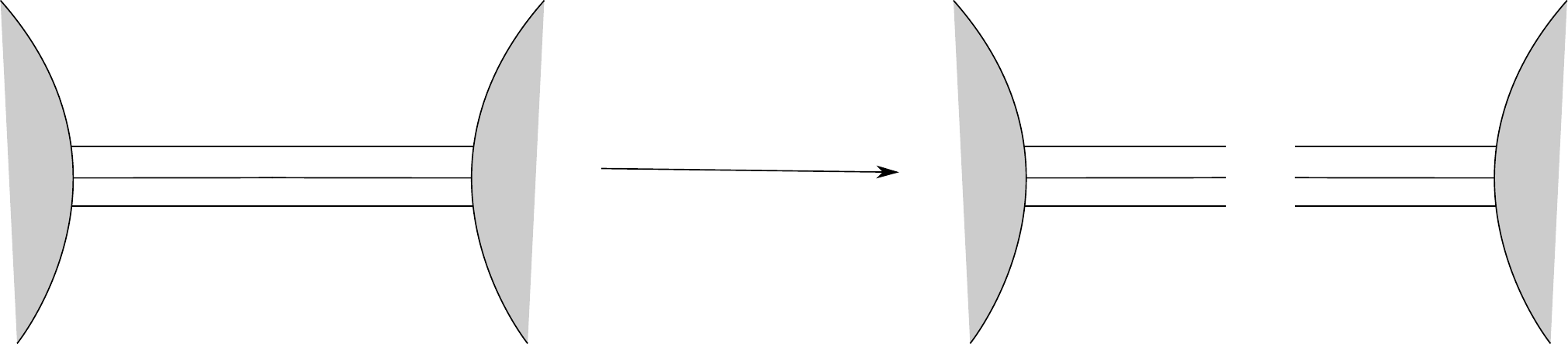}
\caption{ {\small Cutting a rank 3 stranded edge.}}
\label{fig:cuttens}
\end{minipage}
\end{figure}

\begin{definition}[Deletion and cut of a stranded edge]\label{def:cutedtens}
Let $\crrG_{\mf^0}$ be a rank D HESG  
 and $e$ one of its rank $d$ stranded edges, $1\leq d \leq D$.
The \emph{stranded graph} $\crrG_{\mf^0} - e$ is obtained from $\crrG_{\mf^0}$ by \emph{deleting} $e$.
 
The \emph{cut stranded graph} $\crrG_{\mf^0} \vee e$ is obtained from $\crrG_{\mf^0}$ by 
\emph{cutting} $e$ that is by deleting $e$ and attaching two rank $d$ 
sHEs  at the same pre-edge of the stranded vertices where $e$ was incident, one at each of the end vertices of $e$.
 (See Figure \ref{fig:cuttens}.) If $e$ is incident to a unique stranded vertex, the two  sHEs are on the same vertex.  
\end{definition}

The notion of spanning c-subgraph for HESGs follows naturally from Definition \ref{subfla}. 
The spanning c-subgraph inclusion is again denoted by $A \sset \crrtG$. 
In the formulation of \cite{Gurau:2009tz}, the deletion of a stranded edge is in fact the cut of the stranded edge in the sense precisely described above. Thus, the cut operation should be the natural operation for 
HESGs. 

Particular stranded edges as loops, bridges (defined through a cut of a stranded edge  which brings 
an additional connected component in the stranded graph),  ordinary stranded edges and terminal forms extend in the present context because a HESG has an underlying graph. 

Most of the concepts introduced in Definition \ref{faces} are valid for stranded graphs. We give their precise definition
 in the present context.  To proceed, we  realize that 
 the strands and chords in a HESG may be considered as the edges of an abstract graph in which the vertices, namely the vertex or edges points, have degree $1$
 (external points of sHEs) or $2$. In this sense, a HESG is nothing but a collection of disjoint cycles and paths.

\begin{definition}[Faces of a stranded graph]
\label{facestrand}
Consider $\crrG_{\mf^0}$ a HESG. 

$\bullet$ A \emph{face} of $\crrG_{\mf^0}$ is a
maximal alternating sequence of chords and strands, so that the end of
one is incident with the end of the next at a pre-edge point.

$\bullet$  A face is closed if  it is cycle, otherwise it is called open (the sequence
is then a path between two external points of some sHEs).

$\bullet$ By convention, a trivial circle, as a HESG, has a closed face. 

$\bullet$ The set of closed faces is denoted by $\cF_{\inter}$,
 the set open faces is denoted by $\cF_{\ext}$,
and the set of faces by $\cF= \cF_{\inter} \cup \cF_{\ext}$.

\end{definition}

Note that each open face must start at an external point of a sHE and rejoins another external point in the HESG. 
 The convention that a trivial circle has a closed face will be useful to make contact with the rank 2 case, in particular with the contraction of 
simplest trivial loops in ribbon graphs. Furthermore, these will stabilise the number of closed faces during stranded edge operations as explained in the following. 
Several notions (such as open and closed face) which are totally combinatorial in the stranded situation bear a true topological
content in the tensor graph case. For instance, a closed face in a
rank 3 tensor graph is a 2D surface in the bulk (interior) of the dual triangulation. An open face is a surface intersecting the boundary of the simplicial  complex dual to the graph.

\vspi 
\noindent{\bf Structure at the stranded edge/pre-edge connection.}
 In order to introduce the contraction of a stranded edge,
we must detail the connection
between stranded edges and pre-edges at a stranded vertex. 
Consider a stranded edge $e$ in $\crrG_{\mf^0}$,  its set $s$ of strands and the set $c$ of chords that the 
strands meet. The set $s \cup c$ can be regarded as a subgraph of 
the larger graph 
made by (the connection of) all strands and chords of $\crrG_{\mf^0}$. 

\begin{figure}[h]
 \centering
     \begin{minipage}[t]{.8\textwidth}
      \centering
\includegraphics[angle=0, width=12cm, height=4cm]{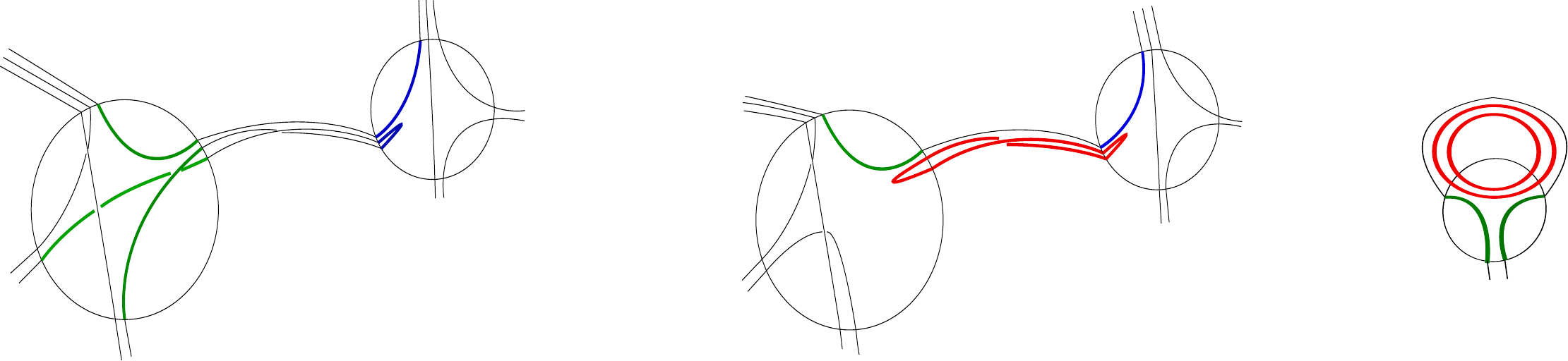}
\vspace{0.1cm}
\caption{ {\small Some rank 3 $p$--inner edges: 0--inner $e_1$ (A), 1--inner 
$e_2$ (B) and 2--inner $e_3$ (C)  edges.}}
\label{fig:neigh}
\end{minipage}
\put(-280,-10){A}
\put(-280,78){$e_1$}
\put(-120,-10){B}
\put(-140,73){$f_1$}
\put(-105,55){$f_2$}
\put(-120,75){$e_2$}
\put(-10,-10){C}
\put(-37,53){$f_1$}
\put(3,53){$f_2$}
\put(-18,87){$e_3$}
\put(-338,-10){$(v_1)$}
\put(-240,40){$(v_2)$}
\put(-170,-10){$(v_1)$}
\put(-80,40){$(v_2)$}
\put(15,30){$(v)$}
\end{figure}

\begin{definition}[Inner face, $p$--inner edge, outer strand]
\label{innerface}
Consider $e$, a rank $d$, stranded edge in a rank $D$ HESG $\crrG_{\mf^0}$,
and consider $s\cup c$ as a subgraph of 
the graph made of the strands of $e$ and the chords that they meet. 
A face is called an \emph{inner} face of $e$ if it is a cycle of $s\cup c$. 
The stranded edge  $e$ is called \emph{$p$--inner} if  $s\cup c$ contains exactly $p$ inner faces, $p\le d$. 
A strand which is not an edge of an inner face of $e$ is called \emph{outer strand} of $e$. 
\end{definition}
We use, for short, ``inner face'' or ``outer strand'' when there is no confusion about the edge
they refer to. 
See Figure \ref{fig:neigh} A, B and C for an illustration 
in the rank 3 situation. In that figure,  A describes a 0--inner edge $e_1$, B a 1--inner 
edge $e_2$, and C a  2--inner edge $e_3$.
The  inner faces are highlighted in red therein, and 
outer strands are the strands of $e_i$, $i=1,2,3$, forming paths with chords highlighted in green and blue. 
The notion of stranded edge contraction can 
be defined at this point. 
  In contracting an edge $e$, we look at the subgraph $s \cup c$ of 
associated with $e$ in the larger graph consisting of all strands
meeting all chords;  $s \cup c$ consists of paths and cycles.
 We replace each cycle by a trivial circle, and each path by a new chord joining its endpoints. 
This operation does not affect the face structure, except that some (inner) faces may be deleted in the former component containing $e$. 
However, the total number of closed and open faces remains constant: to each inner face of $e$, we associate a trivial circle which then introduces back another closed face in the HESG.

 \begin{definition}[Stranded edge   contractions]\label{def:constran}
Let $\crrG_{\mf^0}$ be a rank $D$ HESG and $e$ be a rank $d$ stranded edge with 
$s \cup c$ the subgraph associated with $e$ in 
the larger graph consisting of all strands meeting all chords (in the previous notation). 
Let $p\le d$ be a positive integer. 

If $e$ is a  $p$--inner edge,  $\crrG_{\mf^0}/e$ is the result of contracting $e$ in $\crrG_{\mf^0}$, that 
 is defined from $\crrG_{\mf^0}$ by replacing $e$ and its end vertices $v_{1}$ and $v_2$  (or its end vertex $v$, if $e$ is a loop, respectively) by 
$p$ trivial circles and a new stranded vertex $v'$. The new vertex $v'$  possesses all  pre-edges
except those connected to $e$ and all sHEs as they appear on $v_{1}$ and $v_{2}$ (respectively, on $v$), 
keeping all chords of $v_{1}$ and $v_{2}$ (respectively, of $v$) except those involved in the cycles of $s \cup c$, 
and replacing each open path of $s\cup c$ by a new chord joining its ends. If the vertex $v'$ is disconnected, then we split it. 

If there is no outer strands left after removing the $p$ inner faces
of $e$, then there is no vertex $v'$. 
\end{definition}

\begin{figure}[h]
 \centering
     \begin{minipage}[t]{.8\textwidth}
      \centering
\includegraphics[angle=0, width=12cm, height=3cm]{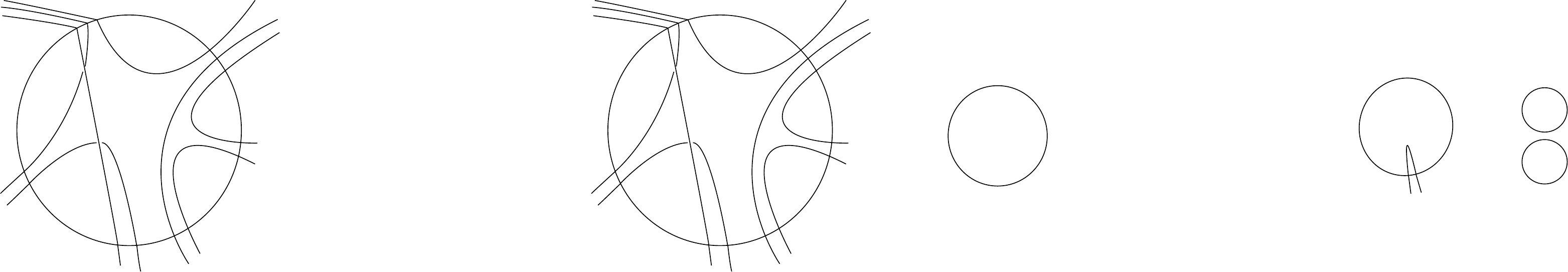}
\vspace{0.1cm}
\caption{ {\small Graphs A', B' and C' obtained after stranded edge contraction of A, B and C of Figure \ref{fig:neigh}, respectively.}}
\label{fig:constra}
\end{minipage}
\put(-310,-10){A'}
\put(-150,-10){B'}
\put(-20,-10){C'}
\end{figure}

Some examples of rank 3 stranded edge contraction are given in 
Figure \ref{fig:constra}A' and B', and an example of a rank 3 loop contraction is given in Figure \ref{fig:constra}C'.

A convenient way, though not necessary, to illustrate this type of contraction is to use a cyclic order for $v_1$ and $v_2$
contract the stranded edge  $e$ incident to these vertices and finally draw the final vertex respecting the cyclic
order of the pre-edges on $v_{1}$ and $v_2$. This choice leads to the representation given in Figure \ref{fig:constra}.

One may directly check that contracting a non loop in a ribbon graph
can be  seen as a rank 2 stranded edge contraction. 
 One may also check that a trivial loop contraction for a ribbon graph coincides with the notion of contraction of a  loop  in the sense of Definition \ref{def:constran},  when the ribbon graph is viewed as a rank 2  stranded graph. Indeed, in a ribbon graph, a trivial  positive  loop can be a 0--,   1-- or 2--inner loop.
 A trivial  negative loop can be a 1-- or 0--inner loop. 
Contracting these trivial loops in a ribbon graph is equivalent
of what is described in Definition \ref{def:constran}. Discussing loops, we  only focus on trivial ones in the following. 

 The following proposition is  straightforward. 
\begin{proposition}\label{prop:rank}
Let $\crrG_{\mf^0}$ be a rank $D$ HESG and $e$ be one of its stranded edge. Then
 $\crrG_{\mf^0}/e$ obtained by contraction of $e$  is a rank $D'$ HESG, with $D'\le D$. 
\end{proposition} 
At this point, one may wonder if the additional circles obtained after  stranded edge contraction are not inessential features of the HESGs. In fact, they are very useful for the preservation of number of faces  during the edge contraction of the stranded 
graph,  in analogy with the case of ribbon graphs. 
In another context related to the topology associated with the stranded graphs, another type of contraction (without trivial circles involved) can be defined  on stranded graphs, see for instance, in \cite{BenGeloun:2011rc}. 

\vspi

\noindent{\bf Colored tensor graphs.} 
  Apart from the stranded structure, the second important feature that we need is the coloring.

\begin{definition}[Colored  tensor graph \cite{Gurau:2009tw,Gurau:2011xp}]\label{def:coltens}
A rank $D\geq 1$ properly colored tensor graph $\crrG$ is a
rank $D$ tensor graph such that the underlying graph 
is bi-partite and the stranded edges of $\crrG$ are colored with 
$D+1$ colors such that no two stranded edges that meet at a stranded vertex
have the same color. The coloring of stranded edges determines the coloring
of pre-edges that they meet. 
The  chord coloring is the following: in each stranded vertex, 
each vertex point has an ordered pair of colors $(i, j)$: 
$i$ is the color of the pre-edge it is in, and $j$ is the color of the pre-edge containing the other end of its chord. 
(So for each $i$, $j$ with $j  \neq i$ there is one point with color $(i, j)$).
The unordered pair $\{i,j\}$ determined by the coloring of the end points of 
some chord yields the chord coloring. 
In a stranded edge of color $i$, the end points of each strand have
the same color pair $(i, j)$; the unordered pair $\{i,j\}$ yields the strand coloring. 
\end{definition}

The color restriction introduced in Definition \ref{def:coltens}
allows us to control the type of graphs 
generated by gluings of stranded vertices and edges. Later, 
the coloring $\{i,j\}$ of strands or chords will be called \emph{bi-coloring} and we will work with an unordered pair  $(ij)$, called
color pair, for strands or chords. Observe that strands and chords that meet have necessarily the same color pair.  For simplicity, we will 
consider that a colored tensor graph does not
have trivial circles from the beginning. The general case should not be hard
to recover with a given assignment of bi-coloring (a prescription by default) 
of the closed faces of these trivial vertices. 

\begin{figure}[h]
 \centering
     \begin{minipage}[t]{.8\textwidth}
      \centering
\includegraphics[angle=0, width=5cm, height=3cm]{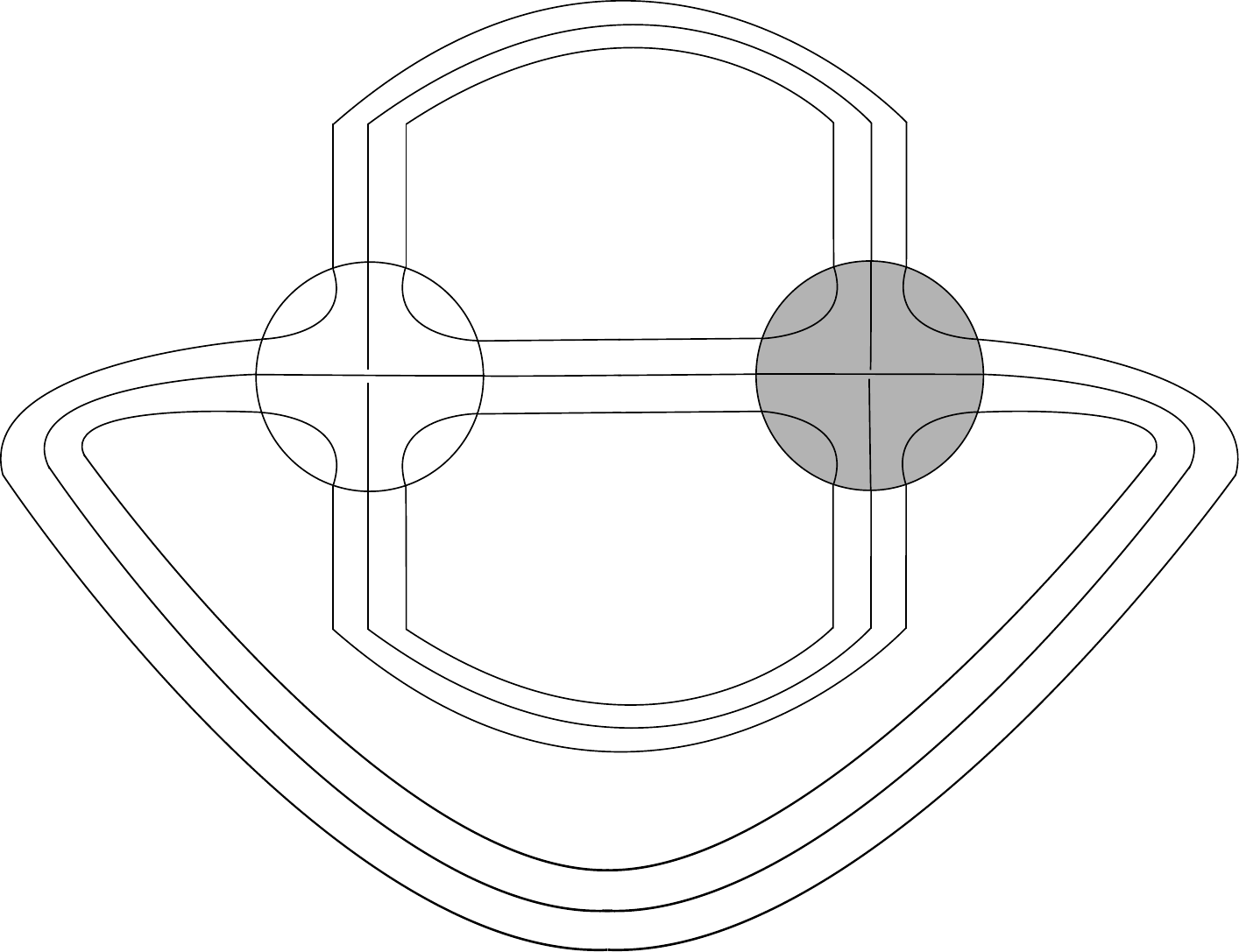}
\vspace{0.1cm}
\caption{ {\small A rank 3 colored tensor graph.}}
\label{fig:colgraph}
\end{minipage}
\put(-175,-10){1}
\put(-175,25){2}
\put(-175,58){3}
\put(-175,90){0}
\end{figure}
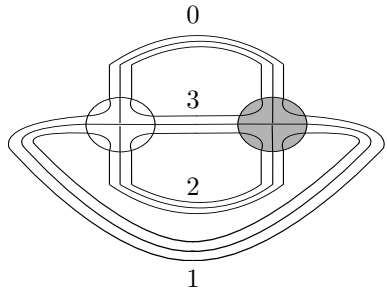

We denote a rank $D$ colored tensor  graph by $\crrG(\cV,\cE)$
as usual.   As an illustration, a rank 3 colored tensor graph is pictured
in Figure \ref{fig:colgraph}. Each vertex (with vertex graph $K_{4}$) 
 is the dual of a tetrahedron and an edge represents a triangle endowed with a color $i \in \{0,1,2,3\}$. The (underlying) graph is 
also bi-partite (white and shaded vertices).

\begin{figure}[h]
 \centering
     \begin{minipage}[t]{.8\textwidth}
      \centering
\includegraphics[angle=0, width=2.5cm, height=1.5cm]{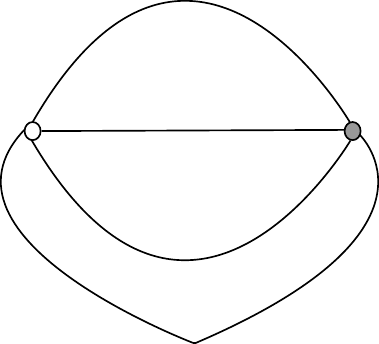}
\vspace{0.1cm}
\caption{ {\small Underlying graph of the graph of Figure \ref{fig:colgraph}.}}
\label{fig:compak}
\end{minipage}
\put(-175,-10){1}
\put(-175,12){2}
\put(-175,28){3}
\put(-175,45){0}
\end{figure}

 It turns out that, in any rank,  the stranded structure of a colored tensor graph $\crrG$ can be  captured at the level 
of its underlying bipartite colored graph.  For instance, the graph of  Figure \ref{fig:colgraph}
can be also drawn in the collapsed form of Figure \ref{fig:compak}. 
It is worth highlighting that 
one representation of the graph unambiguously determines  the other. 
This is not the case for a generic underlying graph of a 
stranded or tensor graph without colors. In order to render this expansion explicit in the colored case,   one must use the bi-coloring of chords and strands.
This is the way to proceed:  consider a $ D+1\geq 2$ edge properly colored graph $G$, then replace each vertex by 
a rank $D$ stranded vertex with coordination 
$D+1$, with vertex graph $K_{{D+1}}$, and pre-edges with fixed cardinality $D$. Now replace each colored
edge in $G$ by a rank $ D+1$ colored edge with the same color and respecting the same incidence relation as in $G$ and
form $\crrG$ a rank $D$ tensor graph.  Assign a color to 
each pre-edge according to the coloring of the edges incident to these. 
A bi-coloring for chords and strands can be deduced from that point. 
This is the rank $ D$ colored tensor graph.

There are certainly 
more data worthwhile to be discussed in such a colored graph.

\begin{definition}[$p$-bubbles \cite{Gurau:2009tw}]
\label{def:pbub}
Let $\crrG$ be a rank $D$ colored  tensor graph.

$\bullet$  A \emph{0-bubble} is a vertex of $\crrG$.

$\bullet$  A \emph{1-bubble} is an edge of $\crrG$.

$\bullet$  For all $p \geq 2$, a \emph{$p$-bubble} of $\crrG$ with colors $i_1 < i_2 < \dots<i_p $, $p\leq D$,
and  $i_k\in \{0, \dots, D\}$ is a  connected rank $p-1$ colored tensor graph the  underlying  graph of which is a   maximal connected subgraph of the  underlying  graph of $\crrG$ made of edges of colors $\{i_1, \dots, i_p\}$. 
 \end{definition}
For $p \geq 2$, a $p$-bubble can be mapped to a rank $p-1$ colored tensor graph because we can associate
to each $p$ colored graph a rank $p-1$ colored tensor graph
according to the procedure explained above.  

 To obtain the set of $p$-bubbles is actually easy in  a colored tensor  
graph $\crrG$. Consider the underlying graph of $\crrG$. 
Take a subset $C$ of cardinality $p$ of the set of colors $\{0, 1, \dots, D\}$, 
delete all edges of colors $\{0, 1, \dots, D\} \setminus C$ in the underlying
graph. Each connnected component of the resulting graph 
is a $p$-regular edge colored graph that determines uniquely 
 a rank $p-1$ colored tensor graph, a $p$-bubble.
We can alternatively perform the deletion at the level of the stranded 
graph: delete all stranded edges of color  $\{0, 1, \dots, D\} \setminus C$,
 all vertex points, and the chords incident to those, if their color pair $(i,j)$ 
 involves $i$ or $j$ in  $\{0, 1, \dots, D\} \setminus C$. Each remaining connected
 component is a rank $p-1$ colored tensor graph that defines a $p$-bubble. 

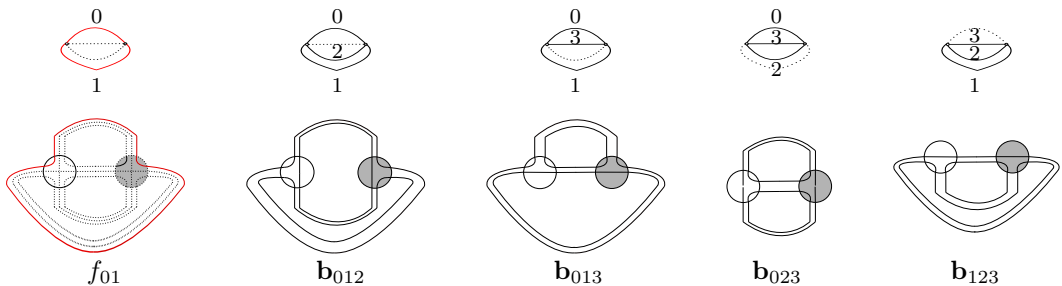
\begin{figure}[h]
 \centering
     \begin{minipage}[t]{.8\textwidth}
      \centering
\includegraphics[angle=0, width=14cm, height=3cm]{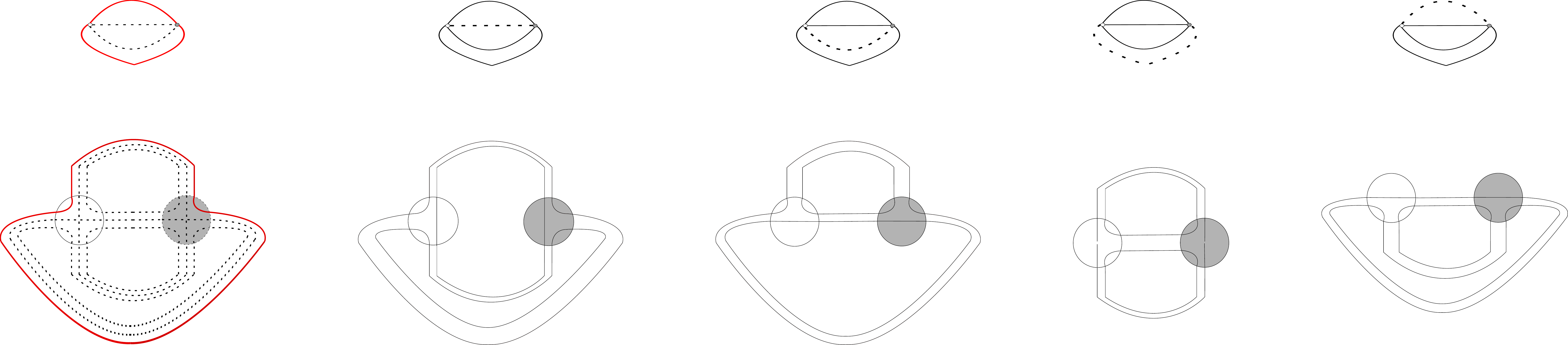}
\vspace{0.1cm}
\caption{ {\small The face  $f_{01}$ (in red) and bubbles of the graph 
of Figure \ref{fig:colgraph}.}}
\label{fig:facebub}
\end{minipage}
\put(-312,-10){ $f_{01}$}
\put(-222,-10){$\bee_{012}$}
\put(-132,-10){$\bee_{013}$}
\put(-57,-10){$\bee_{023}$}
\put(18,-10){$\bee_{123}$}
\put(-307,61){\footnotesize{1}}
\put(-307,87){\footnotesize{0}}
\put(-216,74){\footnotesize{2}}
\put(-216,61){\footnotesize{1}}
\put(-216,87){\footnotesize{0}}
\put(-126,79){\footnotesize{3}}
\put(-126,61){\footnotesize{1}}
\put(-126,87){\footnotesize{0}}
\put(-50,79){\footnotesize{3}}
\put(-50,67){\footnotesize{2}}
\put(-50,87){\footnotesize{0}}
\put(25,80){\footnotesize{3}}
\put(25,73){\footnotesize{2}}
\put(25,61){\footnotesize{1}}
\end{figure}

Restricting to rank $D=3$ colored tensor graphs,  there are 
two types of $p$-bubbles that we shall study in the following:
 2-bubbles coincide with the faces (Definition \ref{facestrand})  of the colored
tensor graph.  Thus, there are 
alternative ways to observe faces of a colored tensor graph. 
Faces can be read from the underlying colored graph as cycles with alternating (edge) colors. 
On the other hand, faces are cycles made with chords and strands with 
a given color pair (see the face $f_{01}$ (in red) in  Figure \ref{fig:facebub}). 
 The alternating colors are precisely the color pair of the strands and chords (in the sense of Definition \ref{def:coltens}) forming
the face. They are also rank 1  colored tensor graphs.
  3-bubbles (or simply bubbles in $D=3$) are in   one-to-one 
correspondence with maximal connected components of the underlying graph which have three colors
 (see Figure \ref{fig:facebub}). These are rank 2 colored tensor
graphs and also colored ribbon graphs with 3-valent vertices.

As mentioned earlier in this section, a 3-bubble, which is a rank 2 tensor graph, can be seen 
as a ribbon graph. In $D=3$, the vertices of a bubble are 3-valent vertices obtained by decomposing  the vertex of the graph
in the way of Figure \ref{fig:facebub}. The edges of a 3-bubble are colored 
ribbon edges generated by the decomposition of the colored stranded edges of the rank $D$ colored graph. Like ordinary 
ribbon graphs, 3-bubbles have faces as well. These faces
are endowed with a pair of colors like the faces of the initial graph. Thus, a bubble is simply a rank 2 colored tensor graph
lying inside the rank 3 colored tensor graph. 
The set of 3-bubbles is denoted by $\cB_3$ and $|\cB_3| = B_3$. The index 3  is omitted when 
discussing $D=3$.

\medskip 

\noindent{\bf Rank $D$ half-edged colored tensor graphs.}
Rank $D$ sHEs can be considered as well on rank $D$ colored tensor
graphs, provided these sHEs possess a color and their gluing respects the graph coloring at each vertex. 
 For any  rank $D$ colored tensor graph with sHEs, we demand
that to each stranded edge and sHE, one  assigns a color $i \in \{0,1,\dots,D\}$ such that no two  stranded edges or sHEs meeting at a vertex share the same color.

The cut of a stranded edge can be understood in the same sense of Definition \ref{def:cutedtens} for colored tensor graphs. 
After cutting a colored stranded edge,  the resulting sHEs possess the same color of that  edge.

\begin{definition}[Half-edged colored tensor graph (HEcTG)]\label{def:hectg}
A rank $D$ HEcTG is a rank $D$ half-edged tensor graph
such that its underlying graph is a bipartite HEG,  its stranded edges and  sHEs are colored with $D+1$ 
colors such that no two stranded edges or half-edges meeting
at a vertex have  the same color. The coloring of pre-edges, vertex points,  
and chords of stranded vertices, and the coloring of strands of stranded edges and sHEs are described by Definition \ref{def:coltens}. 
 \end{definition}

\begin{figure}[h]
 \centering
     \begin{minipage}[t]{.8\textwidth}
      \centering
\includegraphics[angle=0, width=14cm, height=2.5cm]{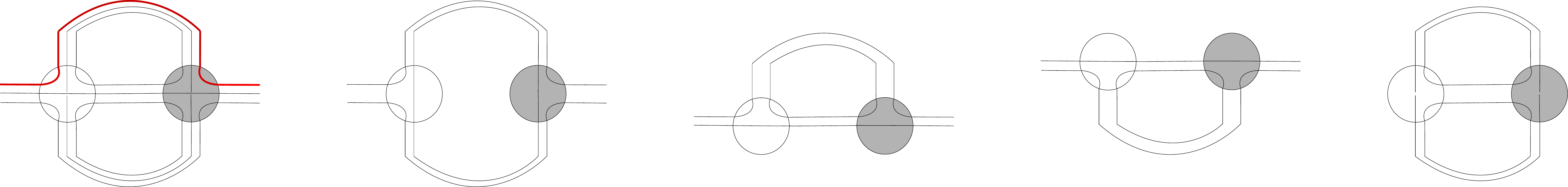}
\vspace{0.1cm}
\caption{ {\small A rank 3 HEcTG  
and its bubbles; $f_{01}$ (highlighted in red) is an open face;
$\bee_{012}$, $\bee_{013}$ and $\bee_{123}$ are open bubbles and $\bee_{023}$ is closed.}}
\label{fig:tenfla}
\end{minipage}
\put(-350,33){1}
\put(-268,33){1} 
\put(-309,-10){2}
\put(-309,45){3}
\put(-309,75){0}
\put(-225,-10){$\bee_{012}$}
\put(-139,-10){$\bee_{013}$}
\put(-50,-10){$\bee_{123}$}
\put(28,-10){$\bee_{023}$}
\end{figure}
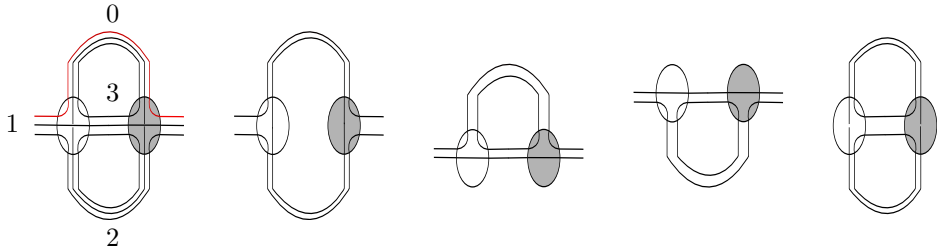

An example of a rank $D=3$ HEcTG is given  on the left in Figure \ref{fig:tenfla} 
(most left, edges and sHEs are colored). 
Spanning c-subgraphs of a HEcTG follow from the similar notion for
HESG and from Definition \ref{subfla}. The only point to be  added is the coloring. 

The following is straightforward:

\begin{proposition}\label{prop:subcol}
Spanning c-subgraphs of a rank $D$  HEcTG are rank $D$  HEcTG. 
\end{proposition}

Let us investigate some properties  of a HEcTG $\crrG_{\mf^0}$.  The notion of bridge is standard. 
 By definition of a proper edge coloring, HEcTGs cannot contain loops. 

As in the case of HESGs, by introducing sHEs on colored tensor
graphs, we distinguish two types of faces. 
Using Definition \ref{facestrand}, faces can be open or
closed connected components if they pass through external points of  sHEs
or not (let us recall that these are also cycles and paths in the abstract graph formed by all strands and chords of the underlying HESG). 
An open face with color pair (01) is highlighted in red in Figure \ref{fig:tenfla}. 
The sets of closed and open faces are denoted by $\cF_{\inter}$
and $\cF_{\ext}$, respectively. 
Hence, for a rank $D$ HEcTG, the set  $\cF$ of faces is the disjoint union $\cF_{\inter} \cup \cF_{\ext}$.

To define $p$-bubbles for HEcTG, we simply replace
``edge'' by ``edge or sHE'',  and ``underlying graph''  by ``underlying HEG''
 in Definition \ref{def:pbub}.
We say that a $p$-bubble in some HEcTG $\crrG_{\mf^0}$ is open if its underlying graph contains half-edges.
Hence, a $p$-bubble is open if it contains open faces, otherwise
it is  closed. Open and closed bubbles for a rank 3 HEcTG have been illustrated in Figure \ref{fig:tenfla}.  
The sets of closed and open bubbles are denoted by $\cB_{\inter}$ and $\cB_{\ext}$, respectively.

The following notion will play a crucial role in our following construction upon 
HEcTG. That notion turns out to be defined at the generic level 
of HESG. 

\begin{definition}[Boundary of an HESG]
\label{def:bgph}
The \emph{boundary graph} $\brrG(\bV,\bE)$  of a rank $D$  HESG $\crrG(\cV,\cE,\mf^0)$ 
is a graph with vertex set $\bV$  in one-to-one  correspondence with $\mf^0$, with edge set $\bE$ in 
one-to-one correspondence with $ \cF_{\ext}$. 
Consider an edge $e \in \bE$, its corresponding open face $f_e\in 
\cF_{\ext}$, a vertex $v \in \bV$, and its corresponding sHE $h_v$. 
Then, $e$ is incident to $v$ if and only if $f_e$ has one end-point
in $h_v$.  The boundary graph of a rank $D$ HESG with 
$\mf^0=\emptyset$ is empty.  
\end{definition}

\begin{figure}[h]
 \centering
     \begin{minipage}[t]{.8\textwidth}
      \centering
\includegraphics[angle=0, width=4cm, height=2.5cm]{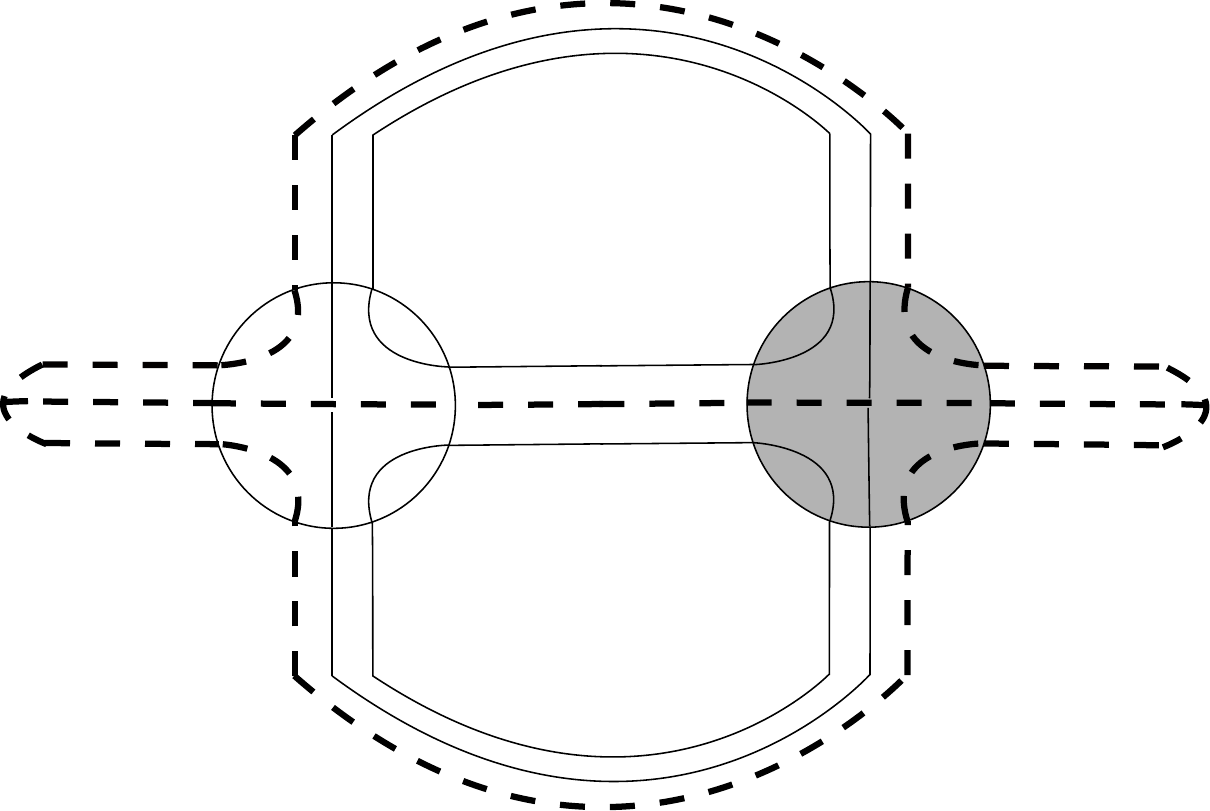}
\vspace{0.1cm}
\caption{ {\small Boundary graph (in dashed lines) of the graph of Figure \ref{fig:tenfla}.}}
\label{fig:bound3}
\end{minipage}
\put(-240,33){1}
\put(-108,33){1}
\put(-174,-11){2}
\put(-174,42){3}
\put(-174,75){0}
\end{figure}
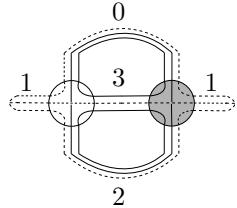

There is a procedure for drawing the boundary graph of a 
rank $D$ HESG described by Gurau as \emph{pinching} in \cite{Gurau:2009tz}. Insert a vertex at each sHE and 
make them incident to the strands of the sHE (they become incident to the open faces).   In the case of a rank $D$ half-edged tensor graph, the vertices of the boundary 
graph must be $D$-regular vertices. 
 In doing so for the rank 3 HEcTG of Figure \ref{fig:tenfla} (most left), 
we associate its boundary  depicted in Figure \ref{fig:bound3}.
 In fact,  the  coloring of  a HEcTG  entails a new coloring on its boundary graph.
Both types of colorings will  allow us to enumerate the different constituents (or to find
bounds on their number)  of these generalized graphs. 
\begin{definition}[Ve-colored graphs]\label{def:2colbi}
  A  $(D+1)$ \emph{ve-colored graph} is a graph with a vertex  coloring 
with colors from $\{0,\dots,D\}$ and a proper edge coloring 
such that each edge is assigned an unordered pair $(ab)$ of colors, $a,b$ in $\{0,\dots,D\}$, $a\neq b$, 
 and such that  the end-vertices of an edge with color (ab) must have color $a$ or $b$. The two ends may or may not have distinct colors.
\end{definition}

\begin{figure}[h]
 \centering
     \begin{minipage}[t]{.8\textwidth}
      \centering
\includegraphics[angle=0, width=7cm, height=2.5cm]{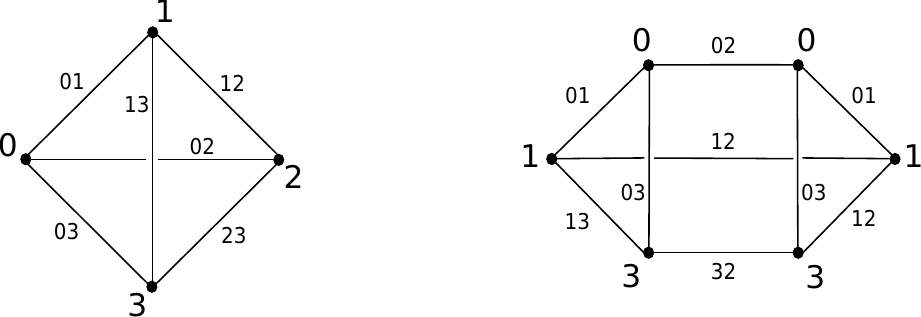}
\vspace{0.1cm}--
\caption{ {\small Examples of 4 ve-colored graphs.}}
\label{fig:vecolor}
\end{minipage}
\end{figure}

Some 4 ve-colored graphs are drawn in Figure \ref{fig:vecolor}.

  One notices that a $(D+1)$ ve-colored graph is a $D(D+1)/2$ edge colored
graph, with edge colors chosen in pairs $ab$, $a<b$, in 
$\{01,02, \dots, 0D, 12,\dots, (D-1)D\}$. Because of the proper edge coloring,
ve-colored graphs have no loops.

\begin{definition}[Ve-colored  tensor graph]\label{def:2coltens}
A rank $D\geq 1$ \emph{ve-colored tensor graph} $\crrG$ is a
rank $D$ tensor graph in which 

1) every stranded vertex has a color from 
$\{0, 1, \cdots , D+1\}$;  within a stranded vertex of color $i$, 
the $D+1$ pre-edges have each a distinct color pair $ij$, where $0 \leq j \leq D+1$, $j \ne i$ 
(vertex coloring); 

2) every stranded edge has a color pair $i j$, where $0 \leq i, j \leq D+1$ are distinct,  such that no two stranded edges that meet at a stranded vertex
have the same color pair (proper edge coloring). 

3) If a stranded vertex $v$ is an end of a stranded edge $e$ of color pair $i j$, then $v$ must have color $i$ or color $j$.  The colorings of the pre-edges and  stranded  edges are consistent in the natural sense: if a  pre-edge and a stranded  edge meet, they must have the same color pair. 

4) Every strand and every chord has a color triple $i j k$, where $0 \leq i, j, k \leq D$ are pairwise distinct, with the following properties:

(i) within a stranded vertex $v$ of color $i$, a chord that is incident 
to two vertex points of two pre-edges of distinct color pair $ij$ and $ik$
has color $ijk$,  $k \notin \{i, j\}$;

(ii) within a stranded edge $e$ with color $i j$, there are exactly $D$ strands, with  $D$ color triples $i j k$, $k \notin \{i, j\}$;

(iii) if a strand and a chord  meet at a vertex point they must have
the same color triple $ijk$. 

\end{definition}

Note that a ve-colored tensor graph is not a colored tensor graph. For instance,  
it may be not bi-partite.

\begin{proposition}[Underlying graph of a ve-colored  tensor graph]\label{def:2coltens0}
The underlying graph of a rank $D$ ve-colored tensor graph  is
$(D+2)$ ve-colored and $D+1$ regular.
\end{proposition}

\proof
The proof is fairly straightforward. 
Let $\crrG$ be a rank $D$ ve-colored tensor graph and $G$ be its corresponding underlying graph. The graph $G$ is $D+1$ regular 
just like $\crrG$ is. $G$ inherits from $\crrG$ the vertex coloring with colors from $\{0, \cdots, D+1\}$. 
Since each edge $e$ of $\cG$ is in one-to-one correspondence with a stranded edge $e'$ of $\crrG$, we assign to $e$
the color pair $i j$ of $e'$
where $0 \leq i, j \leq D+1$ are distinct. This also implies that no two edges in $G$ that meet at a vertex have the same color pair. The last condition stating that the end-vertices of an edge $e$ of $G$ with color pair $ij$ must have color $i$ or $j$ is also fulfilled: the third point in Definition \ref{def:2coltens} holds
for $\crrG$ and  this reflects on $G$. 

\qed

\begin{proposition}
\label{ve-colorisstranded}
A $(D+2)$  ve-colored graph which is $(D+1)$ regular uniquely  determines 
 a rank $D$ ve-colored tensor graph of which it is the underlying
 graph. 
 \end{proposition}
 \proof 
Consider  $G(V,E)$ a $(D+2)$ ve-colored graph,
with colors $\{0,1,\dots, D+1\}$,
that is $(D+1)$ regular. 
   To each vertex $v$ of $V$, 
we assign  a rank $D$ stranded vertex $v'$
with coordination $D+1$ with vertex graph $K_{{D+1}}$, and pre-edges with fixed cardinality $D$. 
The vertex coloring is carried along, from the vertices
of $V$ to the stranded vertices.
In a given stranded vertex $v'$ of color $i$, 
all $D+1$ pre-edges can be given 
a distinct color pair $ij$ in a way compatible with the first point of
Definition \ref{def:2coltens}. 
This fills the point 1). 
The point $4(i)$ of this 
definition can be
also implemented without ambiguity
given the pre-edge bi-coloring. 
Now replace every  
edge $e$ in $G$ of color pair $ij$ by a rank $D$  colored edge $e'$ with the same color pair $ij$.
We keep the incidence relation between all
corresponding stranded edges and vertices. 
Hence the proper edge coloring is carried along. This fulfills condition 2). Note that we have not yet specified the chord-strand connection. 
Let us address 3). 
Since $e$ of color pair $ij$ in $G$ is incident 
to vertices of color $i$ or $j$, 
the same should apply to 
the corresponding stranded edge $e'$: it is incident to stranded vertices of color $i$ or $j$.
A stranded vertex $v'$
of color $i$ has exactly $D+1$ pre-edges with distinct color pair $ij$, on which must 
end a stranded edge of the same pair $ij$. 
There is no ambiguity and therefore 3) holds. It only remains  to specify the connection between strands and chords at the level of a single stranded edge and vertex connection. 
Take a rank $D$ stranded edge $e'$
with color pair $ij$
incident to a stranded vertex $v'$ of color $i$. 
Construct the coloring of the $D$ strands of $e'$ so that 4$(ii)$ is obeyed.  
We immediately see that there is a one-to-one correspondence between the color triples $ijk$ of the strands of $e'$ 
and 
the color triples $ijk$ of the chords 
at the pre-edge of $v'$ incident to 
$e'$.  There
is a single choice for
the connection strand-chords so that
4$(iii)$ holds. 
Finally, it is obvious that the underlying graph of the constructed stranded graph is
our initial graph $G$.  
 
 \qed

\begin{proposition}[Ve-coloring of $\brrG$]
\label{boundaryisve}
The boundary graph $\brrG(\bV,\bE)$ of a rank $D$ HEcTG $\crrG(\cV,\cE,\mf^0)$ 
is $(D+1)$ ve-colored and determines in a unique way  a rank $(D-1)$ ve-colored tensor graph.
\end{proposition}

\proof

The vertex coloring of $\brrG$ is inherited from the coloring of the sHEs of $\crrG_{\mf^0}$, 
that is to each vertex of $\brrG$, one assigns the same color of its
corresponding sHE.
The edge bi-coloring in $\brrG$ coincides with the open face bi-coloring. 
At a sHE, the color pairs of all open faces never coincide. This makes the graph $\brrG$  $(D+1)$ ve-colored. Furthermore the graph $\brrG$ is $D$-regular as $\crrG_{\mf^0}$ is of rank $D$ and so are its sHEs.
The existence of a unique $(D-1)$ ve-colored tensor graph associated with $\brrG$ follows from 
Proposition \ref{ve-colorisstranded}. 

\qed 

With
Proposition \ref{boundaryisve}, we can identify the boundary graph $\brrG$ with the rank $D-1$ ve-colored stranded graph   associated with it.

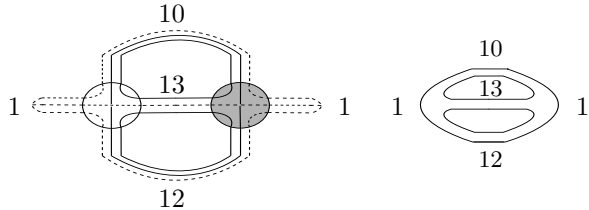
\begin{figure}[h]
 \centering
     \begin{minipage}[t]{.8\textwidth}
      \centering
\includegraphics[angle=0, width=7cm, height=2cm]{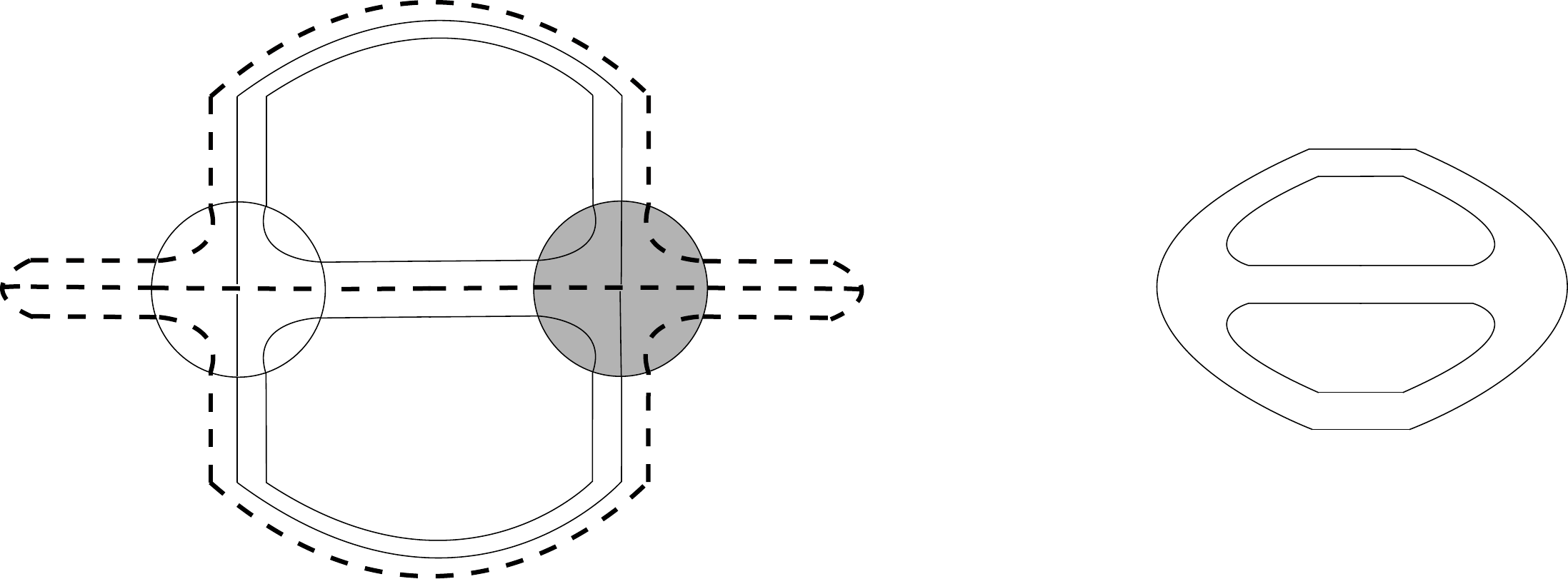}
\vspace{0.1cm}
\caption{ {\small Rank 2 ve-colored stranded structure of the boundary
graph of Figure \ref{fig:tenfla}.}}
\label{fig:boundtens3}
\end{minipage}
\put(-155,25){1}
\put(-280,25){1}
\put(-223,-10){12}
\put(-223,33){13}
\put(-223,60){10}
\put(-135,25){1}
\put(-65,25){1}
\put(-102,5){\small{12}}
\put(-102,32){\small{13}}
\put(-102,47){\small{10}}
\end{figure}

 A boundary graph does not  have sHEs hence all of its faces are closed. 
We take an example in rank $D=3$, see Figure \ref{fig:boundtens3}.
Each vertex of $\brrG$ is colored and is 3 valent. Each edge of $\brrG$ is mapped
to a bi-colored ribbon,  and each face is closed.

\begin{remark}\label{rem:1} 
Let us discuss a feature  in rank $D=3$ that will be important for the next computations. 
Consider a HEcTG $\crrG_{\mf^0}$, its boundary graph $\brrG$, its set of open bubbles, each of them viewed as a HERG (rank 2), 
and the set of their boundary graph. As a rank 1 graph, the boundary graph of any open bubble is only made of cycles that
are still called faces.  A  face of $\brrG$  is a cycle made of strands and chords coming from open faces of $\crrG_{\mf^0}$.
These strands and chords, respectively, necessarily belong to stranded edges and vertices, respectively, of open bubbles in $\crrG_{\mf^0}$ (remember that we decompose $\crrG_{\mf^0}$ in maximal connected 3 colored graphs, the bubbles). 
By taking the boundary graph of such open bubbles, we produce cycles which are in one-to-one correspondence
with the set of (closed) faces of $\brrG$. 
\end{remark}

\begin{figure}[h]
 \centering
     \begin{minipage}[t]{.9\textwidth}
      \centering
\includegraphics[angle=0, width=15cm, height=3cm]{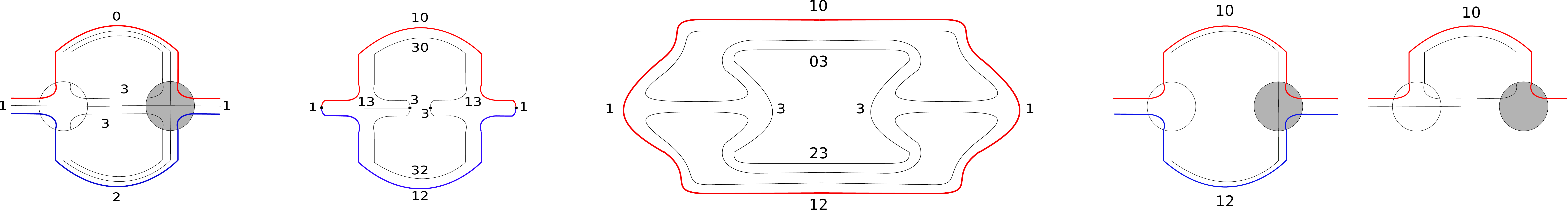}
\vspace{0.1cm}
\caption{ {\small A HEcTG, its boundary, and its expansion 
as a ribbon graph. The face $f_{012}$ is made from the open faces $f_{01}$  and $f_{12}$.
$f_{01}$ appears in 2 open bubbles, $\bee_{012}$ and $\bee_{013}$. 
$f_{012}$ corresponds to a cycle in  $\partial\bee_{012}$.}}
\label{fig:boundC}
\end{minipage}
\put(-355,-10){$\crrG_{\mf^0}$}
\put(-275,-10){$\brrG$}
\put(-165,-10){$\brrG$}
\put(-55,-10){$\bee_{012}$}
\put(12,-10){$\bee_{013}$}
\put(-217,72){$f_{012}$}
\put(-375,71){$f_{01}$}
\put(-375,10){$f_{12}$}
\end{figure}

As an illustration of the above remark, consider the HEcTG $\crrG_{\mf^0}$
given in Figure \ref{fig:tenfla}, with 
boundary $\brrG$ drawn in Figure \ref{fig:bound3}. 
  Let us pick the face $f_{102}$ of $\brrG$ that is a face (obviously closed)
made of the strands with color pairs (12) and (10), and chords of the vertex of color 1. 
Consider the unique open bubble $\bee_{012}$ of $\crrG_{\mf^0}$ as a HERG.   
The open faces $f'_{10}$ and $f'_{12}$ of $\bee_{012}$  
 will generate a cycle in its boundary graph $\partial \bee_{012}$ that is corresponding to 
$f_{102}$. Furthermore, such a closed face cannot
belong to another bubble by color exclusion and the fact that each strand 
of a given stranded edge is exactly used twice to make two stranded edges of two different bubbles
of $\crrG_{\mf^0}$ (for instance the strand of color (01) making $f'_{01}$ is used once  in $\bee_{012}$ and once in $\bee_{013}$). 
We illustrate this fact once again in a slightly more involved situation given by Figure \ref{fig:boundC}. 

\subsection{W-colored stranded graphs}
\label{subset:wgraphs}

It has been discussed earlier that the contraction of a stranded edge in a colored tensor graph yields another
type of graph for which neither proper edge coloring, nor the tensor axioms (Definition \ref{def:tens})  
apply. In order to circumvent such an odd feature of tensor graphs
and thereby find polynomial invariants on these structure, 
some proposals have been made to redefine the notion 
of contraction of an edge or redefine subgraphs
for which the contraction applies \cite{Gurau:2009tz,Tanasa:2010me}. 
In the following, we  use a different scheme. 

Using  Definition \ref{def:constran},  we can now 
contract any rank $D$ edge provided the fact that we are working in the extended
framework of stranded graphs. This definition therefore
applies to a  HEcTG. From now on, edge contraction means always stranded edge contraction
in the sense  Definition \ref{def:constran}.

\begin{definition}[Rank $D$ w-colored graph]\label{wdcolo}
A \emph{rank $D$ weakly colored}  or \emph{w-colored graph} is a rank $D$ HESG obtained by zero or more, successive stranded edge contractions of some rank $D$ HEcTG. 
\end{definition}

\begin{figure}[h]
 \centering
     \begin{minipage}[t]{.8\textwidth}
      \centering
\includegraphics[angle=0, width=7cm, height=2cm]{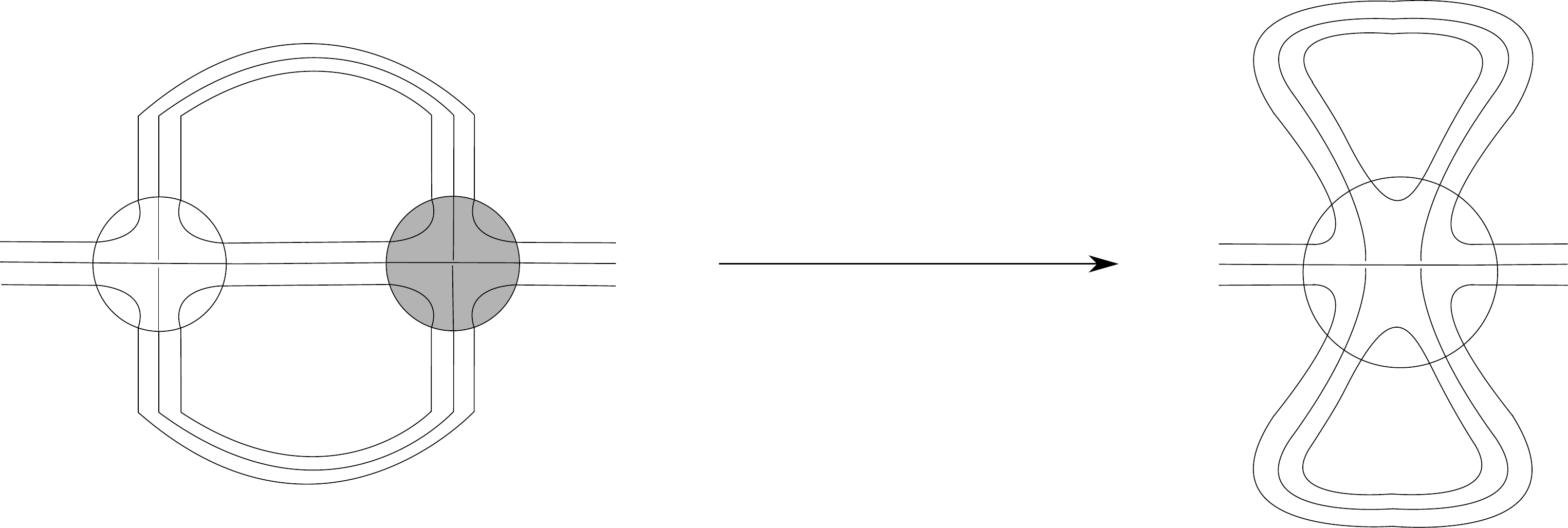}
\vspace{0.1cm}
\caption{ {\small Contraction of an edge in a rank 3 HEcTG.}}
\label{fig:contens}
\end{minipage}
\put(-188,25){1}
\put(-280,25){1}
\put(-233,-5){2}
\put(-233,33){3}
\put(-233,55){0}
\put(-125,25){1}
\put(-68,25){1}
\put(-95,-10){2}
\put(-95,60){0}
\end{figure}

Few remarks must be made at this point.  Any HEcTG is a w-colored graph of the same rank 
when no stranded edge contractions have been performed. For  any colored graph,  
any edge contraction breaks the proper edge coloring 
(see an example, a rank $D=3$
edge contraction in a HEcTG in Figure \ref{fig:contens}).  
However, in rank $D$, there is a color structure on $\crrG_{\mf^0}/e$ defined in a weaker sense. 
Such a weak coloring, that we plan to investigate, is  based on the property that any
stranded edge contraction preserves faces,  the coloring of vertex points, and the bi-coloring of faces. 
 Finally, several HEcTGs may generate the same w-colored graph by edge contractions. 
  
 \begin{lemma}[Stability of the boundary graph under stranded edged contraction]
\label{lem:fullcont}
 \noindent$\bullet$ Contracting a stranded edge in a rank $D$ HESG  $\crrG_{\mf^0}$ 
does not change its boundary graph. 

 \noindent$\bullet$
The contraction in arbitrary order  of all stranded edges of $\crrG_{\mf^0}$ 
yields a HESG $\crrG^0_{\mf^0}$ determined by the boundary graph $\brrG$ of $\crrG_{\mf^0}$ up to additional trivial circles. 
\end{lemma} 
\proof Consider a  HESG  $\crrG_{\mf^0}$ with boundary graph $\brrG$. 
After contracting a stranded edge in $\crrG_{\mf^0}$,  all sHEs remain untouched and all open faces might shrink
by loosing strands and chords but the resulting set of open faces is  in one-to-one correspondence with the former. 
No open face can be created by such operations (we remove all present inner faces
and all remaining faces get shorter but are never cut). Hence, for the resulting stranded graph 
$\crrG_{\mf^0}/e$, one has $\partial (\crrG_{\mf^0}/e)= \brrG$. 
By iteration, contracting all edges of $\crrG_{\mf^0}$ in arbitrary order, the resulting graph $\crrG_{\mf^0}^0$ should contain $\brrG$ as its boundary. 
 Contracting all edges of $\crrG_{\mf^0}$, one obtains a single stranded vertex graph per connected component plus trivial circles.
Then, necessarily, $\crrG_{\mf^0}^0$ is nothing but a graph
made only with stranded vertices (without edges) and sHEs attached to these. Then 
$\crrG^0_{\mf^0}$ is therefore defined by $\partial \crrG^0_{\mf^0}= \brrG$ up to trivial circles. 

\qed

First, we  observe that the previous lemma restricts to half-edged tensor graphs without colors.
Second,  consider $\crrG_{\mf^0}$ a connected HEcTG
with a nonempty set of sHEs, $\mf^0\neq \emptyset$.
After a full contraction of all edges of $\crrG_{\mf^0}$, Lemma 
\ref{lem:fullcont} tells us that the end result $\crrG_{\mf^0}^0$ is totally encoded in the boundary graph $\brrG$ up to some circles.
If $\mf^0 = \emptyset$, then  there is no boundary in $\crrG_{\emptyset}$, 
and $\crrG^0_{\emptyset}$ consists of trivial circles, if non empty. 

 Our next goal is to investigate the properties that
the  contraction and cut operations have on w-colored graphs.

\begin{lemma}\label{lem:internal}
 Let $e$ be an  ordinary stranded edge or a bridge  of  a rank  $3$ w-colored 
graph. Contracting $e$ yields a stranded vertex which is connected. 
\end{lemma}

\proof 

An argument on  the parity of the number of pre-edge points of a given color pair and the fact that two strands of  a colored 
stranded edge in a w-colored graph 
cannot have the same color pair achieve this proof.  The detail of the proof 
follows. 

 Without loss of generality, let us assume that $e$ is of a given color, say 0. Each of its strands is of color pair $(0i)$, $i=1,2,3$. 
Let us concentrate on a single stranded 
vertex $v$, of vertex graph $v'$,
where $e$ is incident at a pre-edge $f_e$. 
If the contraction 
of $e$ disconnects the  resulting stranded pre-vertex, there are at least two
non-empty and distinct sets of pre-edges in $v$, namely $v_1$ and $v_2$, that form two disctinct subgraphs $s_1$
and $s_2$ of $v'$, respectively, that 
are connected only via the vertex corresponding to $f_e$ in $v'$. In other words, there are no edge connections between $s_1$ 
and $s_2$. Note also that $f_e$ does not belong  neither to $v_1$ nor to $v_2$. 
 
There are three strands in $e$. The three of them cannot join the same subset of pre-edges, otherwise it would mean that $v$ was initially disconnected (if $v_1$ and $v_2$ are non empty) or 
$v$ was made only with two pre-edges which also trivializes the proof. Thus, 
at most 2 strands of $e$ could join one subset of pre-edges and the last should connect the other one. 

Consider the strand of color pair $(01)$ in $e$. 
This strand is incident to a chord (at a pre-edge point) incident itself to another pre-edge called $f$. The pre-edge $f$ is necessarily of color 0 or 1. Without loss of generality, let us assume that $f$ is the unique
pre-edge connected to $f_e$ and that also belongs to the set of pre-edges $v_1$. 
Because we can exchange the role played by the colors 0 and 1, 
we can fix the color of $f$ to be 1. Denote $N_i$ the number of pre-edges of color $i$ in $v_1$. 

-  There are $N_i$ pre-edge points of color pair $(i,j)$ in $v_1$, for fixed $i=0,2,3$, and $j\in \{0,1,2,3\}$ but $j\ne i$.   

- There are $N_1$ pre-edge points of color pair $(1,j)$  in $v_1$, for fixed $j=2,3$. 

- There are $N_1-1$ pre-edge points of color pair $(1,0)$ in $v_1$. 

The chords in $v_1$ should connect all these pre-edges points. Denoting $k^{ij}$ the number of chords of color $(ij)$ in $v_1$, we have
\bea
&&
N_0 + N_1 - 1 = 2k^{01} \crcr
&&
N_0 + N_2  =2k^{02} \crcr
&&
N_1 + N_2 = 2k^{12}
\eea
which entails an inconsistent equation: 
$2N_0 + 2(k^{02} - k^{12}) -1 = 2k^{01}$.

\qed

The above lemma should  admit an extension to any rank $D$ by partitioning $D$ and analyzing 
the connection of the strands of $e$ and the chords of the vertices that it meets. We do not need however such a stronger
 result for the following. 

In contrast,  the contraction of a  loop may disconnect the vertex. 
 In any case, since loop graphs are terminal forms, a special treatment is required  for them.

The whole above construction  maintains the consistency of the definition of $p$-bubbles. In a w-colored graph, stranded vertices are still 0-bubbles
(connected objects with 0 color), 
stranded edges are 1-bubbles (connected objects with 1 color), (closed and open) faces are bi-colored maximal connected objects,  (closed and open)
3-bubbles are maximal connected objects with 3 colors. 
It should be however noticed that :

(1) closed faces of the trivial circle vertices naturally inherit of the bi-coloring of
the faces they are coming from; 

\begin{figure}[h]
 \centering
     \begin{minipage}[t]{.8\textwidth}
      \centering
\includegraphics[angle=0, width=8cm, height=2cm]{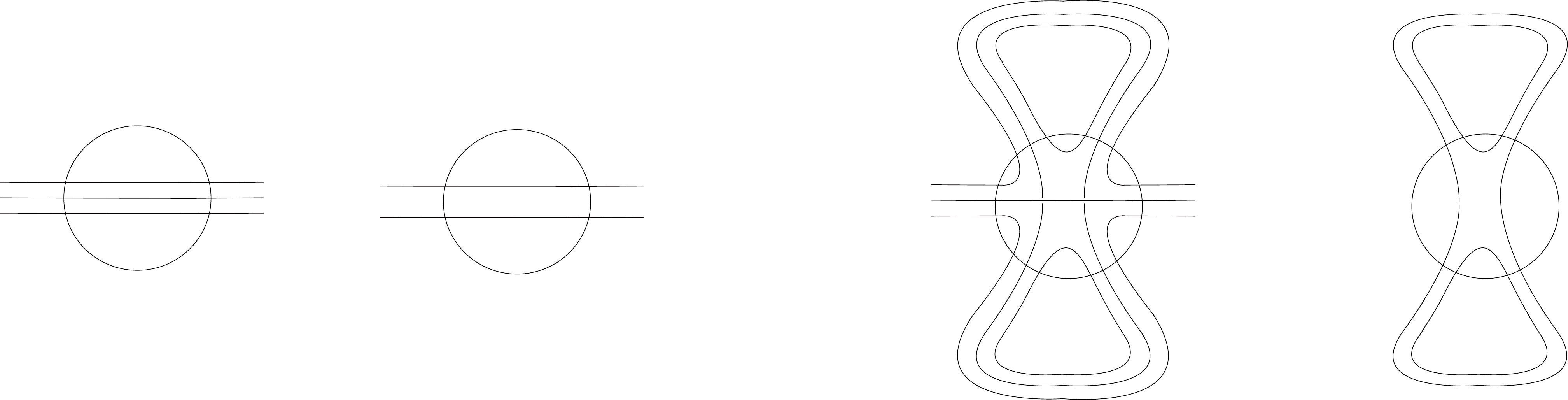}
\vspace{0.1cm}
\caption{ {\small Two bubbles graphs $\bee_1$ and $\bee_2$ of 
$\crrG_{1;\mf^0}$ and $\crrG_{2;\mf^0}$, respectively: $\bee_1$ has a vertex with valence 2 
and $\bee_2$ a vertex with valence 4.}}
\label{fig:npb}
\end{minipage}
\put(-275,0){$\crrG_{1;\mf^0}$}
\put(-213,0){$\bee_1$}
\put(-240,25){$0$}
\put(-295,25){$0$}
\put(-138,-10){$\crrG_{2;\mf^0}$}
\put(-73,-10){$\bee_2$}
\put(-158,25){$0$}
\put(-108,25){$0$}
\end{figure}

(2) 3-bubbles are no longer built uniquely with three valent vertices.  
They can be made with vertices with lower or greater valence as illustrated in Figure \ref{fig:npb}. 
 In any case, Definition \ref{def:pbub} is still valid and we shall focus on this.

 From now on, general discussions  are always valid at fixed rank $D$,  thus we shall omit to mention it. 
 Focusing on  operations on w-colored graphs,  key notions pertain again to the cut/contraction rules.  
\begin{lemma} \label{lem:cutmg}
Let $\crrtG$  be a rank $D$ w-colored graph and $e$ one of its stranded edges. 
$\crrtG \vee e$ the cut graph along $e$,  is a rank $D$ w-colored graph.   $\crrtG/e$, 
called the contraction of  $\crrtG$ along $e$,  is a rank $D$ w-colored graph.  
  \end{lemma} 
\proof  
 Let $\crrtG$ be the  result of a contraction of some HEcTG $\cGcol$. We want to show that both 
 $\crrtG \vee e$ and $\crrtG/e$ come from edge contractions of some HEcTGs. 
The HEcTG $\cGcol \vee e$ obviously contracts to give $\crrtG \vee e$.
There exists a HEcTG $\cGcol^0$ contracting
to $\crrtG/e$. This is nothing but $\cGcol$ on which, before performing all 
contractions yielding $\crrtG$, one performs the contraction
along $e$. 

\qed

To define spanning c-subgraphs of a  w-colored graph, we follow the same procedure as the one dealing with HESGs.  

\medskip 

\noindent{\bf  Rank 3 w-colored graph.}
We henceforth restrict the rank of stranded graphs to $D=3$, and  aim at studying an invariant polynomial
satisfying a contraction/cut relation, and   generalizing the 
Tutte and BR polynomials for these new graphs. The extension for any $D$ should require more work on the subsequent analysis.  
From now rank 3 w-colored graphs  are called more simply w-colored graphs. 

Mostly, we have studied non loop edges so far and their
main properties under contraction and cut can be guessed in any rank. Loops are more subtle
and their behavior under contraction is much more involved.
Restricting to rank 3 simplifies the analysis. 

 The contraction of a loop edge may disconnect a stranded vertex and the resulting HESG. Several cases may be discussed according to the number of connected components into which the vertex graph decomposes upon contraction and whether there are stranded edges linking these components.  

The contraction of a loop edge $e$ may remove some chords within the stranded vertex $v$ on which $e$ is incident. From this removal may result a disconnection the vertex graph of $v$. 
Let us study this in details. 

Given a strand $s$ of a stranded edge $e$
incident to stranded vertex $v$. 
We call \emph{sector} of $s$, the union of
the  set of chords
and set pre-edge points in $v$ that defines a connected 
component vertex graph on which is incident $s$,
when $e$ is contracted. 
If a strand contributes to an inner loop, 
then its sector is empty. 
Stranded edges and sHEs may be incident to some pre-edges of a given sector.
A sector and its incident stranded edges and sHEs will 
be pictured as a black diagram. 
Looking at the set of strands of a fixed stranded edge $e$ incident to $v$, the sectors associated with the strands are called sectors of $v$.  
For a rank $D=3$ HESG, for a given rank $D=3$ stranded edge $e$
incident to $v$, we have at most three sectors in $v$.

\begin{figure}[h]
 \centering
     \begin{minipage}[t]{.8\textwidth}
      \centering
\includegraphics[angle=0, width=10cm, height=2.5cm]{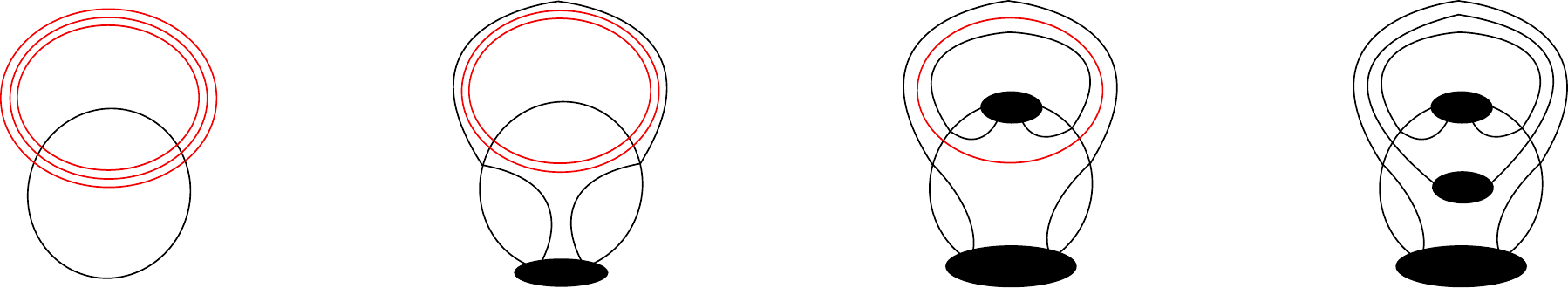}
\vspace{0.4cm}
\caption{ {\small  $p$--inner loops: a 3--inner  (O)
and a 2--inner (A) loop with their unique possible configuration; a 1--inner loop with two potential sectors $v_1$ and $v_2$ (B),
and a 0--inner loop with three  potential sectors $v_1,v_2$ and $v_3$ (C). }}
\label{fig:mpinner}
\end{minipage}
\put(-296,-19){O}
\put(-215,-19){A}
\put(-133,-19){B}
\put(-52,-19){C}
\put(-133,51){$v_1$}
\put(-133,-7){$v_2$}
\put(-52,51){$v_1$}
\put(-52,30){$v_2$}
\put(-52,-7){$v_3$}
\end{figure}
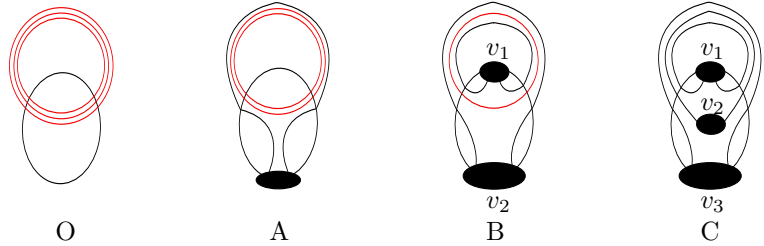

We now restrict to the case
of rank 3 w-colored graphs. 
A loop can be at most a 3--inner edge.
A 3--inner edge determines by itself a graph, see Figure \ref{fig:mpinner}O.   
If it is a 2--inner loop, then $e$ should be of the  
form illustrated in Figure \ref{fig:mpinner}A. Otherwise, $e$ may be a 
1--inner loop (Figure \ref{fig:mpinner}B), or a 0--inner (Figure \ref{fig:mpinner}C).

For a 2--inner loop $e$,  
the outer strand of $e$ should be incident to chords linked to pre-edges
meeting sHEs or edges (otherwise $e$ would be a 3--inner edge)  : it has a unique sector 
and we represent the rest of the stranded vertex possibly with edges   and sHEs attached, by a black diagram in Figure \ref{fig:mpinner}A.

A 1--inner loop $e$ has two outer strands  which should be incident to, potentially, two distinct sectors, $v_1$ and $v_2$ of the vertex,  
as illustrated in Figure \ref{fig:mpinner}B. Any sector, $v_1$ or $v_2$, has at least 
two  sHEs or edges because two chords with the same color
pair are incident on their pre-edges.

A  0--inner loop  $e$ has  potentially three such sectors, $v_1,v_2$ and $v_3$.  
Each sector should contain at least two sHEs or edges. 

The notion of trivial loop in HEcTG can be addressed at this point.   We call a loop $e$ \emph{trivial} if it is a 3--inner or a 2--inner loop, or if it is a 1--inner with exactly 2 sectors 
(upon contraction of $e$ leads to 2 connected component vertex graphs) or 0--inner loop with
exactly 3 sectors (upon contraction of $e$ leads to 3 connected component vertex graphs)
such that there are no stranded edges linking different sectors.

\begin{figure}[h]
 \centering
     \begin{minipage}[t]{.8\textwidth}
      \centering
\includegraphics[angle=0, width=9cm, height=1.5cm]{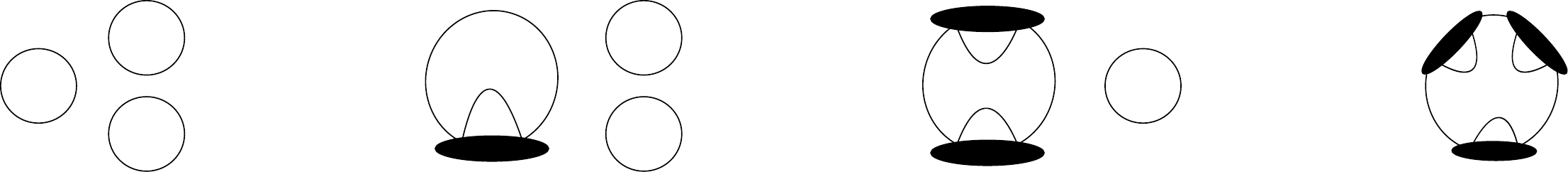}
\vspace{0.3cm}
\caption{ {\small $p$--inner loop contractions
corresponding to O, A, B and C of Figure \ref{fig:mpinner}, respectively. }}
\label{fig:wpcont}
\end{minipage}
\put(-283,-17){O}
\put(-211,-17){A}
\put(-132,-17){B}
\put(-60,-17){C}
\put(-144,45){$v_1$}
\put(-144,-5){$v_2$}
\put(-75,36){$v_1$}
\put(-45,36){$v_2$}
\put(-62,-5){$v_3$}
\end{figure}
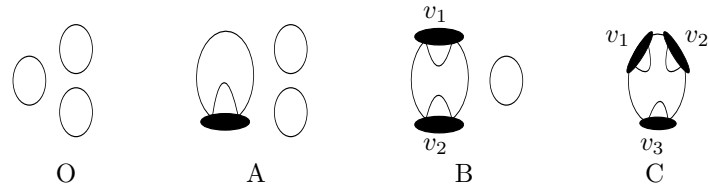

For a 3--inner loop, the contraction gives three trivial circles, see Figure \ref{fig:wpcont}O. 
 The contraction of a trivial 2--inner loop  is again straightforward and yields Figure \ref{fig:wpcont}A. 
For a 1--inner  loop contraction (see Figure \ref{fig:wpcont}B), one gets one additional trivial circle,
 and has two possible configurations: either the stranded vertex 
remains connected or it gets disconnected with two (possibly non trivial) stranded vertices in both situations.  
Concerning the contraction of a trivial 1--inner loop, it is immediate that the vertex gets disconnected in two 
non trivial vertices.  
For a 0--inner loop contraction (see Figure \ref{fig:wpcont}C), we have no additional 
circles but three types of configurations with up to three 
disconnected (and possibly non trivial) stranded vertices.
The contraction of a trivial 0--inner loop yields three disconnected and non trivial stranded vertices. 

\begin{lemma}[Rank 3 trivial loop contraction]
\label{lem:tris}
Let  $\crrtG$ be a w-colored graph  with boundary graph $\brrG$,  $e$ be a trivial loop of $\crrtG$, 
$\crrtG'=\crrtG/e$ and its  boundary graph denoted by $\brrG'$. 
 Let $k$ denote the number of connected components of $\crrtG$, $V$ its number of stranded vertices, $E$ its number of 
 stranded edges, $F_{\inter}$ its number of closed faces, $B_\inter$ its number of closed bubbles, and $B_\ext$ its number of open bubbles, $C_\partial$ the number of connected components of $\brrG$, $f=V_\partial$ its number of sHEs, $E_\partial=F_{\ext}$ the number of edges,  and $F_\partial$
the number of (closed) faces of  $\brrG$, and let $k',V',E',F_{\inter},C'_\partial, B'_{\inter},B'_\ext, C'_\partial,f', E'_\partial$, and $F'_\partial$ denote the analog numbers for $\crrtG'$ and  $\brrG'$. 

If $e$ is a 3--inner loop, then 
\bea
&&
k' = k + 2\,, \quad 
V' = V+2  \,, \quad 
E' = E-1 \, , \quad 
F'_{\inter} = F_{\inter} \,,  \crcr
&&
C'_\partial= C_\partial \,, \quad
f' = f\,, \quad
E'_\partial = E_\partial\, ,
\quad  F'_{\partial } = F_{\partial }\,, \crcr
&& B'_\inter = B_{\inter} -3\,, \qquad 
B'_\ext =  B_\ext\,.
\label{eq:cself3}
 \eea

 If $e$ is a 2--inner loop, then 
\bea
&&
k' = k + 2 \,, \quad 
V' = V+2  \,, \quad 
E' = E-1 \, , \quad 
F'_{\inter} = F_{\inter} \,, \crcr
&&
C'_\partial= C_\partial \,, \quad
f' = f\,, \quad 
E'_\partial = E_\partial \,, \quad
 F'_\partial = F_\partial \,, \crcr
&&
B'_\inter = B_{\inter} -1\,, \qquad 
B'_\ext =  B_\ext\,.
\label{eq:cself2}
 \eea

 If $e$ is a trivial $p$--inner loop such that $p=0,1$,
then 
\bea
&&
k' = k + 2 \,, \quad 
V' = V+2  \,, \quad 
E' = E-1 \, , \quad 
F'_{\inter} = F_{\inter} \,, \crcr
&& 
C'_\partial= C_\partial \,, \quad
f' = f\,, \quad 
E'_\partial = E_\partial \,, \quad
 F'_\partial = F_\partial \,, \crcr
&& 
B'_\inter + B'_\ext = B_{\inter} + B_\ext  + \alpha_p \,,
\label{eq:cselfp}
 \eea
where $\alpha_{p} =3-2p$.

\end{lemma}

\proof 
 The first situation of a 3--inner loop  defines a particular graph without sHEs made of a single vertex with a unique loop. 
The contraction removes the vertex and  the three closed bubbles, and generate three trivial circles. 
The result \eqref{eq:cself3} follows.

For $p$--inner  loops, $p\leq 2$, there should be other  sHEs or edges on the same end vertex. 
The first lines of \eqref{eq:cself2}-\eqref{eq:cselfp} should not cause any difficulty using the definition 
of the contraction which preserve (open and closed) strands. Let us focus on the
variations of the number of bubbles of the graph.

 Consider that $e$ is 2--inner. One first should notice that there is a closed bubble formed by the two inner faces of $e$ and that contracting $e$ destroys this  bubble. In addition, the other (open or closed) bubbles,  which are made of strands of $e$ and containing one of the two inner faces, are simply deformed. They loose one face (the inner) but remain in the contracted graph. The second line of \eqref{eq:cself2} follows. 

 Let us assume now that $e$ is trivial and 1-- or 0--inner, and  of color $a$. 
Consider in the initial  vertex $v$, the two or three distinct sectors  $v_1$, $v_2$ and $v_3$, respectively, determined by the strands of $e$.  
We recall that $v_1$, $v_2$, and $v_3$ should contain each at least two  sHEs or edges. Consider as well the three 
bubbles labeled by colors $(aa_1a_2)$, $(aa_1a_3)$ and $(aa_2a_3)$  containing strands of  $e$. 

 We treat the case of $e$ trivial 1--inner loop, and let us assume that the face $(aa_3)$ 
is the one which is inner, without loss of generality.  The edge $e$ decomposes in three
ribbon edges each determined by a couple of pairs $(aa_i; aa_j)$, 
$i<j$, such that the strand $(aa_1)$ connects to the sector $v_1$ and the strand $(aa_2)$ connects to $v_2$. 
  If $e$ is trivial and 1--inner, its contraction yields two connected vertices plus one circle. 
Consider the bubbles $\bee_{aa_1a_3}$
and $\bee_{aa_2a_3}$ (in the whole graph) which are the only ones containing the face $(aa_3)$. 
 Contracting $e$, the bubbles $\bee_{aa_1a_3}$ and $\bee_{aa_2a_3}$ just loose that face and get deformed. 
They are still present since, in $v_1$ and $v_2$, there are still 
faces that are included in $\bee_{aa_1a_3}$ and $\bee_{aa_2a_3}$. Meanwhile, the last bubble $\bee_{aa_1a_2}$ 
 splits in two parts: these parts are disjoint bubbles with the same triple of colors since there are no
edges in the stranded graph that relate them.  These bubbles can be open or
closed depending on the nature of $\bee_{aa_1a_3}$ in $v_1$ and $v_2$. Hence, $B'_{\inter} + B'_\ext = B_{\inter} + B_\ext+1$.

Finally, let us discuss the trivial 0--inner loop situation  with its three distinct sectors $v_1$, $v_2$ and $v_3$. The discussion
is somehow similar to the 1--inner case.  The contraction of $e$ yields three connected vertices.
Using the similar routine explained above, one finds that
each of the bubbles $\bee_{aa_ia_j}$, $i<j$, gets split into 
two parts after the contraction. Each of these pairs of parts 
are not related at all by any stranded edge. Thus each bubble 
produces two bubbles. The resulting bubbles  may be open or closed depending on the nature of $\bee_{aa_ia_j}$
in each sector. We have $B'_{\inter} + B'_\ext = B_{\inter} + B_\ext+3$.

\qed

 It is clear that the contraction of a 2--inner loop 
also obeys \eqref{eq:cselfp} as well for $p=2$. In fact \eqref{eq:cself2} is a stronger result because it mainly distinguishes the variation of the numbers of closed and open bubbles.  For a general 1--inner loop contraction, 
it turns out that the number $B_\inter + B_\ext$ may vary 
of 1 or not. Similarly, for a non necessarily trivial 0--inner loop contraction, the same number of bubbles  vary from 0, 1, 2, up to 3.

\begin{lemma}[Faces of a bridge]\label{lem:str}
 The faces passing through a bridge of any rank 3 w-colored graph are necessarily open. 
 These faces belong to the same connected component of the boundary graph of the w-colored graph.  
\end{lemma} 

\proof A bridge $e$ with color $a$ has three strands of color pairs $(aa_1)$, $(aa_2)$  and $(aa_3)$ defining three different faces passing through $e$. Clearly, the strand coloring in any stranded edge prevents a face to pass twice through $e$, and so the three faces using strands of $e$ are all open. The next statement follows from the fact that the boundary graph is 3-regular
made with an even number of vertices and has a proper edge coloring.  

\qed

\begin{lemma}[Cut/contraction of special edges]\label{lem:cutbri}
Let $\crrtG$ be a  rank 3 w-colored graph
and $e$ a stranded edge in $\crrtG$. Then, in the above notation,

 if $e$ is a bridge,   we have 
\bea
&&
k(\crrtG\vee e) = k(\crrtG/e)+1 ,\;
V(\crrtG\vee e) = V(\crrtG/e) + 1 ,\crcr
&& 
E(\crrtG\vee e) = E(\crrtG/e),\;
f(\crrtG \vee e)= f(\crrtG/e) +2,   \label{eq:kve} \\
&&
\crcr
&&
F_{\inter}(\crrtG\vee e)=  F_{\inter}(\crrtG/e) \,, \quad
B_{\inter}(\crrtG \vee e)=B_{\inter}(\crrtG/e)  \,,
\label{eq:fint} \\
&&  
C_\partial(\crrtG \vee e) = C_\partial(\crrtG/e)  +1\,,
\quad 
E_\partial(\crrtG \vee e) = E_\partial (\crrtG/e) +3 
\,, \crcr
&& 
F_\partial(\crrtG \vee e) = F_\partial(\crrtG/e) + 3 \,,
\label{eq:cextbord} \\
&&
 B_{\ext}(\crrtG\vee e) = B_{\ext}(\crrtG/e) + 3\,;
\label{eq:bubext} 
\eea

 if $e$ is a trivial $p$--inner loop, $p=0,1,2$, we have  
\bea
&&
k(\crrtG\vee e) = k(\crrtG/e)-2 ,\;
V(\crrtG\vee e) = V(\crrtG/e)- 2 ,\crcr
&& 
E(\crrtG\vee e) = E(\crrtG/e),\;
f(\crrtG \vee e)= f(\crrtG/e) +2, 
\label{eq:kvep} \\
\cr 
&&
F_{\inter}(\crrtG\vee e)+ C_\partial(\crrtG \vee e) 
=  F_{\inter}(\crrtG/e) + C_\partial(\crrtG/e)-2\,, 
\crcr
&& 
 E_\partial(\crrtG \vee e) = E_\partial (\crrtG/e) + 3\,,
\label{eq:fintp} \\
&&
B_{\inter}(\crrtG \vee e) +  B_{\ext}(\crrtG\vee e) 
=B_{\inter}(\crrtG/e)  + B_{\ext}(\crrtG/e) - (3-2p)\,.
\label{eq:cextp} 
\eea
 Moreover, given $\brrG$ the boundary of $\crrtG$ and
 a bridge or a trivial $p$--inner loop $e$,  $p=0,1,2,3$,
\beq
2C_\partial(\crrtG \vee e) - \chi(\partial(\crrtG \vee e)) =
2C_\partial(\crrtG) -  \chi(\crrtG) = 2C_\partial(\crrtG/e) -  \chi(\partial(\crrtG/e)) \,,
\label{eq:chi} 
\eeq
where $\chi(\brrG)$ denote the Euler characteristics of the boundary of $\crrtG$.

 \end{lemma}
\proof We start by the bridge case. The equations in \eqref{eq:kve} are easily found. 
Let us focus on \eqref{eq:fint}. By Lemma \ref{lem:str}, 
we know that, necessarily, the faces passing through $e$ are open. 
All closed faces on each side of the bridge are conserved after cutting $e$. The same are still conserved after
edge contraction. Hence $F_{\inter}(\crrtG\vee e) = F_{\inter}(\crrtG/e)$
and $B_{\inter}(\crrtG \vee e)=B_{\inter}(\crrtG/e)$. 
We now prove \eqref{eq:cextbord}. By the second point of 
Lemma \ref{lem:str}, the three open faces belong to the same 
boundary component. After cutting $e$, this unique component
yields two  boundary components.  It is direct to get 
$C_\partial(\crrtG \vee e) = C_\partial(\crrtG/e)  +1$,
$E_\partial(\crrtG \vee e)= E_\partial (\crrtG/e)+3$
(the cut of $e$ divides each open face into two different open faces) and 
$F_\partial(\crrtG \vee e)= F_\partial (\crrtG/e)+3$
because $C_\partial(\crrtG/e)=C_\partial(\crrtG)$,
$E_\partial (\crrtG/e)=E_\partial (\crrtG)$
and $F_\partial(\crrtG /e)= F_\partial (\crrtG)$
 which are immediate from Lemma \ref{lem:fullcont}.
Concerning the number of open bubbles, there are three bubbles in $\crrtG$,  each with an edge made of strands of the bridge. Each of these is associated with 
two color pairs $(aa_i;aa_j)$, $i<j$.  These bubbles
are clearly in $\crrtG/e$. Cutting the bridge, each of these bubbles splits in two.  This yields \eqref{eq:bubext}. 

We now deal with a trivial $p$--inner loop $e$. 
The relations \eqref{eq:kvep} can be determined without difficulty and so 
we concentrate on the rest of the equations. 
 Consider the faces $f_i$ made with outer strands of $e$. 
For $p=0,1,2$,  we have $f_{i}$, $1\leq i\leq 3-p$. These 
faces can be open or closed. We do a case by case study
according to the number of open or closed faces among the
$f_i$'s. 

- Assume that $3-p$ of $f_i$'s are closed.
Cutting $e$ entails 
\bea
&&
F_{\inter}(\crrtG\vee e) = F_{\inter} (\crrtG) -3
\,,\quad 
C_{\partial}(\crrtG \vee e) = C_\partial (\crrtG) +1
\, , \quad
 E_\partial(\crrtG \vee e) = E_\partial(\crrtG) +3 \,,\crcr
&&
B_{\inter}(\crrtG \vee e)  + B_{\ext}(\crrtG \vee e) 
 = B_{\inter}(\crrtG)  + B_{\ext}(\crrtG) \,,
\label{eq:f3} \\
&& 
 F_\partial(\crrtG \vee e) = F_\partial(\crrtG) + 3 \,.
\label{eq:fpa3} 
\eea
Note that in this situation only the variation of the total number of bubbles can be known.

- Assume that $3-p-1$ of the $f_i$'s are closed
and one is open. Cutting $e$ gives 
\bea
&&
F_{\inter}(\crrtG\vee e) = F_{\inter} (\crrtG) -2
\,,\quad 
C_{\partial}(\crrtG \vee e) = C_\partial (\crrtG) 
\, , \quad 
 E_\partial(\crrtG \vee e) = E_\partial(\crrtG) +3\,,\crcr
&&
B_{\inter}(\crrtG \vee e)  = B_{\inter}(\crrtG)  -1\,,
\quad 
 B_{\ext}(\crrtG \vee e)  = B_{\ext}(\crrtG)  + 1\,,
\label{eq:f31}
\\
&&
 F_\partial(\crrtG \vee e) = F_\partial(\crrtG) +1 \,.
\label{eq:fpa31}
\eea

- Assume that $3-p-2$ of the $f_i$'s are closed
and two are open. Cutting $e$ gives 
\bea
&&
F_{\inter}(\crrtG\vee e) = F_{\inter} (\crrtG) -1
\,,\quad 
C_{\partial}(\crrtG \vee e) = C_\partial (\crrtG) -1
\, ,\quad
 E_\partial(\crrtG \vee e) = E_\partial(\crrtG) +3\,,\crcr
&&
B_{\inter}(\crrtG \vee e)  = B_{\inter}(\crrtG)  \,,
\quad 
 B_{\ext}(\crrtG \vee e)  = B_{\ext}(\crrtG)  \,,\label{eq:f32}\\
&&
 F_\partial(\crrtG \vee e) = F_\partial(\crrtG)- 1\,.
\label{eq:fpa32}
\eea
Note that this case does not apply for $p=2$.

- For $p=0$ an additional situation applies:  assume that all 
 three $f_i$'s are open. Cutting $e$ gives 
\bea
&&
F_{\inter}(\crrtG\vee e) = F_{\inter} (\crrtG) 
\,,\quad 
C_{\partial}(\crrtG \vee e) = C_\partial (\crrtG) -2
\, ,\quad
 E_\partial(\crrtG \vee e) = E_\partial(\crrtG) +3\,,\crcr
&&
B_{\inter}(\crrtG \vee e)  = B_{\inter}(\crrtG)  \,,
\quad 
 B_{\ext}(\crrtG \vee e)  = B_{\ext}(\crrtG)  \,,\label{eq:f33}\\
&&
 F_\partial(\crrtG \vee e) = F_\partial(\crrtG)- 3\,.
\label{eq:fpa33}
\eea
Lemma \ref{lem:tris} relates  the same numbers for $\crrtG/e$ and $\crrtG$, and from this, we can prove \eqref{eq:fintp} and \eqref{eq:cextp}.

Last, we prove \eqref{eq:chi} on the Euler characteristics 
of the boundary  graphs. We first note that, from \eqref{eq:kve}, \eqref{eq:cextbord}, \eqref{eq:kvep} and  \eqref{eq:fintp}, for any special (bridge or trivial $p$--inner, $p=0,1,2$) edge,  
\beq
f(\crrtG \vee e) = f(\crrtG)+ 2 = f(\crrtG/e) + 2\,, \qquad 
E_\partial (\crrtG \vee e) =E_\partial(\crrtG) +3 = E_\partial(\crrtG/e) +3 \,.
\label{eq:fe}
\eeq
For the bridge case, \eqref{eq:chi} follows from the relations
$F_\partial (\crrtG \vee e) =  F_\partial(\crrtG) +3 = F_\partial(\crrtG/e) +3$ 
and $C_\partial (\crrtG \vee e) = C_\partial(\crrtG) +1= C_\partial(\crrtG/e) +1$ 
in \eqref{eq:cextbord}. The result holds also for trivial $0,1,2$--inner loops, after the case by case study giving \eqref{eq:f3}--\eqref{eq:fpa33}.
Last, for the 3--inner loop, in addition to \eqref{eq:fe} which still holds, the following relations are valid
\bea
C_\partial(\crrtG \vee e) = C_\partial(\crrtG)+1 = C_\partial(\crrtG/e)+1\,,
\quad 
F_\partial (\crrtG \vee e) =F_\partial(\crrtG) +3 = F_\partial(\crrtG/e) +3 \,,
\eea
and allow us to conclude.

\qed

\subsection{Polynomial invariant for $3D$ w-colored graphs}
\label{subset:polymt}

We shall define first an invariant for rank 3 w-colored graphs,
check its consistency and then state our main result. 

\begin{definition}[Topological invariant for rank 3 w-colored 
graph]\label{def:topotens}
Let $\crrtG$ be a rank 3 w-colored graph.
The generalized topological invariant associated with $\crrtG$
is given by the following  function 
\bea
&&
\mT_{\crrtG}(x,y,z,s,w,q,t)= \label{ttopfla}\\
 && \sum_{A \sset \crrtG}
 (x-1)^{r(\crrtG)-r(A)}(y-1)^{n(A)}
z^{5k(A)-[3(V- E(A))+2(F_{\inter}(A)-B_{\inter}(A)-B_{\ext}(A))]} \crcr
&& \qquad \times \quad  s^{C_\partial(A)} \,  w^{F_{\partial}(A)} q^{E_{\partial}(A)} t^{f(A)}\,.
\nonumber
\eea
\end{definition}

 A crucial point is to show that the exponent of $z$ 
in \eqref{ttopfla} is always a non negative integer.

\begin{proposition}
\label{prop:euler}
Let $\crrtG$ be a rank 3 w-colored 
graph  without trivial circles. Then 
\beq
\zeta(\crrtG)= 3(E(\crrtG)- V(\crrtG)) + 2[ B_{\inter}(\crrtG) + B_{\ext}(\crrtG) -F_{\inter}(\crrtG) ] \geq 0\,.
\eeq
\end{proposition}  
\proof 
 Let us consider $\cB_{\inter}$ and $\cB_{\ext}$, the sets of closed and  open bubbles of $\crrtG$,
respectively. (In this proof, we drop the notation $\crrtG$ in all
quantities depending on the w-colored graph.) 
Let $B_{\inter}$ and $B_{\ext}$ be their respective cardinality. 
Each open or closed bubble $\bee$ is an HERG 
 with $V_{\bee}$ number of vertices,
$E_{\bee}$ number of edges, $F_{\inter; \bee}$ number
of closed faces, and $C_{\partial}(\bee)$  number of cycles
of the  boundary of $\bee$. We write $\bee_i$ for a closed bubble
and $\bee_x$ for an open bubble. The Euler characteristics 
of $\bee_i$ and $\bee_x$ refer to the same notion for 
their underlying ribbon graphs. 

Any $\bee_i \in \cB_{\inter}$ being a ribbon graph,  its Euler characteristics writes
\beq
 2- \kappa_{\bee_i} = V_{\bee_i} - E_{\bee_i} + F_{\inter; \bee_i}\,,
\eeq
where $\kappa_{\bee_i}$ refers to the genus of $\bee_i$ or twice its genus 
if $\bee_i$ is oriented. Summing over all closed bubbles, we get
\beq
2 B_{\inter} - \sum_{\bee_i \in \cB_{\inter}} \kappa_{\bee_i}
  = \sum_{\bee_i \in \cB_{\inter}}\left[ V_{\bee_i} - E_{\bee_i} + F_{\inter; \bee_i}\right] .
\label{eq:euli}
\eeq
Using the colors, one observes that 
each edge of $\crrtG$ splits into three ribbon edges belonging either
to an open or a closed bubble, and each closed face of $\crrtG$ belongs to  two bubbles which might be open or closed. Thus we have
\beq
\sum_{\bee_i \in \cB_{\inter}} E_{\bee_i} + \sum_{\bee_x
\in \cB_{\ext}} E_{\bee_x}  = 3E\, , \qquad
\sum_{\bee_i \in \cB_{\inter}} F_{\inter;\bee_i} + \sum_{\bee_x
\in \cB_{\ext}} F_{\inter;\bee_x}  =2F_{\inter}   \,.
\label{eq:efi}
\eeq
In addition, each vertex of the graph can be decomposed,
at least, in three vertices (3 vertices is the minimum given by 
the simplest vertex of the form $\crrG_{1;\mf^0}$ in Figure \ref{fig:npb}) which could belong to an open or a closed bubble giving then
\beq
\sum_{\bee_i \in \cB_{\inter}} V_{\bee_i} + \sum_{\bee_x
\in \cB_{\ext}} V_{\bee_x}  \geq 3V \,.
\label{eq:veri} 
\eeq
Combining \eqref{eq:efi} and \eqref{eq:veri}, we re-write \eqref{eq:euli} as
\beq
3V - 3E +2F_{\inter} - 2B_{\inter} -
\sum_{\bee_x \in \cB_{\ext}}\left[ V_{\bee_x} - E_{\bee_x} + F_{\inter; \bee_x}\right] \leq   - \sum_{\bee_i \in \cB_{\inter}} \kappa_{\bee_i} \,.
\label{eq:pres}
\eeq
We complete the last sum involving $\cB_{\ext}$ by adding
$C_{\partial}(\bee_x)$ in order to get
\bea
\sum_{\bee_x \in \cB_{\ext}}\left[ V_{\bee_x} - E_{\bee_x} + F_{\inter; \bee_x} + C_{\partial}(\bee_x)\right] = 
\sum_{\bee_x \in \cB_{\ext}} (2 - \kappa_{\bee_x})\,,
\eea
which, substituted in \eqref{eq:pres}, leads us to
\bea
3V - 3E +2F_{\inter} - 2B_{\inter} - 2B_{\ext} 
 \leq  - \sum_{\bee_i \in \cB_{\inter}} \kappa_{\bee_i} 
- \sum_{\bee_x \in \cB_{\ext}} \big( C_{\partial}(\bee_x) + \kappa_{\bee_x}\big) 
\label{eq:pres2}
\eea
 that  proves the lemma. 
\qed

\vspi

In fact, $ \sum_{\bee_x \in \cB_{\ext}} C_{\partial}(\bee_x)= F_\partial(\crrtG)$ is the number of faces of the boundary graph of
$\crrtG$ (see Remark \ref{rem:1}). The above bound can be refined since
$
F_\partial(\crrtG) \geq B_{\ext}(\crrtG)
$
which merely follows from the fact that each $\bee_x \in \cB_{\ext}$
has at least a boundary component contributing to $F_\partial(\crrtG)$
such that 
\bea
B_{\ext}(\crrtG) \leq \sum_{\bee_x} C_{\partial}(\bee_x) = F_{\partial}(\crrtG)\,.
\eea
Thus we also have 
\beq
3V - 3E +2F_{\inter} - 2B_{\inter} - B_{\ext} 
 \leq   - \sum_{\bee_i \in \cB_{\inter}} \kappa_{\bee_i} 
- \sum_{\bee_x \in \cB_{\ext}} \kappa_{\bee_x}\,.
\eeq
This may lead as well to yet another well defined invariant.
However, we do not use this relation in the following
due to some rich relations that  $-\zeta(A) \geq 0$
entails. We do have other  positive combinations.

\begin{proposition}\label{prop:14}
Let $\crrtG$ be a rank 3 w-colored 
graph  without trivial circles. Then 
\bea
&&
\zeta'(\crrtG)= 3[E(\crrtG)- V(\crrtG)] + 2[ B_{\inter}(\crrtG) + B_{\ext}(\crrtG) -F_{\inter}(\crrtG) -C_\partial(\crrtG)] \geq 0\,,\crcr
&&
\zeta''(\crrtG)= 3[E(\crrtG)- V(\crrtG)-C_\partial(\crrtG)] + 2[ B_{\inter}(\crrtG) + B_{\ext}(\crrtG) -F_{\inter}(\crrtG)] \geq 0\,.
\eea
\end{proposition} 
\proof By the color prescription, each connected component of the boundary of $\crrtG$ has at least three faces. Therefore
\bea
3 C_\partial \leq F_\partial \,. 
\eea
Then we also have $-F_\partial +2C_\partial \leq 0$. The proposition follows from \eqref{eq:pres2} and the fact that $F_\partial = \sum_{\bee_x} C_{\partial}(\bee_x) $ in the proof of Proposition \ref{prop:euler}.

\qed

\begin{proposition}\label{prop:euler2}
Let $\crrtG$ be a rank 3 w-colored graph with $D(\crrtG) \ge 0$ trivial 
circles.  Then 
\beq\label{tildezeta}
\tilde\zeta(\crrtG)= 3(E(\crrtG)- V(\crrtG)) + 2[ B_{\inter}(\crrtG) + B_{\ext}(\crrtG) -F_{\inter}(\crrtG) ] + 5D(\crrtG) \geq 0\,.
\eeq
\end{proposition} 
\proof 
Consider $\crrtG$ a  rank 3 w-colored graph with $D (\crrtG) \ge 0$ trivial circles.  Some steps in the proof of Proposition \ref{prop:euler} should be modified as follows:   in the relation \eqref{eq:efi},  we replace $F_{\inter}(\crrtG) $ by $F_{\inter}(\crrtG)  - D(\crrtG) $, and in \eqref{eq:veri}, $V(\crrtG) $ by $V (\crrtG) - D(\crrtG) $. 
Then, we obtain the desired inequation obeyed by $\tilde\zeta(\crrtG)$. 

\qed

Of course, there exist some modified $\zeta'(\crrtG)$
and $\zeta''(\crrtG)$ that handle the generic case of a rank 3 w-colored graph with trivial circles by simply adding to them $5D(\crrtG)$. 
A second remark is that, as $k(\crrtG) \ge D(\crrtG)$, then  another choice of a positive quantity is still 
given by \eqref{tildezeta} but trading $D(\crrtG)$
for $k(\crrtG)$,  the total number of connected components
of the w-colored graph. From the combinatorial perspective, although
the choice of $\tilde\zeta(\crrtG) \ge 0$ would have been 
enough to conduct the analysis, 
we will use instead the number of connected components
of the w-colored graph which is a global topological quantity. 

Proposition \ref{prop:euler2} ensures 
us that the following statement holds.

\begin{proposition}[Polynomial invariant]
$\mT_{\crrtG}$ is a polynomial. 
\end{proposition}

The quantity $V - E(A) + F_{\inter}(A) - B_{\inter}(A)$ 
is nothing but the Euler characteristics 
for a  colored tensor graph $A$ without sHEs understood as a cellular complex \cite{Gurau:2009tw}. 
 This quantity is bounded by $-\sum_{\bee_i} \kappa_{\bee_i}$.  
 In the case of a HEcTG, the same cellular
complex has a boundary bearing itself a cellular decomposition. 
 We interpret  $\zeta(A)$, in the present situation, as a weighted notion of a Euler characteristics 
 of the cellular complex corresponding to $A$ which also takes into account its boundary. 
 
 As a result of the excess $k(\crrtG)- D(\crrtG)$, 
$\mT_{\crrtG}$ factors by an overall monomial 
that we can keep track. 
The c-spanning subgraphs  $A\sset \crrtG$ have 
the same number of trivial circles as $\crrtG$.  
Using $k(A) \ge k(\crrtG)$, we write 
 $k(A) - D(\crrtG) + D(\crrtG) = (k(\crrtG)  - D(\crrtG)) 
 + (k(A) - k(\crrtG) ) +  D(\crrtG )$ 
 and learn that the monomial in $z$ 
 factors as  $z^{5(k(\crrtG)  - D(\crrtG))}z^{5(k(A) - k(\crrtG)  +  D(\crrtG ))}$.
 Therefore, $\mT_{\crrtG}$ can be factored out by $z^{5(k(\crrtG)  - D(\crrtG))}$. 
This allows us to cover
the case of an invariant based on  $\tilde\zeta$. 

At an even more general level,  
we should emphasize that
any combination $\xi_\alpha(\crrtG) =  5D(\crrtG) + \alpha(k(\crrtG)- D(\crrtG))$, for any $\alpha \ge 0$, 
in the exponent $z^{\xi_\alpha(\crrtG)  + \tilde\zeta(\crrtG)- 5D(\crrtG) }$ 
will determine a valid polynomial invariant. In this work, we restrict 
to the case $\alpha=1$. 

Let us call $\crrtG^0$ a w-colored graph made only with a finite set of stranded vertices with sHEs and no stranded edges,
and possibly with trivial circles. Then $E(\crrtG^0)= B_{\inter}(\crrtG^0)=0$, $F_{\inter}(\crrtG^0)=D$, 
$D$ being the number of trivial circles in $\crrtG^0$, $k(\crrtG^0)= V(\crrtG^0)$, 
$k(\crrtG^0)-D=C_{\partial}(\crrtG^0)\geq 0$, we can check the consistency of \eqref{ttopfla} as the following makes still sense:
 \beq
\label{eq:bounmt}
\mT_{\crrtG^0}(x,y,z,s,w,q,t)= 
z^{2[k(\crrtG^0)-D+B_{\ext}(\crrtG^0)]} 
s^{k(\crrtG^0) - D}  w^{F_\partial(\crrtG^0)} q^{E_\partial(\crrtG^0)}\,t^{f(\crrtG^0)}. 
\eeq

It may exist several possible reductions of the above 
polynomial. We  will focus on the following: 
\bea
&&
\mT_{\crrtG} (x,y,z,z^{-2},w,q,t) = \mT'_{\crrtG} (x,y,z,w,q,t)  \,,\crcr
&&
 \mT_{\crrtG} (x,y,z,z^{-2}s^{2},s^{-1},s,s^{-1}) = \mT''_{\crrtG} (x,y,z,s)\,, \label{eq:mtprim} \crcr
&&
\mT_{\crrtG} (x,y,z,z^{2}z^{-2},z^{-1},z,z^{-1}) = \mT'''_{\crrtG} (x,y,z) \,. \crcr
&&
\eea
Proposition \ref{prop:14} ensures that $\mT'_{\crrtG}$
is a polynomial. Meanwhile, $\mT''_{\crrtG}$ is also a polynomial
with exponent of $s$ the Euler characteristics of the boundary $\brrG$.
$\mT'''_{\crrtG} (x,y,z) $ combines both invariants in a single exponent. These polynomials will be relevant in the next analysis.  
Note that we could have introduced another polynomial
expressed as  $\mT_{\crrtG} (x,y,z,s^{2},s^{-1},s,s^{-1}) = \mT^0_{\crrtG} (x,y,z,s)$.
But this reduction turns out to satisfy the same properties as 
$\mT$ and thus does not provide anything new.

We are now in position to prove our main theorem. 

\begin{theorem}[Contraction/cut rule for w-colored graphs]
\label{theo:contens}
Let $\crrtG$ be a rank 3 w-colored graph. Then, for an   ordinary   edge $e$ of $\crrtG$, we have
\beq
\mT_{\crrtG}=\mT_{ \crrtG\vee e} +\mT_{\crrtG/e}\,,
\label{tenscondel}
\eeq 
for a bridge $e$, we have $\mT_{\crrtG \vee e}= z^8s(wq)^3t^2 \mT_{\crrtG/e}$
and 
\beq
\mT_{\crrtG} =[(x-1)z^8s(wq)^3t^2+1] \mT_{\crrtG/e}\,;
\label{tensbri}
\eeq 
 for a trivial $p$--inner loop $e$, $p=0,1,2$, we have 
\beq
\mT_{\crrtG}=\mT_{\crrtG \vee e} + (y-1)z^{4p-7}\,\mT_{\crrtG/e}\,.
\label{tensself}
\eeq 
\end{theorem}

\proof
Let  $\crrtG$ be a rank 3 w-colored graph. 
 The same  preliminary remarks of the proof of Theorem \ref{theo:BRext} hold also here for $\crrtG$, with adapted consideration, e.g. sHEs replace HRs.  
 Our main concern is the change in the number of closed 
and open bubbles. 

We concentrate first on \eqref{tenscondel}. 
 Consider an ordinary edge $e$ of $\crrtG$  of color $a$, the set of spanning c-subgraphs which do not contain $e$ being the same  
as the set of spanning c-subgraphs of $\crrtG \vee e$, 
 the number of open and closed bubbles on each subgraph is the same, 
it is direct to get 
 \bea\label{sumnotin}
&&
 \sum_{A \sset \crrtG ; e\notin A}
 (x-1)^{r(\crrtG)-r(A)}(y-1)^{n(A)}
z^{5k(A)-[3(V- E(A))+2(F_{\inter}(A)-B_{\inter}(A)-B_{\ext}(A))]} 
\crcr
&&  \qquad \qquad  \times
  s^{C_\partial(A)} \,  w^{F_{\partial}(A)} q^{E_{\partial}(A)} t^{f(A)}=\mT_{\crrtG \vee e} .
\eea 
Let us focus on the second term of \eqref{tenscondel}.
Consider $e$, its end vertices $v_1$ and $v_2$ and
 its 3 strands with color pairs $(aa_1)$, $(aa_2)$ 
and $(aa_3)$ and consider the set of bubbles in $\crrtG$. 
Some bubbles do not use any strand of $e$  and three bubbles can be formed using these
strands (these bubbles are of colors $(aa_ia_j)$, $i<j$). Contracting  $e$, 
the vertex obtained is connected by Lemma \ref{lem:internal}.
 The cycles and paths made of strands of $e$ are clearly preserved after the contraction.  
The bubbles  which use no strands of $e$ are not affected at all 
by the procedure. The  three bubbles using 2 strands of $e$
are also preserved since  the contraction does not delete faces. The faces passing through $e$ getting 
only shortened,  the  result from the point of view of the bubbles 
passing through $e$ is simply an ordinary ribbon edge contraction in the sense of HERGs. 
Using this, we write
\bea\label{sumin}
&&
\sum_{A \sset \crrtG ; e\in A}(x-1)^{r(\crrtG)-r(A)}(y-1)^{n(A)}
z^{5k(A)-[3(V- E(A))+2(F_{\inter}(A)-B_{\inter}(A)-B_{\ext}(A))]} 
\crcr
&&  \qquad \qquad  \times
  s^{C_\partial(A)} \,  w^{F_{\partial}(A)} q^{E_{\partial}(A)} t^{f(A)}
=\mT_{\crrtG /e} .
\eea 
Let us focus now on the bridge case and \eqref{tensbri}. 
Cutting a bridge yields, as in the ordinary case, 
from the sum $\sum_{A \sset \crrtG ; e\notin A}$ the product
$(x-1)\mT_{\crrtG \vee e}$. The second sum remains as it is using the mapping between 
$\{A \sset \crrtG;e\in A\}$ and $\{A\sset \crrtG/e\}$
and provided the fact that the resulting vertex is still connected. That is, once again, ensured  by Lemma \ref{lem:internal}.  
The last stage relates $\mT_{\crrtG \vee e}$ and $\mT_{\crrtG/e}$. 
This can be achieved by using the bijection between $A\sset \crrtG \vee e$ and 
$A' \sset \crrtG/e$ where each $A$ and $A'$ are both 
uniquely related to some $A_0\sset \crrtG$ as $A = A_0\vee e$
and $A'=A_0/e$. Using Lemma \ref{lem:cutbri}, the relation \eqref{tensbri} follows. 

Next, we discuss the trivial $p$--inner loop case and prove \eqref{tensself}. 
The property \eqref{sumnotin}   should be direct. The second term in \eqref{tensself} is now studied. 

The question is whether or not $e$
being a trivial $p$--inner loop in $\crrtG$ remains as such in $A$.
The answer for that question is yes because $A$ contains the end vertex 
of $e$. We may cut some edges in each sectors $v_i$ (see Figure \ref{fig:mpinner})  for defining $A$ 
but the resulting sectors are distinct in $A$. 

Contracting a trivial $p$--inner loop generates $p$ circles and $3-p$ non trivial vertices. 
Now, from Lemma \ref{lem:tris}, we know how  all numbers of components in the graph evolve: the  nullity is again 
$n(A) = n(A') + 1$, and the exponent of $z$ becomes
\bea
&&
5k(A) - (3(V(\crrtG) - E(A)) + 2(F_\inter(A) - B_\inter(A)-B_{\ext}(A)) = \crcr
&&
5(k(A')-2) - \Big[3[V(\crrtG/e)-2) - (E(A')+1)]
  + 2(F_\inter(A') - B_\inter(A')-B_{\ext}(A')+\alpha_p)\Big]\crcr
&&
 = 5k(A') - (3(V(\crrtG/e) - E(A')) + 2(F_\inter(A') - B_\inter(A')-B_{\ext}(A')) 
 -7 +4p \,,
\label{eq:binter}
\eea
where $\alpha_p=3-2p$, $A$ and $A'$ refer to the subgraphs 
related by bijection the usual between spanning c-subgraphs of $\{A \sset \crrtG;e\in A\}$ and $\{A'\sset \crrtG/e\}$.  
At the end, one gets $(y-1)z^{4p-7}\mT_{\crrtG/e}$. 

 \qed 

We realize again that the relations of Theorem \ref{theo:contens} are not complete reduction rules.  General boundary conditions as 
generalized bouquets must be explicitly computed using
Definition \ref{def:topotens}.

\begin{corollary}[Cut/contraction relations for $\mT'$]
\label{coro:trpin}
Let $\crrtG$ be a rank 3 w-colored graph. Then,  for a bridge $e$,
we have $\mT'_{\crrtG \vee e}= z^{6}(wq)^3t^2 \mT'_{\crrtG/e}$ and
\beq
\mT'_{\crrtG}(x,y,z,w,q,t) =[(x-1)z^6(qw)^3t^2+1] \mT'_{\crrtG/e}(x,y,z,w,q,t)\,;
\label{eq:ctprimb}
\eeq 
for a trivial $p$--inner loop $e$ in $\crrtG$, $0\leq p \leq 2$, 
we have
\beq
\mT'_{\crrtG}(x,y,z,1,q,t)=z^{4p-6} [q^3t^2+ (y-1)z^{-1}]\,\mT'_{\crrtG/e}(x,y,z,1,q,t)\,;
\label{tensselfp}
\eeq 

\end{corollary}

\proof Theorem \ref{theo:contens} implies naturally \eqref{eq:ctprimb}
 after the setting $s=z^{-2}$ in $\mT_{\crrtG}$.  
We now work out the contraction/cut of the trivial loops. 

Same arguments as in the proof of Theorem \ref{theo:contens}
should be invoked here. The difference now is that, 
by changing $s\to z^{-2}$, using 

(1) the one-to-one mapping 
between $\{A \sset \crrtG\vee e\}$ and $\{A' \sset \crrtG/e \}$
such that with each $A \sset \crrtG\vee e$ one associates 
$A' \sset \crrtG/e$ defined by $A'=\tilde A/e$,
$\tilde A\sset \crrtG$, and $\tilde A \vee e = A$ and

(2) the relations 
\eqref{eq:kvep}-\eqref{eq:cextp} of  Lemma \ref{lem:cutbri},

\noindent we can map, for a trivial $0,1,2$--inner loop $e$,
$\mT_{\crrtG \vee e}$ on $\mT_{\crrtG/e}$. We compute the variation of the exponent of $z$ between $A$ and $A'$ as
\bea
&&
5k(A) - (3(V(\crrtG) - E(A)) + 2[F_\inter(A) +C_\partial(A)- B_\inter(A)- B_{\ext}(A)]) =\crcr
 &&
5(k(A')- 2 )- (3(V(\crrtG')-2- 
E(A')) + 2[F_\inter(A') +C_\partial(A')-2- 
B_\inter(A') - B_\ext(A') + \alpha_p ])
\crcr
&& 
 = 5k(A') - (3(V(\crrtG/e) - E(A')) + 2[F_\inter(A') - B_\inter(A') - B_\ext(A')])  -2\alpha_p \,.
\label{eq:bin}
\eea
The rest of the variations can be easily identified using the same
lemma. 

\qed

 \begin{corollary}[Cut/contraction rules for $\mT''$ (and $\mT'''$)]
\label{coro:trpin2}
Let $\crrtG$ be a rank 3 w-colored  graph. 

Then, for a bridge $e$, we have $\mT''_{\crrtG \vee e} = z^{6}\mT''_{\crrtG/e}$
and 
\beq
\mT''_{\crrtG} =[ (x-1)z^6 + 1]\,\mT''_{\crrtG/e}\,;
\label{eq:cttrb}
\eeq
for a trivial $p$--inner loop $e$ in $\crrtG$, $0\leq p \leq 2$, 
we have $\mT''_{\crrtG \vee e} =z^{4p-6} \mT''_{\crrtG/e}$ and 
\beq
\mT''_{\crrtG}= z^{4p-6}[1+ (y-1)z^{-1}]\,\mT''_{\crrtG/e}\,.
\label{eq:terself}
\eeq 
The same contraction and cut rules applies for $\mT'''_{\crrtG}$. 
\end{corollary}
\proof The new ingredient to achieve the proof of this statement 
is \eqref{eq:chi} of Lemma \ref{lem:cutbri}. 

\qed

The exponents of  $z^{4p-7}$ or of $z^{4p-6}$ can be negative in some cases.  This simply implies that,
in the polynomial $\mT'_{\crrtG /e}$ or $\mT''_{\crrtG /e}$, all monomials should contain an enough large power of $z$
to make the overall exponent of $z$ positive. 

The disjoint union operation on graph extends naturally in the present formulation.

\begin{lemma}[Disjoint union]\label{lem:disj}

Let $\crrtG$ and $\crrtG'$ two disjoint rank 3 $w$-colored graphs, 
then 
\beq
\mT_{\crrtG \sqcup \crrtG'} = 
\mT_{\crrtG } \, \mT_{\crrtG'} \,.
\eeq
\end{lemma}
\proof This is totally standard as in the ordinary procedure
using additive properties of  exponents in $\mT_{\crrtG}$. 

\qed

\begin{corollary}[3--inner loop contraction]
\label{coro:term}
Given a rank 3 w-colored graph $\crrtG$ containing a 3--inner loop $e$
then 
 \beq
\mT_{\crrtG} = z^5(z^3s(wq)^3t^2 + y-1) \mT_{\crrtG/e}\,.
\eeq
 \end{corollary}
\proof We use the fact that $e$ is a 3--inner and so it forms
a separate subgraph $g$. In order to compute the polynomial
of $\crrtG$, Lemma \ref{lem:disj} can be applied and a direct evaluation of $g$  yields the desired factor.  

\qed

\begin{definition}[Multivariate form]
The multivariate form associated with \eqref{ttopfla} is defined by:
 \bea
&&
\widetilde\mT_{\crrtG}(x,\{\beta_e\},\{z_i\}_{i=1,2,3},s,w,q,t) \label{tmulti}=\\
&&    \sum_{A \sset \crrtG}
 x^{r(A)}(\prod_{e\in A}\beta_e)
z_1^{F_{\inter}(A)}z_2^{B_{\inter}(A)} z_3^{B_{\ext}(A)}\, 
s^{C_\partial(A)} \,  w^{F_{\partial}(A)} q^{E_{\partial}(A)} t^{f(A)}\,,
\nonumber
\eea
for $\{\beta_e\}_{e\in \cE}$ labeling the edges of the graph $\crrtG$. 
\end{definition}
It is direct to prove the following statement by ordinary techniques. 

\begin{proposition} For any   ordinary   edge $e$, 
\bea
 \widetilde\mT_{\crrtG}=
\widetilde\mT_{ \crrtG\vee e} +x\beta_e\,\widetilde\mT_{\crrtG/e}\,.
\eea
\end{proposition}

\begin{remark}\label{rem:2}
We can  compare $\widetilde\mT$ with the Gurau  polynomial
denoted in the following by $G$ \cite{Gurau:2009tz}. Note that we will not use  the full form of $G_{\crrtG}$, denoted $P_{\crrtG}$  in \cite{Gurau:2009tz}, but will  instead introduce two improvements: 
(1) a normalization form  $P_{\crrtG}(\{\beta_e x_1\},\dots)$ of $P_{\crrtG}$,
where $x_1$ is the variable associated with the number of edges
which brings an inessential overall factor of $x_1^{E(A)}$
consistently absorbed by the $\beta_e$, 
and (2) a rank formulation of $G_{\crrtG}(\{\beta_e\},\dots) =P_{\crrtG}(\{\beta_e x_1\},\dots)$, i.e. rather then using two variables $x_0$ for the vertices
and $x_4$ for the number of connected components of the rank 3 colored tensor graph, we  simply use $x^{r(A)}$, 
$A \sset \crrtG$.

For a rank 3 colored tensor graph $\crrtG$, the polynomials $\wt\mT$ and $G$ are related by 
\beq
\widetilde\mT_{\crrtG}(x,\{\beta_e\}, z_1,z_2,z_3=1,s,w,q,t)
 = G_{\crrtG}(x,\{\beta_e\},z_1,z_2,s,q,w,t)\,, 
\eeq
with, according to the convention in  \cite{Gurau:2009tz}, we have
\bea
C_\partial = |\mathfrak{B}_{\partial}^3|, \quad 
F_\partial = |\mathfrak{B}_{\partial}^2|, \quad 
E_\partial = |\mathfrak{B}_{\partial}^1|, \quad 
f = |\mathfrak{B}_{\partial}^0|, \quad 
B_{\inter}= |\mathfrak{B}^3|, \quad 
F_{\inter} = |\mathfrak{B}^2|.
\eea
As expected, for this set of graphs the polynomial $\wt\mT$ is more general because it contains one additional variable ($z_3$ fixed to 1)
  associated with the number of open bubbles. This variable can be introduced
by hand in $G$ making it a bit more general. This additional variable does not lead to any new features for the multivariate form
however, as seen in our developments, it plays an important role in the non-multivariate form $\mT$.  

\end{remark}

\section{Conclusion}
\label{concl}

We have generalized the Tutte and BR polynomials to rank 3 w-colored graphs. 
The rank $D$  w-colored graphs are new combinatorial objects obtained from the (stranded) edge contraction 
of rank $D$ half-edged, colored and stranded graphs. Colored stranded and tensor graphs have appeared in several contexts
in physics, in particular, as Feynman graphs of field theories describing quantum geometry and gravity. 
The new polynomial invariant $\mT$ presented in this work satisfies a contraction/cut recursion relation on rank 3 w-colored graphs, with 
the cut operation which seems natural in this context.  
We have  evaluated some boundary cases or terminal forms of the contraction/cut recursion relation for these graphs. 
The multivariate form of the invariant has been given, and its relation with the Gurau polynomial, see \cite{Gurau:2009tz},  has been established. 

Studying the limit cases, the invariant $\mT$ 
reduces to the Tutte polynomial on its underlying graph. 
The connection with the BR polynomial
on  half-edged ribbon graphs might be also achieved after reduction of the present analysis to rank 2. 

The perspectives of this work are numerous.  
 First, we must mention 
 that it is not excluded that  rank $D$ w-colored
graphs and all other graph species studied in this work 
find other definitions. 
For instance, there is another interesting underlying graph
of an HESG  $\crrG_{\mf^0}$. The {\emph string graph} $S(\crrG_{\mf^0})$ of $\crrG_{\mf^0}$ is the graph whose vertex set is the set of vertex points and SHe external points of $\crrG_{\mf^0}$, and whose edge set is the
set of strands and chords in $\crrG_{\mf^0}$, the incidence
relation between the vertices and edges of $S(\crrG_{\mf^0})$ being obvious. 
 Appearing as a watermark in Definition \ref{innerface}, such a graph that might reveal a new point of view and perhaps brings simplification of the present results. 
 This is certainly worth investigating. 

For instance, the universality property of $\mT$ must be addressed, once the notion of  the one-vertex
w-colored graphs is well mastered.  The introduction of colors might help in classifying  the bouquets of w-colored graphs  
 needed in this proof. In fact, before addressing the universality issue on stranded structures, one first needs to understand to how
the universality property can be extended for the extension of the BR polynomial to HERGs. 
This will amount to generalizing the notion of chord diagrams 
$\mathcal{D}$, their equivalence relation under rotations about chords, and finally 
their associated basic canonical diagrams $\mathcal {D}_{ijk}$ (counting 
the nullity $i$, genus $j$ of the diagram, $k$ is associated with
the orientation of the surface associated with diagram) \cite{bollo,bollo3}
to \emph{chord diagrams with half-edges}. Such a preliminary task 
will certainly imply the existence of a new canonical diagram 
$\mathcal {D}_{ijkl}$, where $i$, $j$, $k$ keep their meaning, and 
$l$ defines the number of connected components associated with the boundary of the HERG.  
It is also known that the Tutte and  BR polynomials can be expressed  in terms of spanning tree expansions and
 enumerate specific trees.  We may ask: what type of tree expansion does $\mT$ satisfy, 
and which objects does it enumerate? These new lines of investigation will 
certain yield interesting results relating combinatorics and topology 
in higher dimension \cite{Brighton}.

Finally, the study of tensor graphs have been motivated
by their importance in quantum gravity. A natural follow up question might be what are the new information brought by the new
class of w-colored graphs or the new polynomial invariant $\mT$ in that context? An answer of this question is not obvious.  
However, in the search of new tools for classifying the 
Feynman graphs which are viewed as random manifolds, the following
investigation tracks might be relevant to address. First, each color in the field theoretic setting represents a different field, and the path integral 
in a field theory integrates over field configurations. Integrating over a color
makes the resulting Feynman graph independent of it. 
Practically, we use an edge contraction to represent the new Feynman
graph. This procedure has been used by Bonzom et al. \cite{Bonzom:2012hw} to reveal a new class of objects called uncolored graphs. Uncolored Feynman graphs are colored tensor graphs on which we perform edge contractions, for all edge colors  except for   one. 
It could be therefore timely to investigate in which sense 
w-colored graphs are interpolating configurations between a fully colored
tensor theory by Gurau and an uncolored one in the sense of Bonzom et al., 
and  making perhaps the formalism developed therein a bit
more general. Furthermore, to find a practical use of the new invariant
$\zeta(\crrtG)$ found in the present work,  one must first interpret,
at the level of the dual simplicial complex, the meaning of 
the contraction, and the cut procedure performed at the
level of the w-colored graph. Thus, 
another possible investigation is to understand dualities at the level
of graphs and simplicial complexes (with boundaries), and consequently dualities satisfied by the polynomial $\mT$ itself.   This  hopefully could lead to  a better understanding of simplicial complexes generated
by path integrals.

\begin{center}
{\bf Acknowledgements}
\end{center}

{\footnotesize Discussions with Vincent Rivasseau at early stage
of this work are gratefully acknowledged. RCA thanks the Perimeter Institute for its hospitality and support in the time that this work had been initiated. 
This work was partially supported by Perimeter Institute for Theoretical Physics. Research at Perimeter Institute is supported by the Government of Canada through Industry Canada and by the Province of Ontario through the Ministry of Research and Innovation.
This work was also partially supported by TWAS Research Grant RGA No. 17-542 RG/MATHS/AF/AC\_G -FR3240300147.The ICMPA-UNESCO Chair is in partnership with Daniel Iagolnitzer Foundation (DIF), France, and the Association pour la Promotion Scientifique de l'Afrique (APSA), supporting the development of mathematical physics in Africa.
}

\section*{ Appendix: Examples}
\label{app}

\appendix

\renewcommand{\theequation}{\Alph{section} \arabic{equation}}
\setcounter{equation}{0}

We carry out explicitly two examples in order to illustrate our results in the present appendix. 

\vspi

\begin{figure}[h]
 \centering
     \begin{minipage}[t]{.8\textwidth}
      \centering
\includegraphics[angle=0, width=6cm, height=3cm]{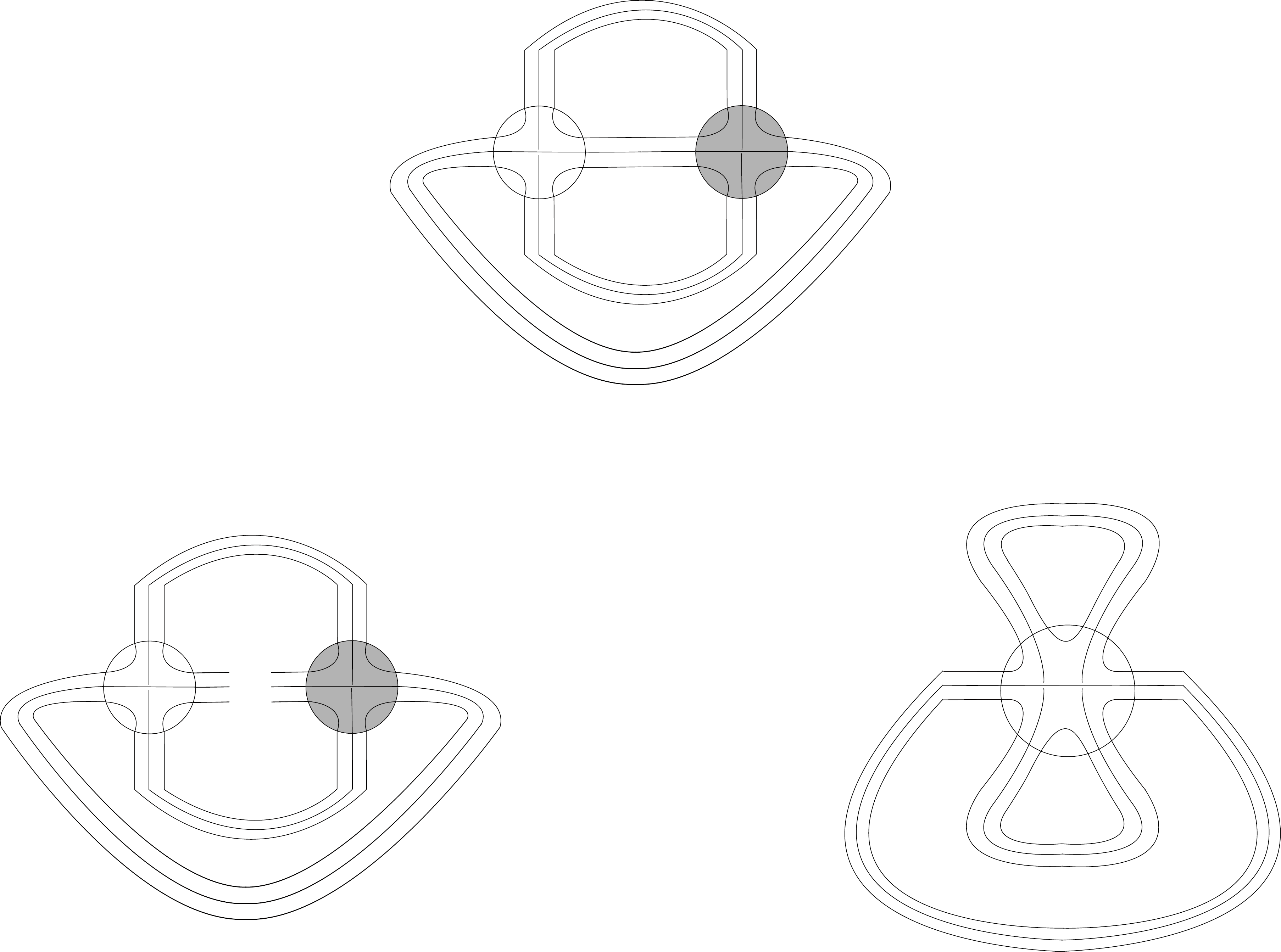}
\vspace{0.1cm}
\caption{ {\small A rank 3 colored tensor graph $\crrG$ and
its cut $\crrG \vee e$ and contraction $\crrG/e$ for an    ordinary edge $e$
of color 2.}}
\label{fig:cg1}
\end{minipage}
\put(-235,-10){$\crrG \vee e$}
\put(-120,-10){$\crrG/e$}
\put(-175,75){$e$}
\put(-170,60){$3$}
\put(-175,90){$0$}
\put(-200,75){$1$}
\put(-175,40){$\crrG$}
\end{figure}
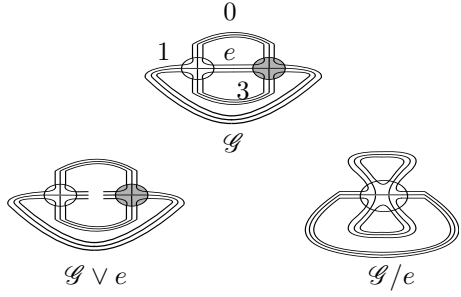

\noindent{\bf Example 1: A colored graph.}
Consider the graph $\crrG$ given by Figure \ref{fig:cg1}. Computing the multivariate form of the polynomial 
by the spanning c-subgraph summation in \eqref{tmulti} with $\beta_i$ is associated with the edge of color $i$,
we obtain
\bea
\mT_{\crrG}(x,\{\beta_e\},z,s,q,t)&=& \beta_0\beta_1\beta_2\beta_3\, xz_1^6z_2^4 \cr &+& 
( \beta_0\beta_1\beta_2 +  \beta_0\beta_2\beta_3 +  \beta_1\beta_2\beta_3 +  \beta_0\beta_1\beta_3)\, x z_1^3 z_2z_3^3sw^{3}q^3t^2 \cr 
&+& ( \beta_0\beta_1 +  \beta_0\beta_2 +  \beta_0\beta_3 +  \beta_1\beta_2 + \beta_1\beta_3 +  \beta_2\beta_3)\, xz_1z_3^{4} sw^{4}q^6t^4 \cr &+& ( \beta_0 + \beta_1 + \beta_2 + \beta_3)\,xz_3^{5}s w^{5} q^9t^6 + z_3^{8}s^2w^8q^{12}t^8\,.
\eea
Then this polynomial coincides exactly with the normalized Gurau polynomial $G_{\crrG}$ after setting to 1 the variable $z_3$ associated with open bubbles. Note that reversely, introducing a new variable associated with the same component in   $G_{\crrG}$,
we infer that both polynomials match for the present example. Note also that the  exponents of $z_3$  and $w$ always coincide. It is, however, not always true that 
each open bubble would have a unique boundary component. 

Cutting one edge, say the one of color 0, yields $\crrG \vee e$, so that  we
 evaluate
\bea
\mT_{ \crrG\vee e}(x,\{\beta_e\},z,s,q,t)&=&  \beta_1\beta_2\beta_3\, xz_1^3 z_2 z_3^{3} sw^3q^{3} t^2 \cr &+&(   \beta_1\beta_2 + \beta_1\beta_3 +  \beta_2\beta_3)xz_1z_3^{4} sw^4q^{6}t^4  \cr &+& ( \beta_1 + \beta_2 + \beta_3)xz_3^{5}sw^5q^{9}t^6 + z_3^{8}s^2w^{8} q^{12}t^8
\eea
and contracting the same edge gives
 \bea
\mT_{\crrG/e}(x,\{\beta_e\},z,s,q,t)&=&\beta_1\beta_2\beta_3\, z_1^6 z_2^4 \cr &+& ( \beta_1\beta_2 +  \beta_2\beta_3 + \beta_1\beta_3)z_1^3 z_2z_3^{3} sw^{3} q^3t^2 \cr &+& (\beta_1 + \beta_2 +  \beta_3)z_1z_3^{4}sw^4q^{6}t^4 + z_3^{5}sw^5q^{9}t^6\,.
\eea
 Thus, we get
\bea
	 \mT_{\crrG}=\mT_{ \crrG\vee e} +x\beta_0\mT_{\crrG/e}\,.
\eea 
 It should be also emphasized that the contraction of an edge is performed in the sense of Definition \ref{def:constran}. This definition 
 allows us to improve the active/passive contraction scheme as used in \cite{Gurau:2009tz}. 

\vspi

\begin{figure}[h]
 \centering
     \begin{minipage}[t]{.8\textwidth}
      \centering
\includegraphics[angle=0, width=12cm, height=6cm]{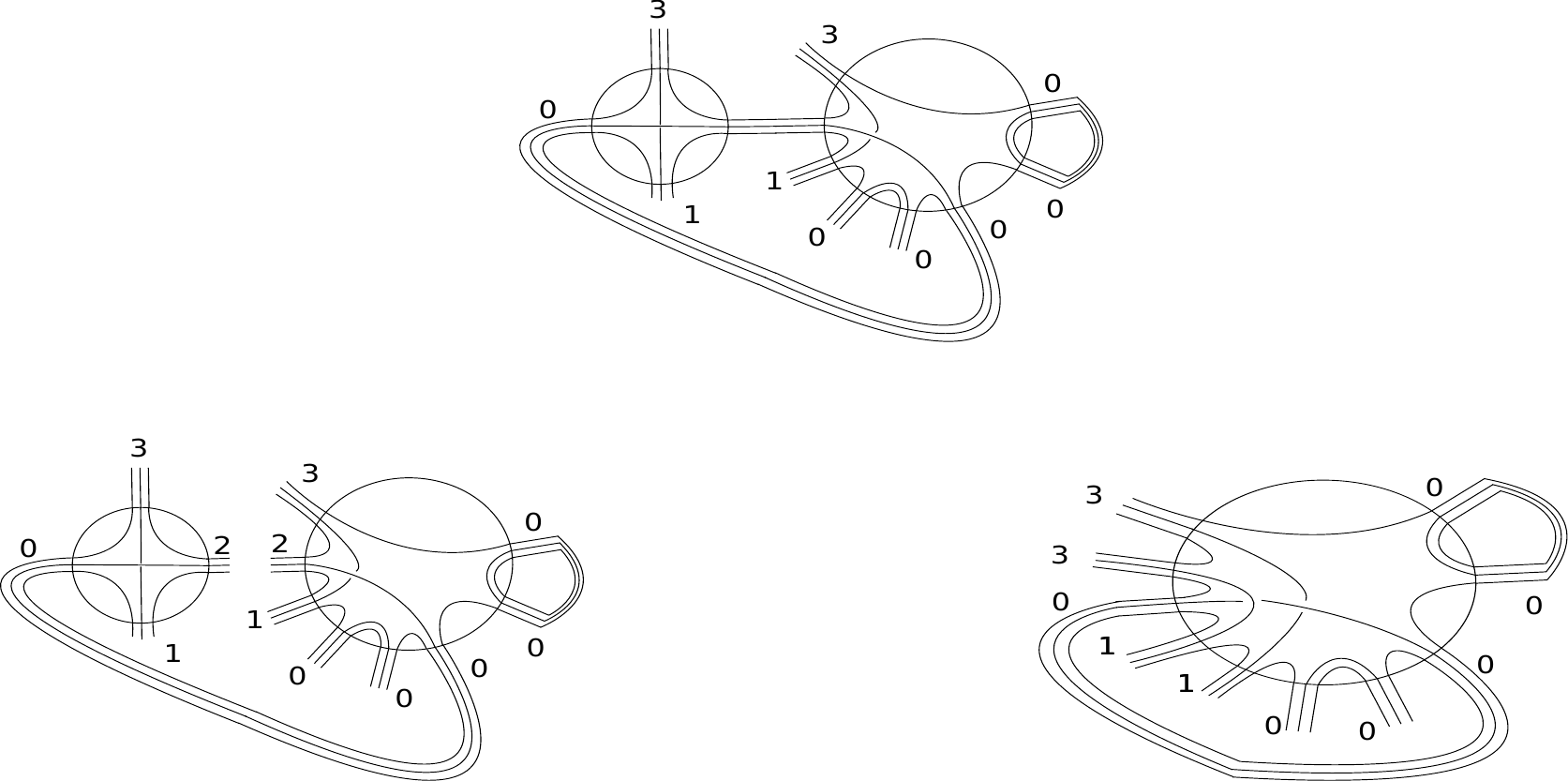}
\vspace{0.1cm}
\caption{ {\small A w-colored graph $\crrtG$, the cut graph $\crrtG\vee e_2$ and the contracted graph $\crrtG/e_2$
with respect to $e_2$.}}
\label{fig:cg2}
\end{minipage}
\put(-177,145){\small{$e_2$}}
\put(-96,136){\small{$e_1$}}
\put(-215,114){\small{$e_0$}}
\put(-325,15){\small{$e_0$}}
\put(-120,15){\small{$e_0$}}
\put(-285,-10){$\crrtG \vee e_2$}
\put(-70,-10){$\crrtG/e_2$}
\put(-210,42){\small{$e_1$}}
\put(5,52){\small{$e_1$}}
\put(-175,85){$\crrtG$}
\end{figure}

\noindent{\bf Example 2: A ``planar'' w-colored graph.}
We can compute $\mT$ for other types of graphs
which are not colored tensor graphs.  
In a specific instance, consider the graph $\crrtG$ of Figure \ref{fig:cg2}. 
It combines both one colored vertex and another type
vertex. For simplicity, we change the variables $(x-1) \to x$
and $(y-1) \to y$.

By the spanning c-subgraph summation, we get 
\bea
\mT_{\crrtG}(x,y,z,s,q,t) = [ x z^{10} s w^5q^{9} t^{6}
+  xyz^{9} s w^{4}q^{6} t^{4}
 +   2z^2w^2q^{6} t^{4} 
+ 3yzw q^{3} t^{2} + y^2 ] z^{14}sw^5q^9t^{6} \,.
\label{eq:mtinit}
\eea

We want to compare \eqref{eq:mtinit} with the result obtained by contraction and cut procedure 
using a much as possible results on terminal forms. Using the notation of Figure \ref{fig:cg2}, we must check
that
\beq
\mT_{\crrtG} = \mT_{\crrtG\vee e_2} +\mT_{\crrtG/e_2} \,.
\label{eq:interm0}
\eeq
Evaluating $ \mT_{\crrtG\vee e_2}$ and $\mT_{\crrtG/e_2}$, one finds 
\bea
&&
\mT_{\crrtG\vee e_2}= [x z^8s(wq)^3t^2 + 1 ]\mT_{(\crrtG\vee e_2)/e_0} 
\,,\quad 
\mT_{(\crrtG\vee e_2)/e_0} = \mT_{((\crrtG\vee e_2)/e_0)\vee e_1}
+ y z \, \mT_{((\crrtG\vee e_2)/e_0)/e_1}\,,\crcr
&& \label{eq:interm1}\\
&&
\mT_{\crrtG/e_2} = \mT_{(\crrtG/e_2) \vee e_1} + yz\,\mT_{(\crrtG/e_2)/e_1}\,.
\label{eq:interm2}
\eea
We used in \eqref{eq:interm1} the fact that $e_0$ is a bridge in $\crrtG \vee e_2$
and that $e_1$ is a 2--inner loop in $(\crrtG \vee e_2)/e_0$. 
Meanwhile, in \eqref{eq:interm2}, we used the fact that 
$e_1$ is a 2--inner loop in $\crrtG/e_2$. 

A straightforward calculation yields
\bea
&&
 \mT_{((\crrtG\vee e_2)/e_0)\vee e_1}= z^{16}sw^7q^{15} t^{10}\,,
\qquad \mT_{((\crrtG\vee e_2)/e_0)/e_1}= z^{14}sw^{6}q^{12} t^{8}\,,\crcr
&&
 \mT_{(\crrtG/e_2) \vee e_1}  = z^{16}sw^{7}q^{15} t^{10} + yz^{15}sw^6q^{12} t^{8}\,, 
\qquad 
\mT_{(\crrtG/e_2)/e_1} = z^{14}sw^{6}q^{12} t^{8} + yz^{13}sw^{5}q^{9} t^{6}\,.\crcr
&&
\eea
Plugging these results on \eqref{eq:interm1} and \eqref{eq:interm2},
and summing their contributions in \eqref{eq:interm0}, allow us to recover \eqref{eq:mtinit}.

\vspace{0.3cm}

\begin{center}
\rule{3cm}{0.01cm}
\end{center}
\end{document}